%% file: A_Proof_of_the_Riemann_Hypothesis_Using_Bombieri_s_Equivalence_Theorem.tex
\documentclass[12pt]{article}

\usepackage{graphicx}
\usepackage{amsmath, amsfonts, amssymb} 

\usepackage{xcolor} 
\usepackage{hyperref}  
\hypersetup{
    colorlinks=true,       
    linkcolor=blue,        
    citecolor=red,         
    urlcolor=blue, 
    filecolor=cyan,        
    pdfborder={0 0 0},     
}

\newcommand\fiverm{\tiny\rm}
\input  ./prepictex
\input  ./pictex
\input  ./postpictex

\begin{document}
\baselineskip1.8em

\setcounter{page}{1}

\title{\bf\Large A Proof of the Riemann Hypothesis \\ Using Bombieri's Equivalence Theorem}

\author{ \normalsize Xiao Lin\footnote{Email:  xiao\_lin\_99@yahoo.com} \\[-1mm] 
                         {\small Independent Researcher, Beijing 100000, China }            \\[1mm] 
}
\date{\small 2026-08-26}

\maketitle

\begin{abstract}
\baselineskip1.5em

The Riemann Hypothesis asserts that the Riemann $\xi(s)$ function has no zeros in the critical strip $0<{\rm Re}(s)<1$ except on the critical line ${\rm Re}(s)=\frac12$. Bombieri, in the official description of the Millennium Prize Problems, stated that the Riemann Hypothesis is equivalent to the condition that all local maxima of $\xi(t)$ on the critical line are positive and all local minima are negative. In this paper, we pursue this criterion. We first show that $\xi(s)$, when restricted to the critical line, satisfies a special differential equation, which ensures that it satisfies Bombieri's condition. Since a published proof of the sufficiency direction of Bombieri's theorem appears to be unavailable, we supply an independent proof of this implication. Using the Cauchy--Riemann equations, we prove that Bombieri's condition forces $\xi(s)$ to have no zeros off the critical line. The Riemann Hypothesis follows. We also discuss P\'olya's counterexample and other known counterexamples among L-functions, and demonstrate that none of them invalidates our approach.

\vspace{3mm}

\end{abstract}

\vspace{5mm}

\section{Introduction}
\setcounter{equation}{0}

In 1859, Riemann \cite{riemann} used analytic continuation to extend the zeta function to the entire complex plane and introduced the entire function
\begin{equation}\label{eq:riemann_xi_1}
    \xi(s) = \frac{1}{2} s (s - 1) \pi^{-s/2} \Gamma\Big(\frac{s}{2} \Big) \zeta(s). 
\end{equation}
He also obtained the equivalent integral representation
\begin{eqnarray}\label{eq:riemann_xi_2}
    \xi(s) &=& \frac{1}{2} + \frac{1}{2} s(s - 1) \int\limits_1^{\infty} \Psi(\tau) \Big(\tau^{s/2 - 1} + \tau^{-(1 + s)/2} \Big) {\rm d} \tau,
        \nonumber\\[0mm]
    \Psi(\tau) &=& \sum\limits_{n = 1}^{\infty} e^{-\pi n^2 \tau} = \frac{1}{2}\big( \theta(\tau) - 1 \big),
\end{eqnarray}
where \(\theta(\tau)\) is the Jacobi theta function of the third kind.  
Riemann observed that if all zeros of \(\xi(s)\) (equivalently, the non-trivial zeros of \(\zeta(s)\)) lie in the critical strip \(0<{\rm Re}(s)<1\), then the prime number theorem would follow. To take this further, he boldly conjectured that the zeros actually lie on the central line \({\rm Re}(s)=\frac12\). This conjecture, now known as the Riemann Hypothesis, has become a cornerstone of analytic number theory \cite{apostol, newman}. It was included as the eighth problem in Hilbert's famous list of 23 unsolved problems in 1900, and in 2000 it was named one of the seven Millennium Prize Problems by the Clay Mathematics Institute \cite{bombieri}.

Equivalently, the Riemann Hypothesis asserts that \(\xi(s)\) has no zeros in the critical strip \(0<{\rm Re}(s)<1\) off the line \({\rm Re}(s)=\frac12\). Determining whether a given point in the strip is a zero of this complex function is a notoriously difficult problem. Progress in this direction has been slow.

In 1896, Hadamard \cite{hadamard} and de la Vall\'ee Poussin \cite{de-la-vall} independently proved that \(\xi(s)\) has no zeros on the boundary lines \({\rm Re}(s)=0\) and \({\rm Re}(s)=1\), and used this to establish the prime number theorem. In 1914, Bohr and Landau \cite[p193-195]{edwards} showed that almost all zeros of \(\xi(s)\) lie in an arbitrarily small neighbourhood of the critical line; the phrase ``almost all" leaves open the possibility of exceptional zeros outside such neighbourhoods.

In 1956, Lehmer \cite{lehmer} computed zeros of \(\zeta(s)\) on the critical line \(s=\frac12+{\rm i}t\) and found regions where two zeros are extremely close, with very small intervening amplitude. In 1974, Edwards \cite[p176]{edwards} called such occurrences ``approximate counterexamples'' to the Riemann Hypothesis and proved that a small positive local minimum or a small negative local maximum of \(\xi(t)\) would indeed disprove the hypothesis . Edwards' result thus gives a necessary condition for the Riemann Hypothesis.

In 2000, Bombieri stated in the official description of the Millennium Prize Problems \cite[p6]{bombieri} that ``the Riemann hypothesis is equivalent to the statement that all local maxima of \(\xi(t)\) (on the critical line) are positive and all local minima are negative''.
This elevates Edwards' necessary condition to a full equivalence. If Bombieri's assertion is correct, it offers a more tractable route to the Riemann Hypothesis: on the critical line, \(\xi(s)\) becomes a real function of a single real variable, whose geometric behaviour is simpler than that of a complex function in the two-dimensional strip, and to which standard analytical tools can be applied. The viability of this approach is supported by Hardy's theorem that \(\xi(s)\) has infinitely many zeros on the critical line \cite{hardy}, and by the numerical verification of Platt and Trudgian \cite{platt}, which confirms the Riemann Hypothesis up to \({\rm Im}(s) \le 3 \times 10^{12}\).

We follow Bombieri's approach in this paper.

In Section~2, we simplify the complex function \(\xi(s)\) using Jensen's method, adopting notation that will be convenient for the subsequent analysis. We then collect several key properties of the Jensen function that will be needed for the study of \(\xi(s)\) on the critical line.

In Section~3, we prove that on the critical line — where \(\xi(s)\) becomes a real function of a single real variable — it satisfies a special differential equation. As a consequence, the curve near each zero is either increasing and concave down, or decreasing and concave up; this forces all zeros to be simple, with all local maxima positive and all local minima negative. Thus Bombieri's equivalence condition is satisfied.

In Section~4, since a published proof of the sufficiency direction of Bombieri's theorem appears to be unavailable, we supply an independent proof of this implication. Using the Cauchy--Riemann equations, we show that Bombieri's condition on the critical line forces \(\xi(s)\) to have no zeros off the line. The Riemann Hypothesis follows.

Finally, in Section~5 and~6, we discuss P\'olya's counterexample, which has sometimes been thought to undermine Jensen-type approaches, as well as other known counterexamples among L-functions. We show that neither P\'olya's example nor these L-function counterexamples invalidate Bombieri's equivalence theorem or the arguments of this paper, and hence they pose no conflict with our results.
\\

\section{Integral form of $\xi(s)$ and Jensen's function}
\setcounter{equation}{0}

The Riemann \(\xi(s)\) function has two equivalent representations, \eqref{eq:riemann_xi_1} and \eqref{eq:riemann_xi_2}. Our analysis focuses primarily on the latter, though we shall also use \eqref{eq:riemann_xi_1} in Lemma~1 for a convenient asymptotic estimate. In his original paper \cite{riemann}, Riemann transformed \eqref{eq:riemann_xi_2} into another integral form, still fairly complicated. In 1911, at the Copenhagen Congress, Jensen \cite{jensen} recast it as an elegant Fourier cosine transform. This form was later used by many authors, including Titchmarsh \cite{titchmarsh30}, Wintner \cite{wintner}, Haviland \cite{haviland}, Spira \cite{spira}, Matiyasevich \cite{matiyasevich}, and Pólya \cite{polya26}.

We shall use the Jensen form of the \(\xi(s)\) function, with a slight change of notation that will facilitate the application of elementary calculus tools—such as differentiation and the Cauchy--Riemann equations—to the study of \(\xi(s)\). Starting from \eqref{eq:riemann_xi_2}, we introduce a new complex variable \(z = x + {\rm i}y\) by
\begin{equation} \label{3}
        s = \frac{1}{2} (z + 1).
\end{equation} 
Then \(\xi(s)\) takes the form
\begin{equation}\label{4}
    \xi(s)  = \Xi(z) = \frac{1}{2} + \frac{1}{8} \big(z^2 - 1\big) \int\limits_1^{\infty} \Psi(\tau)\tau^{-3/4} \big(\tau^{z/4} + \tau^{-z/4} \big) {\rm d} \tau.
\end{equation}
This parametrization differs from those used by Jensen \cite{jensen, polya27} and Landau \cite{landau}, who set \(s = \frac{1}{2} + {\rm i}z\). Under our change of variables, the original strip \(0 < {\rm Re}(s) < 1\) becomes \(-1 < {\rm Re}(z) < 1\), and the critical line \({\rm Re}(s) = \frac{1}{2}\) becomes \({\rm Re}(z) = x = 0\).

We next change the integration variable from \(\tau\) to \(t\) by setting
\begin{equation} \label{5}
    \tau = e^{4t},
\end{equation}
which transforms \eqref{4} into
\begin{equation} \label{6}
    \Xi(z)  =  \frac{1}{2} + \frac{1}{2} \big(z^2-1\big) \int\limits_0^{\infty} \Psi(e^{4t}) e^{t} \big(e^{zt} + e^{-zt} \big) {\rm d}t.
\end{equation}
Using \(\cosh(w) = \frac{1}{2}(e^w + e^{-w})\), we obtain
\begin{equation} \label{7}
    \Xi(z)  =  \frac{1}{2} + \big(z^2-1\big) \int\limits_0^{\infty} F(t) \cosh(zt) {\rm d}t,
\end{equation}
where
\begin{equation} \label{8}
        F(t) = \Psi(e^{4t}) e^t = \sum\limits_{n=1}^{\infty} \exp\big(t - \pi n^2 e^{4t}\big).
\end{equation}

We continue by differentiating \(F(t)\). Its first and second derivatives are
\begin{equation} \label{F(0)}
    F'(t) =  \sum_{n = 1}^{\infty} \exp \big(t - \pi n^2 e^{4t} \big) \; \big(1 - 4\pi n^2 e^{4t}\big),
\end{equation}
\begin{eqnarray} \label{F'(0)}
    F''(t) &=& \sum_{n = 1}^{\infty} \exp \big(t - \pi n^2 e^{4t} \big) \;\Big[ \big(1 - 4\pi n^2 e^{4t}\big)^2 - 16\pi n^2 e^{4t} \Big]
            \nonumber\\[3mm]
            &=& F(t) + G(t),
\end{eqnarray}
where
\begin{equation} \label{G_func}
    G(t) = \sum_{n = 1}^{\infty} \exp \big(t - \pi n^2 e^{4t} \big) \;\Big[16\pi^2 n^4 e^{8t} - 24\pi n^2 e^{4t}\Big]
\end{equation}
is the Jensen function \cite{jensen}, up to a constant factor of \(8\) relative to Jensen's original definition. At \(t = 0\), we have \cite{lin}:
\begin{equation} \label{13}
    F(0) = \sum_{n = 1}^{\infty} \exp \big( - \pi n^2 \big) = \frac{1}{2}\big(\theta(1) - 1\big),
\end{equation} 
\begin{equation} \label{14}
    F'(0) = \sum_{n = 1}^{\infty} \exp \big( - \pi n^2 \big)\;\big(1 - 4\pi n^2\big) = -\frac{1}{2},
\end{equation}
and as \(t \to \infty\),
\begin{equation} \label{15}
        \lim_{t\to \infty} F(t)\sinh(zt) = 0, \quad \lim_{t\to \infty} F'(t)\cosh(zt) = 0,
\end{equation} 
where \(\theta(1)\) in \eqref{13} is a finite constant (see Whittaker and Watson \cite[§21]{whittaker}); the identities \eqref{14} and \eqref{15} go back to Riemann and are essential for handling the integral representation of \(\xi(s)\).

With these preparations, we apply integration by parts twice to \eqref{7}:
\begin{eqnarray} \label{16}
        \int\limits_{0}^{\infty} F(t) \, \cosh(zt) {\rm d}t &=&
        \frac{1}{z^2}F'(0) + \frac{1}{z^2} \int\limits_{0}^{\infty} [F(t)+G(t)] \, \cosh(zt) {\rm d}t.
\end{eqnarray} 
Thus,
\begin{equation} \label{17}
        (z^2 - 1) \int\limits_{0}^{\infty} F(t) \, \cosh(zt) {\rm d}t = F'(0) + \int\limits_{0}^{\infty} G(t) \, \cosh(zt) {\rm d}t.
\end{equation}
Substituting \eqref{17} into \eqref{7} and using \(F'(0) = -1/2\), we obtain the integral representation
\begin{equation} \label{xi_G}
        \Xi(z) = \int\limits_{0}^{\infty} G(t) \, \cosh(zt) {\rm d}t.
\end{equation}
This formula agrees with Jensen's original form \cite{jensen} up to notation and a constant scaling of \(G(t)\). We shall use it exclusively in the remainder of the paper.

We now recall some properties of the Jensen function. It is immediate from \eqref{G_func} that \(G(t) > 0\) for \(t > 0\), and that \(G(t) \to 0\) exponentially fast as \(t \to \infty\). Wintner \cite{wintner} proved in a short note that \(G(t)\) is strictly decreasing (\(G'(t) < 0\)); this was later rediscovered independently by Spira \cite{spira}. Wintner also noted an earlier observation of Jensen and Hurwitz: if \(G(t)\) is extended to negative \(t\), it becomes even, i.e. \(G(-t)=G(t)\), and in particular \(G'(0)=0\).

Since this evenness is not immediately obvious from \eqref{G_func}, we briefly recall the argument. With \(\tau = e^{4t}\), we have from \eqref{5} and \eqref{8}:
\begin{equation}
        F(t) = \Psi(\tau) e^t, \quad \tau(t) = e^{4t}, \quad \tau'(t) = 4\tau(t), \quad \tau(-t) = \frac{1}{\tau(t)},
\end{equation}
and a direct differentiation gives
\begin{equation} \label{G_tau}
    G(t) = F''(t) - F(t) = \Big[ 16\tau^2 \Psi''(\tau) + 24\tau \Psi'(\tau) \Big] e^t.
\end{equation}
Differentiating the Jacobi identity
\begin{equation}
    2\Psi(\tau) + 1 = \tau^{-1/2} \Big[ 2\Psi\Big(\frac{1}{\tau}\Big) + 1 \Big],
\end{equation}
which is due to Riemann \cite{riemann}, yields
\begin{equation}
     \Big[ 16\tau^2 \Psi''(\tau) + 24\tau \Psi'(\tau) \Big] e^t =
     \Big[ \frac{16}{\tau^2} \Psi''\Big(\frac{1}{\tau}\Big) + \frac{24}{\tau} \Psi'\Big(\frac{1}{\tau}\Big) \Big] e^{-t} = G(-t).
\end{equation}
Thus \(G(t)\) is even and has a bell-shaped graph (see Figure~1, computed in \cite{lin}). This Gaussian-like behaviour has motivated several constructions of analogous functions in the literature; see, for example, Pólya \cite{polya26}.

\input {./fig_01.tex}

Since \(G(t)\) is smooth and even, Jensen \cite{polya27} showed that its derivatives at \(t=0\) and at infinity satisfy
\begin{equation} \label{G^(2k)}
            G^{(2k - 1)}(0) = 0, \hspace{1cm} |G^{(2k)}(0)| < \infty,
            \hspace{1cm} \lim\limits_{t\to\infty}  G^{(k)}(t) = 0.
\end{equation}
Although Jensen did not give an explicit formula for the even-order derivatives, it follows from \eqref{eq:riemann_xi_2} and \eqref{G_tau} that \(G^{(2k)}(0)\) is expressible in terms of \(\theta^{(2k+2)}(1), \theta^{(2k+1)}(1), \ldots, \theta'(1)\). A general formula for the derivatives of \(\theta(\tau)\) at \(\tau=1\) was  found by Romik \cite{romik}; using his result, we computed the numerical values \cite{lin}
\begin{equation} \label{G(0)G''(0)}
    G(0) \approx 3.5736, \hspace{5mm} G''(0) \approx -267.6880,  \hspace{5mm} G^{(4)}(0) \approx 51978.4213, \ \cdots.
\end{equation}

Jensen \cite{polya27} and others discovered many further properties of \(G(t)\) and its derivatives. These were later collected and systematically summarized by Gélinas in a reading note, where he also derived the following useful estimates \cite[§2.27, §2.29]{gelinas21a}. For \(t \geq 0\), with \(\tau = e^{4t}\),
\begin{equation} \label{Jensen_1}
        (2\pi \tau - 3) \pi \tau e^{t-\pi \tau} \ <\ \frac{1}{8} G(t)\ <\ (2\pi \tau - 3) \pi \tau e^{t-\pi \tau}  + 32\pi^2 \tau^2 e^{t-4\pi \tau},
\end{equation}
and
\begin{equation} \label{Jensen_2}
        4G(t) \  < \ \frac{\ - G'(t)\ }{\pi \tau - \pi} \ < \ 4G(t) + 8 (8\pi - 18) \pi \tau e^{t-\pi \tau}.
\end{equation}
These inequalities will be essential in the arguments that follow.

We will use \eqref{xi_G} to study the Riemann Hypothesis. Writing \(z = x + {\rm i}y\), we separate \(\Xi(z)\) into real and imaginary parts:
\begin{eqnarray} \label{u+iv}
\Xi(z) &=& u(x,y) + {\rm i} v(x,y) \nonumber \\[2mm]
&=& \int\limits_{0}^{\infty} G(t) \cosh(xt) \cos(yt) {\rm d} t + 
{\rm i} \int\limits_{0}^{\infty} G(t) \sinh(xt) \sin(yt) {\rm d} t,
\end{eqnarray}
Thus, for \(z\) in the critical strip \(-1 < x < 1\), we have
\[
\Xi(z) = 0 \quad \Longleftrightarrow \quad u(x,y) = 0 \;\text{and}\; v(x,y) = 0.
\]

\section{Bombieri's equivalence condition}
\setcounter{equation}{0}

As mentioned in the introduction, we shall approach the Riemann Hypothesis via Bombieri's equivalence theorem. In this section we restrict attention to the critical line \(x=0\), where \(\Xi(z)\) reduces to a real-valued function of a single real variable:
\begin{align} \label{u(0,y)}
        u(0,y) &= \Xi({\rm i} y) =  \int\limits_{0}^{\infty} G(t) \cos(yt) {\rm d} t =: U(y),
                \nonumber\\[0mm]
        v(0,y) &= 0.
\end{align}
We will prove that \(U(y)\) satisfies Bombieri's condition: all of its local maxima are positive and all of its local minima are negative.

\vspace{3mm}
\textbf{Lemma 1.} {\it For sufficiently large \(y\) with \(U^{(n)}(y)\neq 0\), and for each fixed \(k\ge 1\), we have
\begin{equation} \label{Deriv U^n(y)}
\frac{U^{(n+k)}(y)}{U^{(n)}(y)} = (-1)^k \left(\frac{\pi}{8}\right)^k \left(1 + \mathcal{O}(y^{-1})\right).
\end{equation}
}

{\bf Proof.} Titchmarsh \cite[p29]{titchmarsh1986} proved that,  for sufficiently large $|s|$,  \(\xi(s)\) is an entire function of order 1 with the estimation
\begin{equation} \label{Titch}
       \log \xi(s) = \frac{1}{2} s \log s  + \cdots.
\end{equation}
Putting \(s = (1 + {\rm i}y)/2\) with \(y\) large and taking the appropriate real part gives
\begin{eqnarray}
   \log U(y) &=& \log \xi\Big(\frac{1+{\rm i}y}{2}\Big)
            =  - \frac{y}{4} \arctan y + \frac{1}{4} \log \sqrt{1+y^2}  + \cdots \nonumber \\[2mm]
            &=& -\frac{\pi}{8} y +\frac{1}{4} \log y + \cdots.
\end{eqnarray}
Differentiating with respect to \(y\) gives
\begin{equation}
        \frac{U'(y)}{U(y)} = -\frac{\pi}{8} + \frac{1}{4y} + \cdots,
\end{equation}
Now differentiate both sides of this equation and use the identity
\begin{equation}
\frac{\rm d}{{\rm d}y}\left(\frac{U'}{U}\right) = \frac{U''}{U} - \left(\frac{U'}{U}\right)^2
\end{equation}
to obtain
\begin{equation}
\frac{U''}{U} = \Big(-\frac{\pi}{8} + \frac{1}{4y} + \cdots\Big)^2 - \frac{1}{4y^2} + \cdots
= \left(\frac{\pi}{8}\right)^2  - \frac{\pi}{16y}  + \cdots 
\end{equation}
Repeated differentiation then yields
\begin{equation} 
\frac{U^{(j)}(y)}{U(y)} = (-1)^j \left(\frac{\pi}{8}\right)^j  \big(1 + \mathcal{O}(y^{-1})\big).
\end{equation}
Since  \(U^{(n)}(y)\neq 0\), setting $j = n$ and $n+k$  to get two estimates, then their ratio yields  \eqref{Deriv U^n(y)}.  

\hfill $\square$
\vspace{2mm}

{\bf Remark.} Lemma~1 is stated in terms of a ratio of two derivatives, with the condition that the denominator \(U^{(n)}(y)\neq 0\). This formulation avoids the need for \(U(y)\neq 0\). If \(U(y)\) itself vanishes, one can apply the lemma to a neighbourhood where the relevant derivatives do not vanish, and then obtain the limiting value by continuity.

\vspace{3mm}
{\bf Lemma 2.} {\it All zeros of the function \(U(y)\) defined in \eqref{u(0,y)} are simple.}

{\bf Proof.} The proof is somewhat lengthy, so we first outline the strategy. 
Suppose \(U\) has a zero of order \(n \ge 1\) at \(y_0\). Its Taylor expansion is
\begin{equation} \label{U^{(n)}}
    U(y) = \frac{1}{n!} U^{(n)}(y_0) (y - y_0)^n \Big(1 + \mathcal{O}(y-y_0)\Big),
\end{equation}
with \(U^{(n)}(y_0) \neq 0\). Then, on the right interval \(y_0 < y < y_0 + \delta\), 
all derivatives \(U^{(k)}(y)\) (\(1 \le k \le n\)) have the same sign as \(U(y)\). 
We will use this fact to derive a contradiction for every \(n \ge 2\), 
thereby proving that every zero of \(U\) is simple.

Differentiating \(U(y)\) in \eqref{u(0,y)} and integrating by parts gives
\begin{equation} \label{yU'(y) + U(y)}
        yU'(y) + U(y) = -\int\limits_0^\infty G'(t)t \cos(yt) \;\!{\rm d} t.
\end{equation}
Differentiating twice yields the second derivative
\begin{equation} \label{U''(y)}
        U''(y)  = -\int\limits_0^\infty G(t)t^2 \cos(yt) \;\!{\rm d} t.
\end{equation} 
Thus, to relate \(U(y)\) to its derivatives, we need a connection between \(G'(t)t\) and \(G(t)t^2\). Using the left inequality in \eqref{Jensen_1}, we have
\begin{equation}
        \pi \tau e^{t - \pi \tau}  <  \frac{1}{\ 8 (2\pi \tau - 3)\  }  G(t).
\end{equation}
Substituting this into the right inequality in \eqref{Jensen_2} gives
\begin{equation} \label{Jensen_3}
        4G(t)  < \frac{\ - G'(t)\ }{\pi \tau - \pi} < 4G(t) + \frac{\ 8\pi - 18\ }{2\pi \tau - 3} G(t).
\end{equation}
Since \(\tau = e^{4t}\), we rewrite \eqref{Jensen_3} as
\begin{equation} \label{Jensen_4}
        16\pi  \;\!G(t) <  - \frac{G'(t)}{t} < 16\pi \frac{(e^{4t} -1)}{4t}\bigg[1 + \frac{\ 4\pi - 9\ }{4\pi e^{4t} - 6} \bigg]G(t).
\end{equation}
This indicates that \(-G'(t)/t\) is at least \(16\pi\) times \(G(t)\), a fact that will be used later in Section~5 when we compare the present setting with P\'olya's counterexample. Figure~2 (from \cite{lin}) compares the three curves and confirms the estimates \eqref{Jensen_1} and \eqref{Jensen_2}.

\input {./fig_02.tex}

We now introduce the function \(a(t)\) by
\begin{equation}  \label{Taylor_1}
       \frac{-G'(t)}{t\, G(t)} =  a(t) =a_0 + a_2 t^2 + a_4 t^4 + a_6 t^6 + \cdots.
\end{equation}
Since both \(G(t)\) and \(-G'(t)/t\) are even, \(a(t)\) is even as well. Taking \(t \to 0\) in \eqref{Taylor_1} gives (see \cite{lin})
\begin{equation}  \label{Taylor_3}
        a_0  =\lim_{t\to 0} \frac{-G'(t)}{t G(t)} = - \frac{G''(0)}{ G(0)} \approx 1.4902 \times 16\pi \approx 74.9076 > 0.
\end{equation} 
Using the values of \(G(0)\), \(G''(0)\), and \(G^{(4)}(0)\) from \eqref{G(0)G''(0)}, we obtain (see \cite{lin})
\begin{eqnarray} \label{a_2_limit}
            a_2 & = &
            \lim_{t\to 0} \frac{1}{t^2} \Big(
      - \frac{G'(t)}{tG(t)} -a_0
                    \Big)
                    =     \frac{ 3[G''(0)]^2 - G(0) G^{(4)}(0) } { 6 [G(0)]^2 }
                    \nonumber\\[2mm]
                     & \approx& 381.3732 = 5.0912a_0.
\end{eqnarray}
Continuing this process to compute higher coefficients would be tedious. Since our main interest lies in the zeros of \(U(y)\) for large \(y\), however, an approximate function \(\hat{a}(t)\) suffices to capture the essential relation between \(G'(t)t\) and \(G(t)t^2\). We therefore construct the coefficients \(\{\hat{a}_{2k}\}\) by a different method.


Let us analyze the inequality on the right-hand side of \eqref{Jensen_4}. The comparison in Figure~2 shows that the two curves involved are very close. This suggests that the inequality can be expressed as an equality through the introduction of a smooth function \(b(t)\):
\begin{equation} \label{Jensen_5}
        a(t) = \frac{-G'(t)}{t\, G(t)} = 16\pi\left[1 + b(t) \frac{\,4\pi - 9\,}{4\pi e^{4t} - 6} \right]\frac{(e^{4t} -1)}{4t}.
\end{equation}
The discrepancy between the curves is greatest at \(t=0\), suggesting that \(b(0) < b(t) < 1\). To determine the value of \(b(0)\), we take the limit as \(t \to 0\) in \eqref{Jensen_5} and use the known value of \(a_0\):
\begin{equation}
        a_0 =  16\pi \left[1 + b(0) \frac{\,4\pi - 9\,}{4\pi  - 6} \right].
\end{equation}
Solving this gives \(b(0) \approx 0.9026\), so \(b(t)\) varies only slightly. 

Then, from \eqref{Jensen_5} and the fact that \(b(t)\) varies only slightly, we observe that the behavior of \(a(t)\) is primarily governed by the term \((e^{4t}-1)/(4t)\). Since this term is bounded between the two even functions \(\sinh(4t)/(4t)\) and \(\cosh(4t)\), we obtain the bounds
\begin{equation}  \label{lele-0}
        a_0 \frac{\sinh(4t)}{4t} \le a(t) \le a_0 \cosh(4t),
\end{equation}
as illustrated in Figure~3 (source: \cite{lin}).

\input {./fig_03.tex}

Therefore, performing Taylor expansions for \(\sinh(4t)/(4t)\) and \(\cosh(4t)\) and comparing coefficients in \eqref{Taylor_1} yields
\begin{equation} \label{lele}
    a_0 \sum_{k=0}^{\infty} \frac{4^{2k}}{(2k+1)!} t^{2k} \le
    \sum_{k=0}^{\infty} a_{2k} t^{2k} \le
    a_0 \sum_{k=0}^{\infty} \frac{4^{2k}}{(2k)!} t^{2k}.
\end{equation}
Guided by this inequality, we shall use the following approximation for the coefficients:
\begin{equation} \label{lele-1}
    a_0 \frac{4^{2k}}{(2k+1)!} \le
    \hat{a}_{2k} \le
    a_0 \frac{4^{2k}}{(2k)!}.
\end{equation}
This allows us to choose \(\hat{a}_{2k}\) as a weighted average of the two bounding sequences:
\begin{equation} \label{a_2k<}
    \hat{a}_{2k} = (1 - c_{2k}) a_0 \frac{4^{2k}}{(2k+1)!} + c_{2k} a_0 \frac{4^{2k}}{(2k)!},
    \qquad (k = 0, 1, 2, \ldots),
\end{equation}
where \(0 \le c_{2k} \le 1\).

We now choose the coefficients \(c_{2k}\) to suit our purpose. The three series in \eqref{lele} share the same value \(a_0\) at \(t=0\), so \(c_0\) is arbitrary; we set \(c_0 = 1/2\). As \(t \to \infty\), from \eqref{Jensen_5} and \eqref{lele-0},
\begin{equation} \label{a_2k_lim}
        \lim_{t\to\infty} \frac{ a(t) } { a_0 \sinh(4t) / (4t)} = \frac{32\pi}{a_0}, \qquad
        \lim_{t\to\infty} \frac{ a(t) } { a_0 \cosh(4t)} = 0.
\end{equation}
Thus the approximating series \(\sum \hat{a}_{2k} t^{2k}\) should closely follow \(a_0 \sinh(4t)/(4t)\) for large \(t\), which means that only finitely many \(c_{2k}\)'s need to be nonzero. We choose
\begin{equation} \label{c_coef}
    c_0 = \frac12, \qquad c_2 = \frac14, \qquad c_4 = \frac{1}{10}, \qquad c_{2k}=0 \quad (k \ge 3).
\end{equation}
which gives \(\hat{a}_0 = a_0, \, \hat{a}_2 = 4a_0, \cdots\). The resulting fitted curve (from~\cite{lin}) is also shown in Figure~3. 

To quantify the accuracy of this approximation, define the fitting error
\begin{equation} \label{a(t)-Delta}
\Delta_T = \max_{0\le t < T} \bigg| a(t) \Big( \sum_{k=0}^{\infty} \hat{a}_{2k} t^{2k} \Big)^{-1} - 1 \bigg|.
\end{equation}
Since $G(t)$ decays faster than any exponential for $t>1$ (indeed $G(1)/G(0) \sim 10^{-70}$), the ratio $a(t)/\hat{a}(t)$ outside $[0,1]$ is negligible for the Fourier integrals; only the neighborhood of $t=0$ matters for the limit $y\to\infty$. We therefore set $T=1$ and denote $\Delta:=\Delta_1$. A direct numerical comparison (see the supplementary material in \cite{lin}) gives $\Delta \le 0.03$ for the coefficients $c_{2k}$ chosen in \eqref{c_coef}. Consequently, defining
\begin{equation} \label{A(t)}
\hat{A}(t) = a(t) \Big( \sum_{k=0}^{\infty} \hat{a}_{2k} t^{2k} \Big)^{-1},
\end{equation}
we have
\begin{equation} \label{A(y)O}
1-\Delta \le \hat{A}(t) \le 1+\Delta, \qquad (0\le t \le 1).
\end{equation}

Using the function \(a(t)\) and the relation \(-G'(t)t = G(t) t^2 a(t)\), we can combine \(U(y),\) \(U'(y),\) and \(U''(y)\). From \eqref{yU'(y) + U(y)} and \eqref{A(t)}, we have
\begin{eqnarray} \label{A(y)}
      yU'(y) + U(y) &=& -\int\limits_0^\infty G'(t)t \cos(yt) \;\!{\rm d} t =
      \int\limits_0^\infty G(t)\;\! t^2 \;\!a(t)\cos(yt) \;\!{\rm d} t
      \nonumber\\[1mm]
      &=& \int\limits_0^\infty G(t)\;\! t^2 \;\!  \hat{A}(t)
      \Big( \sum\limits_{k = 0}^\infty \hat{a}_{2k} t^{2k}\Big) \cos(yt) \;\!{\rm d} t
            \nonumber\\[1mm]
      &=:& A(y) \int\limits_0^\infty G(t)\;\! t^2 \;\!  
      \Big( \sum\limits_{k = 0}^\infty \hat{a}_{2k} t^{2k}\Big) \cos(yt) \;\!{\rm d} t,
\end{eqnarray}
where \(A(y)\) is a weighted mean of \(\hat{A}(t)\). For any \(y\) such that \(yU'(y) + U(y) \not= 0\), \(A(y)\) is well defined.  Since \(\hat{A}(t)\) is continuous and the integrand \(G(t)t^2( \sum \hat{a}_{2k} t^{2k} )> 0\), the standard asymptotic theory of Fourier-type integrals with smooth, rapidly decaying amplitudes (see, e.g., \cite{erdelyi}) gives
\begin{equation} \label{A(y)lim}
       A(y) =  \frac{ \displaystyle \int\limits_0^\infty G(t)\;\! t^2 \;\!  \hat{A}(t)
      \Big( \sum\limits_{k = 0}^\infty \hat{a}_{2k} t^{2k}\Big) \cos(yt) \;\!{\rm d} t}
                { \displaystyle  \int\limits_0^\infty G(t)\;\! t^2 \;\!  
      \Big( \sum\limits_{k = 0}^\infty \hat{a}_{2k} t^{2k}\Big) \cos(yt) \;\!{\rm d} t} \to \hat{A}(0) =: A, \hspace{5mm} 
            ({\rm as~} y \to \infty),
\end{equation}
Thus, for sufficiently large \(y\) with \(yU'(y) + U(y) \not= 0\), 
\begin{equation}
            A(y) = A + {\mathcal O}(y^{-1}),
\end{equation}
where \(1-\Delta \le A  \le 1+\Delta\). Substituting this into \eqref{A(y)} yields
\begin{eqnarray}
      yU'(y) + U(y) = \big(A+ {\mathcal O}(y^{-1})\big)\int\limits_0^\infty G(t) \Big( \hat{a}_0t^2 + \hat{a}_2 t^4 + \hat{a}_4 t^6 +\hat{ a}_6 t^8 + \cdots \Big) \cos(yt) \;\!{\rm d} t.
\end{eqnarray}
Differentiating \eqref{U''(y)} repeatedly gives, for higher-order derivatives,
\begin{equation} \label{U^{(2k)}(y)}
              U^{(2k)}(y) = (-1)^k \int\limits_0^\infty G(t) t^{2k} \cos(yt)\;\! {\rm d} t.
\end{equation}
Using this in the previous equation, we obtain the differential relation
\begin{equation} \label{ODE}
       yU'(y) + U(y) = - (A+ {\mathcal O}(y^{-1}))\Big[ \hat{a}_0 U''(y) - \hat{a}_2 U^{(4)}(y) + \hat{a}_4 U^{(6)}(y) - \cdots\Big].
\end{equation}

According to Hardy \cite{hardy}, \(U(y)\) has infinitely many zeros in \(0 < y < \infty\). Let \(y_0\) be one of them. We assume \(y_0\) is sufficiently large; indeed, Platt and Trudgian \cite{platt} have already verified the Riemann Hypothesis for \(y \le 6 \times 10^{12}\), so we only need to consider the range \(y > 6 \times 10^{12}\), where the asymptotic estimates of Lemma~1 are valid. We further assume, for contradiction, that \(y_0\) is a zero of order two. Then, near \(y_0\),
\begin{equation} \label{U(y)_A1}
    U(y)= \frac{1}{2} U''(y_0) (y-y_0)^2 \Big( 1 + \mathcal{O}(y-y_0) \Big),
\end{equation}
with \(U''(y_0) \neq 0\). 
Hence, in a sufficiently small right-neighborhood \(y_0 < y < y_0 + \delta\), we have \(U(y) \neq 0\), \(U'(y) \neq 0\), \(U''(y) \neq 0\). In this interval, \(U'(y)\) has the same sign as \(U(y)\), so \(yU'(y)+U(y) \neq 0\).

Since \(y\) is large and \(U''(y)\neq 0\) in this neighbourhood, Lemma~1 (with \(n=2\)) gives
\begin{equation} \label{U^{(2k+2)}(y)}
            U^{(2k+2)}(y) = \Big(\frac{\pi}{8}\Big)^{2k} U''(y) \Big(1 + {\mathcal O}(y^{-1})\Big),
            \hspace{5mm}  (k = 1, 2, 3, \cdots).
\end{equation}
Substituting \eqref{U^{(2k+2)}(y)} into \eqref{ODE}, we obtain
\begin{equation} \label{ODE_1}
       yU'(y) + U(y) \ = \ -A c\, U''(y) + {\mathcal O}\Big( \frac{ U''(y)}{y} \Big),
\end{equation}  
where
\begin{equation} \label{c}
    c = \hat{a}_0 - \hat{a}_2 \Big(\frac{\pi}{8}\Big)^2
            + \hat{a}_4 \Big(\frac{\pi}{8}\Big)^4
            - \hat{a}_6 \Big(\frac{\pi}{8} \Big)^6 + \cdots.
\end{equation}
The series in \eqref{c} is alternating. From \eqref{lele-1}, it satisfies the Leibniz conditions:
\begin{eqnarray}
\label{a_2<a_0}
            \hat{a}_2 \Big(\frac{\pi}{8} \Big)^2 &=& 4 a_0 \Big(\frac{\pi}{8}\Big)^2 \ < \ \hat{a}_0,
\\[3mm]
        \hat{a}_{2k+2} \Big(\frac{\pi}{8}\Big)^{2k+2}
            & \le & a_0 \frac{4^{2k+2}}{(2k+2)!} \Big(\frac{\pi}{8} \Big)^{2k+2}
            \nonumber \\[2mm]
             & = &  a_0 \frac{4^{2k}}{(2k+1)!} \Big(\frac{\pi}{8} \Big)^{2k} \frac{1}{2k+2} \Big(\frac{\pi}{2}\Big)^2
             \nonumber \\[2mm]
             &\le& \hat{a}_{2k} \Big(\frac{\pi}{8} \Big)^{2k},
             \hspace{10mm} (k = 1,2,3,\cdots),
\end{eqnarray}
and
\begin{eqnarray}
        \lim_{k\to\infty} \hat{a}_{2k} \Big(\frac{\pi}{8}\Big)^{2k} &\le&
            \lim_{k\to\infty} a_0\frac{4^{2k}}{(2k)!} \Big(\frac{\pi}{8}\Big)^{2k} \ =\  0.
\end{eqnarray}
Thus the alternating series \eqref{c} satisfies both conditions of the Leibniz test, converges to a finite value \(c\), with
\begin{equation}
        \hat{a}_0 - \hat{a}_2 \Big(\frac{\pi}{8} \Big)^2 \le c \le \hat{a}_0.
\end{equation}
From \eqref{a_2<a_0}, we have \(c > 0\).

If some of the higher-order derivatives vanish, the corresponding terms in the alternating series are absent. The worst case for the positivity of \(c\) is the deletion of all positive terms except the first, while retaining all negative terms. Then
\begin{equation}
            c = \hat{a}_0 - \sum_{j=0}^{\infty} \hat{a}_{4j+2} \Big(\frac{\pi}{8}\Big)^{4j+2}.
\end{equation}
Using the coefficients from \eqref{c_coef}, we have \(\hat{a}_0 = a_0\), \(\hat{a}_2 = 4a_0\), and for \(k\ge 3\), \(\hat{a}_{2k} = a_0 4^{2k}/(2k+1)!\). Hence
\begin{align}
            c &= a_0 \bigg(1 - 4\Big(\frac{\pi}{8}\Big)^2 - \sum_{j=1}^{\infty} \frac{4^{4j+2}}{(4j+3)!} \Big(\frac{\pi}{8}\Big)^{4j+2}      \bigg)
            \nonumber \\[2mm]
            &> a_0 \Big(1 - \frac{\pi^2}{16} - 0.004\Big) > 0.38 a_0 > 0.
\end{align}
Thus \(c>0\) is stable under deletion of any terms. Since \(\Delta \le 0.03\), we also have \(A \ge 1 - \Delta > 0\); consequently \(Ac > 0\).

Therefore, for \(y = y_0 + \delta\) with \(\delta > 0\) sufficiently small, \eqref{U(y)_A1} shows that if \(U(y) > 0\), then \(U'(y) > 0\) and \(U''(y) > 0\); if \(U(y) < 0\), then \(U'(y) < 0\) and \(U''(y) < 0\). 
Both cases contradict \eqref{ODE_1}, since the left-hand side of \eqref{ODE_1} has the same sign as \(U(y)\) (and hence as \(U''(y)\)), while the dominant term on the right-hand side has the opposite sign because \(Ac > 0\). Thus \(y_0\) cannot be a second-order zero.

Next, we assume that \(y_0\) is a third-order zero. Near \(y_0\), \(U(y)\) has the expansion
\begin{equation} \label{60}
    U(y)= \frac{1}{3!} U'''(y_0)  (y-y_0)^3\Big( 1 + {\mathcal O} (y-y_0) \Big),
\end{equation}
with \(U'''(y_0) \neq 0\), so \(U'''(y)\) is also non-zero in a sufficiently small neighborhood of \(y_0\). One checks similarly that \(yU'(y) + U(y) \neq 0\) and \(yU''(y) + 2U'(y) \neq 0\) there.
Differentiating \eqref{ODE} with respect to \(y\) gives
\begin{equation} \label{ODE_A}
       yU''(y) + 2U'(y) \ = \ - (A+ {\mathcal O}(y^{-1})) \Big[  \hat{a}_0 U'''(y) - \hat{a}_2 U^{(5)}(y) + \hat{a}_4 U^{(7)}(y) -  \cdots \Big].
\end{equation}
From Lemma~1, for \( n = 3\) and \(k \ge 1\),
\begin{equation} \label{U^{(2k+3)}(y)}
            U^{(2k+3)}(y) = \Big(\frac{\pi}{8}\Big)^{2k} U'''(y) \Big(1 + {\mathcal O}(y^{-1})\Big),
\end{equation}
 Substituting \eqref{U^{(2k+3)}(y)} into \eqref{ODE_A} yields
\begin{equation} \label{ODE_3}
       yU''(y) + 2U'(y) \ = -A c\, U'''(y) + {\mathcal O}\Big( \frac{ U'''(y)}{y} \Big),
\end{equation}
where \(c > 0\) is the same constant as in \eqref{c}.

Therefore, for \(y = y_0 + \delta\) with \(\delta > 0\) sufficiently small, \eqref{60} shows that if \(U(y) > 0\), then \(U'(y) > 0\), \(U''(y) > 0\), and \(U'''(y) > 0\); if \(U(y) < 0\), then all three derivatives are negative. Both cases contradict \eqref{ODE_3}, since the left-hand side of \eqref{ODE_3} has the same sign as \(U'''(y)\), while the dominant term on the right-hand side has the opposite sign because \(Ac > 0\). Hence \(y_0\) cannot be a third-order zero.

\input {./fig_04.tex}

If \(U(y)\) has a zero of order \(n \ge 4\) at \(y_0\),  the expansion \eqref{U^{(n)}} can be used. Then, differentiating \eqref{ODE} \(n-2\) times, we obtain
\begin{equation} \label{ODE_n}
       yU^{(n - 1)}(y) + (n - 1)U^{(n - 2)}(y) = - (A+ \mathcal{O}(y^{-1}))\Big[ \hat{a}_0 U^{(n)}(y) - \hat{a}_2 U^{(n + 2)}(y) + \cdots \Big].
\end{equation}
On the right neighborhood of \(y_0\), the left-hand side of \eqref{ODE_n} has the same sign as \(U^{(n)}(y)\), while the dominant term on the right-hand side has the opposite sign, since \(Ac > 0\). This contradiction rules out \(n \ge 4\). Therefore, every zero of \(U(y)\) is simple.
\hfill $\square$
\vspace{2mm}

{\bf Remark.} Lemma 2 can be understood intuitively from Figure~4 (taken from \cite{lin}), which shows a portion of the \(U(y)\) curve near the 5th through 7th zeros. At the zero near \(y=75\), the curve is strictly increasing and concave down (\(U'>0,\; U''<0\)); near adjacent zeros, it is strictly decreasing and concave up (\(U'<0,\; U''>0\)). This geometric behavior is incompatible with the presence of higher-order zeros.

\vspace{2mm}

{\bf Theorem 1.} {\it The function \(U(y):=u(0,y)\) defined in \eqref{u(0,y)} satisfies Bombieri's equivalence conditions. Specifically, between each pair of adjacent zeros, \(U(y)\) has exactly one stationary point \(y_m\) with \(U'(y_m)=0\), and at this point \(U(y_m)\) is either a positive local maximum or a negative local minimum.}

{\bf Proof.} By Lemma 2, all zeros of \(U(y)\) are simple. Rolle's theorem then guarantees at least one point \(y_m\) between each pair of adjacent zeros such that \(U'(y_m)=0\). Near \(y_m\), we have
\begin{equation} \label{U'(y)_m}
    U'(y)= \frac{1}{n!} U^{(n + 1)}(y_m) (y - y_m)^n\Big( 1 + {\mathcal O} (y - y_m) \Big),
\end{equation}
where \(n\ge 1\), \(U^{(n+1)}(y_m)\neq 0\), and \(yU'(y)+U(y)\neq 0\) also holds. If \(n=2\), then \(U'''(y_m)\neq 0\). Assume without loss of generality that \(U'''(y_m)>0\). Then, by \eqref{U'(y)_m}, on the right interval \(y_m<y<y_m+\delta\), we would have \(U'(y)>0\), \(U''(y)>0\), and \(U'''(y)>0\), which contradicts \eqref{ODE_3}. The case \(U'''(y_m)<0\) is analogous. Hence \(n\neq 2\). The same argument rules out every \(n\ge 3\).

Thus the only possible case is \(n = 1\). Near \(y_m\), \eqref{ODE} reduces to
\begin{equation} \label{ODE_m}
       yU'(y) + U(y) \ = -Ac \, U''(y) + {\mathcal O}\Big( \frac{ U''(y)}{y} \Big),
\end{equation}
which is the same as \eqref{ODE_1}. Letting \(y \to y_m\) and using \(U'(y_m)=0\), we obtain
\begin{equation}
U(y_m) = -Ac \, U''(y_m) + \mathcal{O}\!\Big(\frac{U''(y_m)}{y_m}\Big).
\end{equation}
If \(U(y_m)>0\), then \(U''(y_m)<0\), so \(y_m\) is a positive local maximum. If \(U(y_m)<0\), then \(U''(y_m)>0\), so \(y_m\) is a negative local minimum. Hence \(U(y)\) has no positive local minima or negative local maxima.

Moreover, between any two adjacent zeros, the number of stationary points must be odd. If there were three or more, an intermediate stationary point would necessarily be a positive local minimum or a negative local maximum, contradicting what we have just proved. Hence exactly one stationary point lies between each pair of adjacent zeros.
\hfill $\square$
\vspace{3mm}

\begin{figure}
  \centering
  \includegraphics[width=11.1cm]{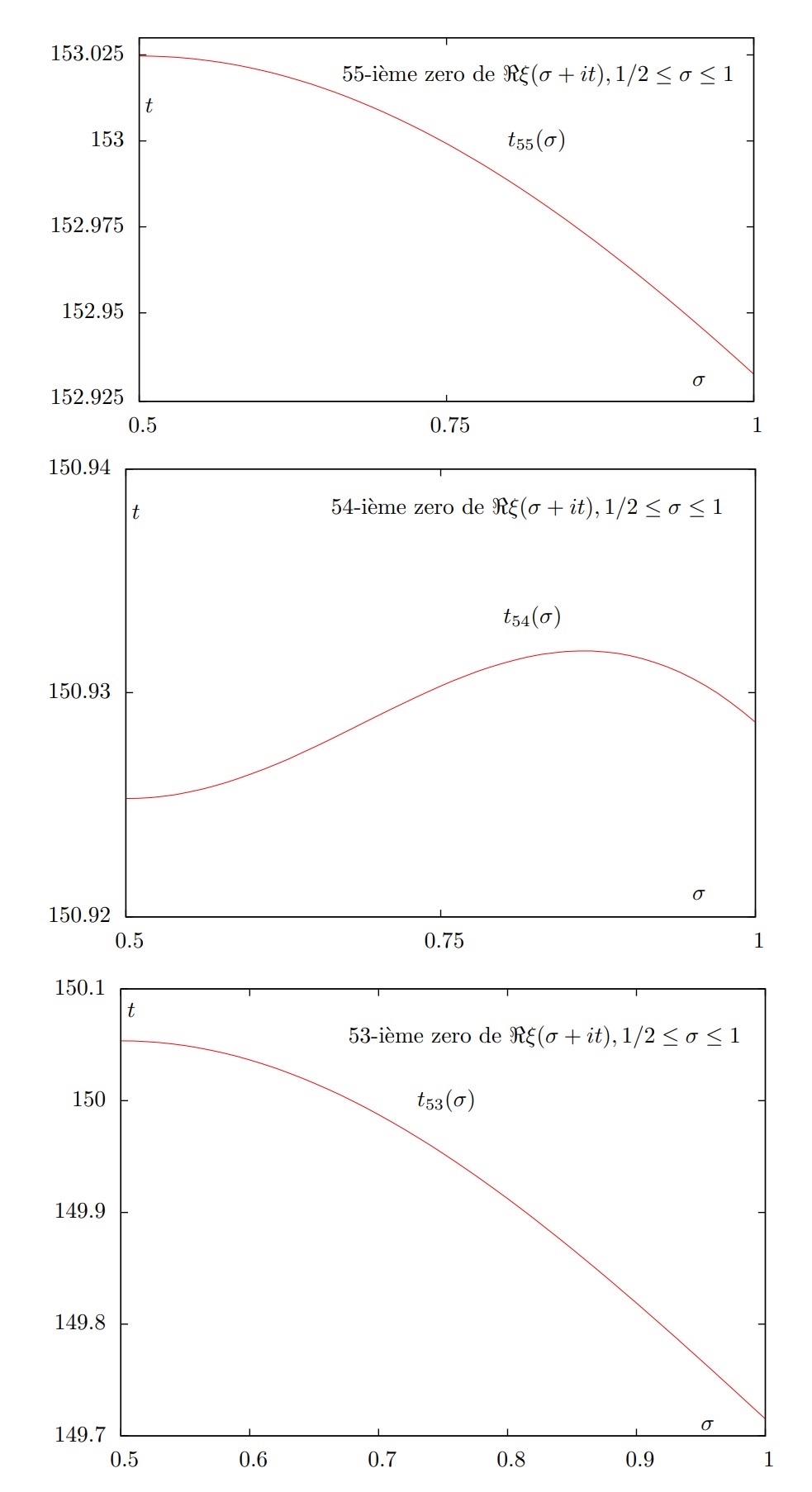}\\
  \caption{The zero contours of ${\rm Re} (\xi(s))$  in the critical strip, obtained by G\'elinas using Pari/GP and Gnuplot system in 2017. Picture permission from G\'elinas.}
\end{figure}

\section{Zeros off the critical line}
\setcounter{equation}{0}

By Theorem 1, the function \(\Xi(z)\) satisfies Bombieri's equivalence conditions on the critical line. If Bombieri's equivalence theorem is accepted, the Riemann Hypothesis follows immediately: \(\Xi(z)\) has no zeros off the critical line \( {\rm Re}(z)=0\). However, although Bombieri's theorem is stated in the official problem description \cite{bombieri}, we have not been able to locate a published proof of its sufficiency direction. For completeness, we therefore supply an independent proof of this implication in the present section. Should Bombieri's original proof come to light later, he of course deserves full credit for the theorem.

We now seek to rule out zeros of \(\Xi(z)\) off the critical line. A standard approach, described by Pólya \cite{polya26,lin}, is to construct a function \(A(z)\) such that \(|A(z)\Xi(z)-1| < 1\), which would immediately imply that \(\Xi(z)\) has no zeros. For the actual Riemann \(\Xi(z)\), however, finding such an \(A(z)\) appears to be out of reach; Pólya himself was able to construct such an \(A(z)\) only for a certain approximation to \(\Xi(z)\).

Our method is motivated by numerical experiments on the zero-level contours of the real and imaginary parts of \(\Xi(z)\) in the critical strip. Figure~5 displays some such contours of \(u(x,y)\), computed by Gélinas in 2017 using Pari/GP and Gnuplot system~\cite{lin, gelinas19}. These contours partition the critical strip into regions where \(u\neq 0\); on the contours themselves, \(u=0\). Since \(\Xi=0\) requires both \(u=0\) and \(v=0\), it suffices to show that \(v\neq 0\) on every such contour. Then the zero sets of \(u\) and \(v\) do not intersect, so \(\Xi\) has no zeros in the critical strip.

The same strategy applies with the roles of \(u\) and \(v\) interchanged: the zero contours of \(v\) also partition the critical strip, and it remains to show that \(u\neq 0\) on each such contour. As we shall see, this follows readily from Theorem 1 and the Cauchy--Riemann equations. These arguments will establish the sufficiency of Bombieri's equivalence theorem and thereby complete the proof of the Riemann Hypothesis.

\vspace{3mm}
{\bf Lemma 3.} {\it There exists a positive function \(\varepsilon(y)\) such that the two strips
\[
\omega = \{(x,y) : 0 < |x| < \varepsilon(y)\}
\]
on either side of the critical line \(x=0\) contain no zeros of \(\Xi(z)\).}

{\bf Proof.} Bohr and Landau \cite[p193-195]{edwards} proved in 1914 that almost all zeros of \(\Xi(z)\) lie within any prescribed neighbourhood \(|x| \le \varepsilon\) of the critical line. Their result, however, does not determine whether any of those zeros might actually lie off the line. The present lemma strengthens this by showing that, in fact, the regions on either side of the critical line contain no zeros at all.

\input {./fig_06.tex}

From \eqref{u+iv}, we have \(v(0,y)=0\) on the critical line. If \(v_x(0,y)\neq 0\) at some point \((0,y)\), then by continuity there is an \(\varepsilon_1>0\) such that
\begin{equation} \label{v(x,y) not=0}
v(x,y) = v_x(0,y)x + \mathcal{O}(x^2) \neq 0
\end{equation}
for all nearby points with \(0<|x|<\varepsilon_1\). The size of \(\varepsilon_1\) depends on the distance from \((0,y)\) to the nearest point \((0,y_m)\) where \(v_x(0,y_m)=0\); as long as \(y\neq y_m\), we have \(\varepsilon_1=\varepsilon_1(y)>0\).

If \((x,y)\) approaches such a point \((0,y_m)\) with \(v_x(0,y_m)=0\), then the Cauchy--Riemann equation
\begin{equation} \label{cr}
    \frac{\partial u}{\partial y} = - \frac{\partial v}{\partial x}
\end{equation}
gives \(u_y(0,y_m)=0\). By Theorem 1, \(u(0,y)\) has a local extremum at \(y_m\), so \(u(0,y_m)\neq 0\). By continuity of \(u\), there exists \(\varepsilon_2=\varepsilon_2(y_m)>0\) such that \(u(x,y)\neq 0\) whenever
\[
\sqrt{x^2+(y-y_m)^2} < \sqrt{2}\,\varepsilon_2.
\]
These two cases are illustrated in Figure~6.

Since $U(y)$ is real-analytic and not identically zero, its zeros are isolated and countable. By Theorem~1, the stationary points of \(U(y)\) are also isolated. Hence, along the critical line, consecutive stationary points \(y_m\) and \(y_{m+1}\) satisfy \(|y_{m+1}-y_m|>0\). We may therefore define a piecewise function
\begin{equation}
\varepsilon(y) = \left\{ \begin{array}{ll}
\varepsilon_2(y_m), & {\rm if\ } |y - y_m| < \varepsilon_2(y_m), \\[3mm]
\varepsilon_1(y), & {\rm otherwise.}
\end{array} \right.
\end{equation}
Then the two regions
\begin{equation}
\omega = \{(x,y) : 0 < |x| < \varepsilon(y) \}
\end{equation}
on either side of the critical line contain no zeros of \(\Xi(z)\), since in these regions either \(u\neq 0\) or \(v\neq 0\). 
\hfill $\square$
\vspace{2mm}

{\bf Lemma 4.} {\it The function \(\Xi(z)\) has no zeros with \(|x|>0\), i.e. no zeros off the critical line.}

{\bf Proof.} This is the final step in the proof of the Riemann Hypothesis. 
We shall use the zero contours of \(v(x,y)\) to partition the half-plane \(x>0\) into subregions. In each such subregion, \(v\neq 0\); it remains to show that on the contours themselves, where \(v=0\), we have \(u\neq 0\). By the symmetry of \(\Xi(z)\), it suffices to consider the right half-plane \(x>0\) (which contains the strip \(0<x<1\)).

By Lemma 3, there is an \(\varepsilon\)-region \(\omega = \{(x,y): 0 < x < \varepsilon(y)\}\) adjacent to the critical line. Since \(v_x(0,y)\) changes sign at the points where \(v_x(0,y_m)=0\), the region \(\omega\) is partitioned into subregions \(\omega_k\) (\(k = \cdots, -1, 0, 1, \cdots\)) on which \(v\) has constant sign; adjacent subregions have opposite signs. The boundary between \(\omega_k\) and \(\omega_{k+1}\) is a curve \(y=\varphi_k(x)\) on which \(v(x,\varphi_k(x))=0\). This curve emanates from a point \((0,y_m)\) where \(v_x(0,y_m)=0\), equivalently \(u_y(0,y_m)=0\) by the Cauchy--Riemann equations. At such a point, Theorem~1 shows that \(u(0,y_m)\) is either a positive local maximum or a negative local minimum.

Specifically, we write
\begin{equation}
        \omega_k =  \{(x,y) : 0<x<\varepsilon(y), \ \varphi_{k - 1}(x) < y < \varphi_k(x), \ (-1)^k v(x,y) >0 \}.
\end{equation}
Since \(v\) is analytic for \(x>0\), each \(\omega_k\) extends naturally to a full subregion of the half-plane:
\begin{equation}
    \Omega_k =  \{(x,y) : x > 0, \ \varphi_{k - 1}(x) < y < \varphi_k(x), \ (-1)^k v(x,y) >0 \}.
\end{equation}
The boundary between \(\Omega_k\) and \(\Omega_{k+1}\) is again the curve \(y=\varphi_k(x)\), on which \(v=0\). These curves are the zero contours of \(v\). 

For clarity, we recall that the zero set of \(v\) consists of curves
\begin{equation}
\gamma_k : y = \varphi_k(x), \qquad v(x,\varphi_k(x)) = 0.
\end{equation}
These are the zero-level contours of \(v\). Except at points where \(\nabla v=0\), the implicit function theorem guarantees that each contour is a smooth curve without intersections or branches. At points where \(\nabla v=0\), the local structure  is determined by the higher-order Taylor expansion of \(v\); in all cases, the zero set is a finite union of smooth arcs. Thus the following three geometric possibilities exhaust all configurations: (i) the contours neither intersect nor branch; (ii) some contours intersect; (iii) a contour branches.

We shall show that cases (ii) and (iii) cannot occur, and that in case (i), \(u\neq 0\) on every contour. Hence \(\Xi(z)\neq 0\) throughout the half-plane.

\input {./fig_07.tex}

In case (i), suppose that none of the curves \(y=\varphi_k(x)\) intersect or branch in the half-plane \(x>0\) (see Figure~7 for a typical configuration). In a subregion \(\Omega_{2k}\), we have \(v>0\). On its left boundary \(x=0\), the Cauchy--Riemann equation gives \(-u_y=v_x>0\), so \(u(0,y)\) is decreasing along this boundary. Thus \(u(0,y)\) attains a minimum value \(u_{\min}<0\) at the endpoint where \(\Omega_{2k}\) meets \(\Omega_{2k+1}\), i.e. at the starting point of the curve \(y=\varphi_{2k}(x)\).

Along the curve \(y=\varphi_{2k}(x)\), define the tangent and normal vectors
\begin{equation}
\mathbf{e}_1 = \bigl(1,\varphi_{2k}'(x)\bigr), \qquad
\mathbf{e}_2 = \bigl(-\varphi_{2k}'(x),1\bigr).
\end{equation}
Since \(v>0\) in \(\Omega_{2k}\) and \(v<0\) in \(\Omega_{2k+1}\), we have the directional derivative of \(v\) in the \(\mathbf{e}_2\) direction  $\nabla v \cdot \mathbf{e}_2 \le 0$. Using the Cauchy--Riemann equations, we compute
\begin{eqnarray}
0 \ge \nabla v \cdot \mathbf{e}_2 &=& \left( \frac{\partial v}{\partial x}, \frac{\partial v}{\partial y} \right) \cdot \bigl(-\varphi_{2k}'(x), 1\bigr)
= -\varphi_{2k}'(x) v_x + v_y
\nonumber\\[3mm]
&=& \varphi_{2k}'(x) u_y + u_x
= \nabla u \cdot \mathbf{e}_1.
\end{eqnarray}
Hence \(u\) is  non-increasing in the \(\mathbf{e}_1\) direction, so \(u(x,\varphi_{2k}(x)) \le u_{\min} < 0\). Similarly, along the next curve \(y=\varphi_{2k+1}(x)\), we obtain \(u(x,\varphi_{2k+1}(x)) \ge u_{\max} > 0\).

In summary, for any \(z=x+{\rm i} y\) with \(x>0\), either \((x,y)\) lies in some \(\Omega_k\), where \(v\neq 0\), or it lies on a boundary curve \(y=\varphi_k(x)\), where we have just shown that \(u\neq 0\). In either case, \(\Xi(z)\neq 0\).

For cases (ii) and (iii), we first prove that the inequalities above are in fact strict,  so that \(\nabla u \not= 0\) on every contour.   
Take any $0 < x_1 < x_2 $ and construct  a curvilinear quadrilateral \(Q\) with vertices
\[
        \begin{array}{cc}
                P_2=(x_1,\varphi_{2k+1}(x_1)),  & P_3=(x_2,\varphi_{2k+1}(x_2)), \\[2mm]
                P_1=(x_1,\varphi_{2k}(x_1)), & P_4=(x_2,\varphi_{2k}(x_2))
        \end{array}
\]
in the subregion \(\Omega_{2k+1}\). The left and right boundaries are vertical segments, while the upper and lower boundaries are arcs of the contours \(\varphi_{2k+1}\) and \(\varphi_{2k}\), respectively. On the lower boundary \(P_1P_4\), we have \(u\le u_{\min}<0\); on the left boundary \(P_1P_2\), \(u_{\min}\le u\le u_{\max}\); on the upper boundary \(P_2P_3\), \(u\ge u_{\max}>0\). Since \(u\) is harmonic and not constant, the maximum on \(\partial Q\) must occur on the upper boundary, and by the uniqueness principle it must be strictly greater than \(u_{\max}\) at every interior point of \(P_2P_3\). Hence \(u\) is strictly increasing along \(\varphi_{2k+1}\) from left to right. Similarly, the minimum occurs on the lower boundary, and \(u\) is strictly decreasing along \(\varphi_{2k}\). Thus \(\nabla u\cdot {\bf e_1} \neq 0\), i.e.  \(\nabla u\neq 0\) on every zero contour of \(v\).

Then, by the Cauchy--Riemann equations, \(\nabla v\neq 0\) on these contours as well. 
 Therefore, 
by the implicit function theorem (see, e.g., \cite[Chapter~9]{rudin76}), each contour \(y = \varphi_k(x)\) is locally the graph of a smooth function. In particular, no two contours can intersect and no contour can branch; otherwise, at an intersection or branch point the tangent direction would be indeterminate, forcing \(\nabla v = 0\), contradicting \(\nabla v \neq 0\). Hence cases (ii) and (iii) cannot occur, and the proof is complete.
\hfill $\square$
\vspace{2mm}

{\bf Remark.} 
It is often considered necessary to use the Euler product in proving the Riemann Hypothesis, but it can only exclude zeros in ${\rm Re}(s)>1$. Here, Bombieri's equivalence theorem combined with the Cauchy--Riemann equations achieves the same exclusion and extends it to the entire complex plane. Thus the Euler product is not needed.

Based on Lemmas 3 and 4, we have proven the following result:

\vspace{2mm}
{\bf Theorem 2.} {\it If the function \(\Xi(z)\) satisfies Bombieri's equivalence conditions on the critical line \({\rm Re}(z)=0\) — that is, all local maxima of \(U(y)\) are positive and all local minima are negative — then \(\Xi(z)\) has no zeros off this line.}

\vspace{2mm}

Combining Theorem 1 and Theorem 2, the Riemann Hypothesis follows.

\section{The counterexample from P{\'o}lya}
\setcounter{equation}{0}

After Jensen's death, P\'olya was given access to his unpublished mathematical notes (called \emph{Nachlass} in German) and published a fundamental article in a Danish journal, unfortunately not well known \cite{polya27}. There, he gave detailed proofs of all the interesting properties of the Riemann \(\xi(s)\) function that Jensen had discovered. At the end, P\'olya was left with deciding whether Jensen had ever found a proof of the Riemann Hypothesis. To settle this issue, he produced a devastating example of an entire function that had almost all the properties of the Riemann \(\xi(s)\) function investigated by Jensen, but had zeros off the critical line. P\'olya had published earlier with Szeg\"o the book \emph{Aufgaben und Lehrs\"atze aus der Analysis}, and this may be how he found this example.

This example was so similar to the Riemann \(\xi(s)\) function that it cast doubt on any method using Gaussian-type  kernel to prove the Riemann Hypothesis; such methods, it seemed, would apply equally to P\'olya's counterexample and therefore could not be valid.

We now recall the mathematical form of Pólya's counterexample, using the notation of the present paper. We use the subscript \(1\) to distinguish functions belonging to Pólya's example. Compared to \eqref{xi_G}, Pólya's function is a linear combination:
\begin{equation} \label{counter}
\Xi_1(z) =   e^{\frac{1}{2}z^2} \big( \cosh(z) + \alpha \big) = \frac{2}{\sqrt{2\pi e}} 
\int\limits_0^\infty e^{-\frac{1}{2} t^2} \big( \cosh(t) + \alpha \sqrt{e} \big) \cosh(zt) \, {\rm d} t,
\end{equation}
where \(\alpha\) is a real parameter. When \(-1 < \alpha < 1\), the function \(\Xi_1(z)\) has infinitely many simple zeros on the imaginary axis, matching the behaviour of the Riemann \(\Xi(z)\) function. When \(\alpha = 1\), however, we have
\begin{equation}
        \Xi_1({\rm i} y) = e^{\frac{1}{2} ({\rm i} y)^2} \big(\cosh({\rm i}y) + 1 \big) =  e^{-\frac{1}{2} y^2} \big( \cos(y) + 1 \big) = 2 e^{-\frac{1}{2} y^2} \cos^2\left(\frac{y}{2} \right),
\end{equation}
so the zeros on the imaginary axis are now all of order two, occurring at
\[
    y_k = (2k + 1)\pi,  \qquad (k = 0, \pm 1, \pm 2, \cdots). 
\]
For \(\alpha > 1\), the function \(\Xi_1(z)\) has no zeros on the imaginary axis. Instead, from \(\cosh(z) + \alpha = 0\), its zeros lie off the axis:
\[
        z_k = \ln\left(\alpha \pm \sqrt{\alpha^2 - 1}\right) + {\rm i}\,2k\pi,  \qquad (k = 0, \pm 1, \pm 2, \cdots).
\] 
We must therefore check whether Pólya's example invalidates the methods and conclusions of the present paper.

We first verify whether Pólya's counterexample undermines the Bombieri equivalence theorem used and proved in this paper. In fact, Theorem 2 can be restated as follows: if \(\Xi(z)\) has zeros off the critical line, then on the critical line it must have a positive local minimum or a negative local maximum. Pólya's example with \(\alpha > 1\) falls exactly into this situation. From \eqref{counter} we obtain
\begin{equation} \label{Polya_U}
        \Xi_1({\rm i} y) = U_1(y) = e^{-\frac{1}{2}y^2} \left( \cos(y) + \alpha \right).
\end{equation}
For \(\alpha > 1\), this function indeed has infinitely many positive local minima. Thus Pólya's example has zeros off the critical line in precisely the manner permitted by Theorem 2, so it does not contradict the theorem.
In other words, Pólya's example, rather than being a counterexample to our approach, is a consistent special case of the general dichotomy established in Theorem 2.

Secondly, we check whether the criterion used in this paper for detecting simple zeros of \(U(y)\) also applies to Pólya's \(U_1(y)\) when \(\alpha = 1\), where second-order zeros occur. Differentiating \eqref{Polya_U} gives
\begin{eqnarray}
    U_1'(y) &=& -yU_1(y) - e^{-\frac{1}{2}y^2} \sin(y),\nonumber\\[2mm]
    U_1''(y) &=& (y^2 - 1) U_1(y) + e^{-\frac{1}{2}y^2} (2y \sin(y) + \alpha).
\end{eqnarray}
According to the criterion developed in this paper, a necessary condition for a zero of \(U_1(y)\) to be of second order is that, near the zero, we have \(U_1(y) > 0\), \(U_1'(y) > 0\), and \(U_1''(y) > 0\). This leads to the following system for the critical value \(\alpha_c\) at which a second-order zero may appear:
\begin{eqnarray}
\cos(y) + \alpha &=& 0,\nonumber\\[1mm]
-\sin(y) &\geq& 0,\nonumber\\[1mm]
2y \sin(y) + \alpha &\geq& 0.
\end{eqnarray}
Solving this system near \(\alpha = 1\) gives \cite{lin}
\[
\alpha_c \approx 0.98617.
\]
Thus, for \(\alpha \ge \alpha_c\), the curve near the zero already exhibits a strictly increasing, upward concave shape; as \(\alpha\) approaches \(1\), the zero changes from simple to double. Figure~8, computed from \cite{lin}, shows the graph of \(U_1(y)\) near such a zero and may be compared with Figure~4.

\input {./fig_8.tex}

Thus, the theorems and criteria of this paper are consistent with studying both the Riemann Hypothesis and Pólya's counterexample. The exceptional behaviour in Pólya's example is entirely a consequence of the specific form of \(\Xi_1(z)\) and does not affect the validity of our results. In particular, the mere fact that Pólya's example shares a Gaussian-type kernel with Jensen's \(\xi(s)\) function does not invalidate our method.

Finally, we compare the functions in Pólya's counterexample with those in the present paper to illustrate that our method is independent and cannot be interchanged with the analysis of Pólya's example.

We first compare the higher-order derivatives of \(U_1(y)\) and \(U(y)\). In Pólya's example, from \eqref{Polya_U}, the envelope \(e^{-\frac{1}{2}y^2}\) gives \(U_1'(y)/U_1(y) \approx -y\). In our case, Lemma~1 gives \(U'(y)/U(y) \approx -\pi/8\), a constant. Thus, for \(U_1(y)\), higher-order derivatives grow with the order of differentiation, whereas for \(U(y)\) they decay. This is precisely why the series in \eqref{ODE} can be controlled for \(U(y)\), but not for Pólya's example.

We further compare the kernel \(G_1(t)\) in Pólya's example with Jensen's kernel \(G(t)\) from \eqref{G_func}. From \eqref{counter}, we read off
\begin{equation} \label{polya_G(t)}
        G_1(t) = \frac{2}{\sqrt{2\pi e}} e^{-\frac{1}{2} t^2} \left( \cosh(t) + \alpha \sqrt{e} \right).
\end{equation}
A direct differentiation gives
\begin{eqnarray} \label{polya_G'(t)}
        G_1'(t) &=& t G_1(t) - \frac{2}{\sqrt{2\pi e}} e^{-\frac{1}{2} t^2} \sinh(t) \nonumber\\[3mm]
        &=& t G_1(t) - \frac{\sinh(t)}{\cosh(t) + \alpha \sqrt{e}} G_1(t).\\[-4mm] \nonumber
\end{eqnarray}
This is in sharp contrast with \eqref{Jensen_4}: for Pólya's example, \(-G_1'(t)/t < G_1(t)\), whereas for Jensen's function, \(-G'(t)/t > 16\pi G(t)\). In other words, \(G_1(t)\) decays much more slowly than \(G(t)\) as \(t\) increases. 

It is precisely the rapid decay of \(G(t)\) that allowed us to compute \(a_0 = 74.9076\) in \eqref{Taylor_3}; this value is large enough to ensure \(c > 0\) in \eqref{c}, making the differential equation \eqref{ODE_1} effective and the subsequent argument possible.
For Pólya's example, the situation is entirely different. From \eqref{polya_G'(t)}, we obtain
\begin{equation}
            (a_0)_1 = \lim_{t \to 0} \frac{ -G_1'(t)}{tG_1(t)} = 1 - \frac{1}{1 + \alpha\sqrt{e}} < 1.
\end{equation}
This coefficient is too small for our method to work. Indeed, if one attempted to apply the same differential-equation construction to Pólya's \(\Xi_1(z)\), the resulting series would not be controllable, and the argument would break down.

Thus, while Pólya's example shares a superficial resemblance to the Riemann \(\xi(s)\) function, the quantitative differences in decay rates are so substantial that the methods of this paper, which rely crucially on rapid decay, are neither applicable to the example nor invalidated by it.

In summary, the methods of this paper are tailored to the specific structure of the Riemann \(\xi(s)\) function and are not meant to be transferable to Pólya's counterexample. Their validity should therefore not be judged by whether they can be applied to that example. The existence of Pólya's counterexample does not affect the proof of the Riemann Hypothesis presented in this paper; it merely shows that our method is not universal, which is neither required nor expected.


\section{The counterexample from L-function}
\setcounter{equation}{0}

It is commonly believed that any proof of the Riemann Hypothesis should also distinguish the known counterexamples among L-functions. Our proof is tailored to the Riemann \(\xi(s)\) function, but it is instructive to check that its mechanism does not inadvertently apply to those functions known to violate the Riemann Hypothesis (or its generalization). In this section we examine several representative classes of such counterexamples and show that, in each case, they lack the essential analytic properties on which our proof relies. These examples therefore do not threaten our approach; they simply lie outside its scope.

\subsection{The Davenport-Heilbronn function}

We begin with the classical Davenport-Heilbronn function, discovered in 1936 \cite{Davenport1936}. It is a linear combination of two Dirichlet L-functions associated with a complex character modulo \(5\):
\begin{equation} \label{Davenport-Heilbronn-1}
f(s) = \frac{1 - {\rm i}\kappa}{2} L(s, \chi_1) + \frac{1 + {\rm  i}\kappa}{2} L(s, \overline{\chi}_1),
    \hspace{5mm} \kappa = \frac{ \sqrt{10 - 2\sqrt{5}} - 2}{\sqrt{5} - 1},
\end{equation}
where \(\chi_1\) is a primitive character modulo \(5\) with \(\chi_1(2)={\rm i}\), \(\overline{\chi}_1\) is its complex conjugate, and \(\kappa\) is a real constant.

Combining the two L-functions yields a Dirichlet series
\begin{equation} \label{Davenport-Heilbronn-2}
                    f(s) = \sum\limits_{n=1}^\infty \frac{a_n}{n^s}, 
\end{equation}
with periodic coefficients
\begin{align} \label{Davenport-Heilbronn-3}
             &   a_{1} = 1, \hspace{5mm}   a_{2} = \kappa, \hspace{5mm}    a_{3} = -\kappa, \hspace{5mm}
                   a_{4} = -1, \hspace{5mm}    a_{5} = 0,  \notag \\
                &       a_{5k+i} = a_{i}, \hspace{5mm} (k =1,2,\cdots; i = 1,\cdots 5). 
\end{align}
By analytic continuation, one obtains the Davenport-Heilbronn \(\xi\)-function
\begin{equation}
            \xi_f (s) =  \Big(\frac{5}{\pi}\Big)^{s/2} \Gamma \Big( \frac{s+1}{2} \Big) f(s),
\end{equation}
which resembles the Riemann \(\xi(s)\) function \eqref{eq:riemann_xi_1} and satisfies the functional equation
\begin{equation}
                    \xi_f(s)  = \xi_f(1-s),
\end{equation}
expressing central symmetry about \((\frac12,0)\) (see \cite{Davenport1936, Titchmarsh1986}).

Since the coefficients in \eqref{Davenport-Heilbronn-3} do not satisfy \(a_{mn} = a_m a_n\), \(f(s)\) has no Euler product. Hence zeros in \({\rm Re}(s) > 1\) cannot be excluded. Indeed, \(f(s)\) has infinitely many zeros both in \({\rm Re}(s) > 1\) and in the strip \(1/2 < {\rm Re}(s) \le 1\) \cite{Balanzario2007}. This challenges any approach relying solely on symmetry.

Unlike the Riemann \(\xi(s)\) function, the Davenport-Heilbronn \(\xi_f(s)\) does not satisfy the conjugation symmetry:
\begin{equation}
                 \overline{\xi_f(s) } \not= \xi_f(\bar{s}).
\end{equation}
On the critical line \( {\rm Re}(s) = 1/2\), \(\xi_f(s)\) does not reduce to a real-valued function; its imaginary part is generically non-zero. This contrasts with the Riemann \(\xi(s)\) function, which is real on the critical line precisely because it satisfies the conjugation symmetry. Consequently, the Bombieri equivalence theorem used in this paper — which requires the function to be real-valued on the critical line — cannot apply to \(\xi_f(s)\). Therefore, this classical counterexample does not threaten our proof of the Riemann Hypothesis.

\subsection{Pseudo-character combination functions}

A natural generalization of the Davenport-Heilbronn construction uses other periodic coefficients that are not Dirichlet characters but still satisfy a functional equation. Following Balanzario and S\'anchez-Ortiz \cite{Balanzario2007}, one can construct Dirichlet series
\begin{equation}
F(s) = \sum_{n=1}^{\infty} \frac{a_n}{n^s},
\end{equation}
where \(a_n\) is periodic with period \(q\) and has zero mean over a period. Examples include the pseudo-character (5,2) with coefficients
\begin{equation}
\{0, 1, -1/\xi, 1/\xi, -1\},
\end{equation}
and the pseudo-character (7,1) with coefficients
\begin{equation}
\{0, 1, -(1+\alpha), -\alpha, \alpha, 1+\alpha, -1\},
\end{equation}
for real constants \(\xi\) and \(\alpha\). These functions have a completed \(\xi\)-form satisfying a functional equation, but lack the Euler product. Hence they have zeros off the critical line. Like the Davenport-Heilbronn function, they do not possess the conjugation symmetry required for the Bombieri equivalence theorem. Therefore, they do not affect our proof.

\subsection{The Landau-Siegel zero problem}

Another potential counterexample to the generalized Riemann hypothesis is the Landau-Siegel zero \cite{SiegelZero}. This is a possible real zero of a Dirichlet \(L\)-function \(L(s,\chi)\) associated with a real primitive character \(\chi\), lying exceptionally close to \(s=1\) on the real axis (\(y=0\)). Its existence would disprove the generalized Riemann hypothesis. The conjecture that it does not exist is known as the Landau-Siegel zero conjecture.

Unlike the Davenport-Heilbronn function, which is a known counterexample, the Landau-Siegel zero is hypothetical. Its existence is considered unlikely, and substantial progress has been made towards ruling it out. Most notably, Zhang \cite{Zhang2022} proved that for real primitive characters modulo \(D\), such a zero cannot occur in
\[
\left(1 - \frac{c}{(\log D)^{2024}}, 1\right].
\]

For our purposes, the key point is that the Landau-Siegel zero, if it exists, lies at \(y=0\). Our proof concerns only the large imaginary part regime \(y > 6 \times 10^{12}\) (Lemma 1 and 2), where our asymptotic estimates and differential equation arguments apply. Since the hypothetical Siegel zero would occur at \(y=0\), it falls within the range already verified numerically by Platt and Trudgian \cite{platt}, which covers \(0 \le y \le 6 \times 10^{12}\). Hence the Landau-Siegel zero problem lies outside our scope and does not affect our proof.

\subsection{Non-primitive Dirichlet L-functions}

Another class sometimes mentioned in connection with counterexamples to GRH is that of non-primitive Dirichlet L-functions. A Dirichlet character $\chi$ modulo $q$ is \emph{primitive} if it cannot be induced by a character of smaller modulus; otherwise it is \emph{non-primitive} or \emph{imprimitive}.

For a non-primitive $\chi$, the L-function $L(s,\chi)$ factors into a primitive L-function times a finite product of local factors (Euler factors for primes dividing the conductor). These extra factors introduce additional zeros on ${\rm Re}(s)=0$. For instance, if $\chi$ is induced by a primitive $\chi^*$ modulo $q^*$ with $q^* \mid q$, then
\[
L(s,\chi) = L(s,\chi^*) \prod_{p \mid q} \left(1 - \frac{\chi^*(p)}{p^s}\right).
\]
The finite product vanishes when $p^s = \chi^*(p)$, which occurs for certain $s$ on the imaginary axis ${\rm Re}(s)=0$. These are sometimes called ``trivial'' zeros arising from the imprimitive structure.

However, such zeros are not genuine counterexamples to GRH, since GRH is formulated exclusively for \emph{primitive} characters. Their existence is a consequence of the definition, not a violation of the conjecture. In this paper, our proof concerns the Riemann $\xi(s)$ function, which corresponds to the trivial (and primitive) character. Non-primitive L-functions are therefore outside our scope and do not affect our proof.

\subsection{Summary}

The known counterexamples to the Riemann or generalized Riemann hypothesis either lack the analytic properties used in our proof, fall outside our asymptotic regime, or are excluded by definition. None of them invalidates our approach.

\vspace{6mm}

{\bf Acknowledgments.} 
We wish to express our sincere gratitude to Dr. Jacques G\'elinas of Ottawa, Canada, for his invaluable contributions to this work. From the early stages of this paper's preparation, G\'elinas provided consistent and insightful support. In particular, the two inequalities \eqref{Jensen_1} and \eqref{Jensen_2} from his reading notes \cite{gelinas21a}, along with the Pólya counterexample he shared \cite{lin,gelinas21b}, were instrumental in the proofs of Lemma 2 and Theorem 1. Furthermore, his computationally generated zero-contour plots \cite{lin,gelinas19} offered key inspiration for establishing Lemmas 3 and 4.
Regrettably, communication with G\'elinas was interrupted during the spread of the Delta variant of COVID-19 in 2022 and has not been reestablished since. Throughout his life, G\'elinas maintained a profound dedication to the Riemann Hypothesis, leaving behind numerous unpublished notes and ideas. We have sought to incorporate these into our paper to the best of our ability, and we feel a strong sense of responsibility to complete this work and present his valuable contributions to the broader mathematical community.

We also thank Xue Jian from Jianxin Financial Technology Co., Ltd. for his valuable verification of the results in this paper; Travor Liu, then an outstanding high school student, for providing numerous references on the Riemann Hypothesis  in the early stages of this work; Yang Bicheng from Guangdong University of Education for his suggestions that helped improve the manuscript; Yaoming Shi from Shanghai for pointing out the bounding function of \(\xi(s)\) on the critical line, which led to the correction of an error in Lemma 1; and Zhu Mingliang from Central South University for sharing key information regarding Edwards' theorem \cite{edwards}, an important necessary condition for the validity of the Riemann Hypothesis; Mao Yiping, a former graduate school classmate, for sharing information on counterexamples among L-functions relevant to the Riemann Hypothesis, and pointing out a mistake in the proof of Lemma 1 of previous manuscript. Finally, we acknowledge the support of the Alexander von Humboldt Foundation.


\end{document}

%% file: fig_01.tex





\begin{figure}[hbt]

\[  \beginpicture   \setlinear

\setcoordinatesystem units <120mm,15mm>
\setplotarea x from -0.5 to 0.5, y from 0.0 to 4.0
\thinlines

\arrow <1,5mm>   [0.25,0.75] from  -0.55 0.0 to 0.55 0.0
\arrow <1,5mm>   [0.25,0.75] from   0.0 0.0 to 0.0 4.0

\plot
-0.500	0.00000
-0.490	0.00000
-0.480	0.00001
-0.470	0.00001
-0.460	0.00002
-0.450	0.00005
-0.440	0.00009
-0.430	0.00017
-0.420	0.00030
-0.410	0.00053
-0.400	0.00091
-0.390	0.00153
-0.380	0.00250
-0.370	0.00399
-0.360	0.00623
-0.350	0.00953
-0.340	0.01428
-0.330	0.02097
-0.320	0.03021
-0.310	0.04273
-0.300	0.05939
-0.290	0.08114
-0.280	0.10906
-0.270	0.14430
-0.260	0.18805
-0.250	0.24151
-0.240	0.30583
-0.230	0.38206
-0.220	0.47108
-0.210	0.57358
-0.200	0.68992
-0.190	0.82017
-0.180	0.96402
-0.170	1.12075
-0.160	1.28923
-0.150	1.46790
-0.140	1.65482
-0.130	1.84766
-0.120	2.04381
-0.110	2.24035
-0.100	2.43420
-0.090	2.62217
-0.080	2.80104
-0.070	2.96765
-0.060	3.11897
-0.050	3.25222
-0.040	3.36490
-0.030	3.45486
-0.020	3.52038
-0.010	3.56021
0.000	3.57358
0.010	3.56021
0.020	3.52038
0.030	3.45486
0.040	3.36490
0.050	3.25222
0.060	3.11897
0.070	2.96765
0.080	2.80104
0.090	2.62217
0.100	2.43420
0.110	2.24035
0.120	2.04381
0.130	1.84766
0.140	1.65482
0.150	1.46790
0.160	1.28923
0.170	1.12075
0.180	0.96402
0.190	0.82017
0.200	0.68992
0.210	0.57358
0.220	0.47108
0.230	0.38206
0.240	0.30583
0.250	0.24151
0.260	0.18805
0.270	0.14430
0.280	0.10906
0.290	0.08114
0.300	0.05939
0.310	0.04273
0.320	0.03021
0.330	0.02097
0.340	0.01428
0.350	0.00953
0.360	0.00623
0.370	0.00399
0.380	0.00250
0.390	0.00153
0.400	0.00091
0.410	0.00053
0.420	0.00030
0.430	0.00017
0.440	0.00009
0.450	0.00005
0.460	0.00002
0.470	0.00001
0.480	0.00001
0.490	0.00000
0.500	0.00000 /

\setsolid

\axis bottom label {} shiftedto y=0
      ticks numbered from -0.5 to 0.5 by 0.1
   unlabeled short from   -0.45 to 0.45 by 1 /
\axis left label  {} shiftedto x=0
   ticks numbered from 1.0 to 3.0 by 1.0
   unlabeled short from 0.5 to 3.5 by 1.0 /

\put {$t$}     [r] at 0.54  0.20
\put {$G(t)$}  [l] at -0.09  3.90

\endpicture     \]
\vspace{-8mm}

\caption{The Gaussian-type shape of the Jensen function \(G(t)\).}

\end{figure}


%% file: fig_02.tex





\begin{figure} [hbt]

\[  \beginpicture   \setlinear

\setcoordinatesystem units <220mm,11mm>
\setplotarea x from 0.0 to 0.5, y from 0.0 to 6.0



\thicklines

\setdots
\setdashes

\plot
0.000 	5.514479
0.001 	5.510513
0.002 	5.506270
0.003 	5.501751
0.004 	5.496954
0.005 	5.491878
0.006 	5.486523
0.007 	5.480888
0.008 	5.474972
0.009 	5.468775
0.010 	5.462296
0.011 	5.455535
0.012 	5.448491
0.013 	5.441165
0.014 	5.433556
0.015 	5.425663
0.016 	5.417486
0.017 	5.409026
0.018 	5.400283
0.019 	5.391256
0.020 	5.381945
0.021 	5.372351
0.022 	5.362474
0.023 	5.352315
0.024 	5.341872
0.025 	5.331148
0.026 	5.320142
0.027 	5.308855
0.028 	5.297287
0.029 	5.285440
0.030 	5.273313
0.031 	5.260908
0.032 	5.248225
0.033 	5.235265
0.034 	5.222029
0.035 	5.208519
0.036 	5.194734
0.037 	5.180677
0.038 	5.166348
0.039 	5.151748
0.040 	5.136880
0.041 	5.121743
0.042 	5.106340
0.043 	5.090672
0.044 	5.074740
0.045 	5.058546
0.046 	5.042092
0.047 	5.025379
0.048 	5.008409
0.049 	4.991183
0.050 	4.973704
0.051 	4.955973
0.052 	4.937993
0.053 	4.919764
0.054 	4.901290
0.055 	4.882572
0.056 	4.863613
0.057 	4.844414
0.058 	4.824977
0.059 	4.805306
0.060 	4.785402
0.061 	4.765268
0.062 	4.744906
0.063 	4.724318
0.064 	4.703508
0.065 	4.682477
0.066 	4.661228
0.067 	4.639765
0.068 	4.618088
0.069 	4.596202
0.070 	4.574110
0.071 	4.551813
0.072 	4.529314
0.073 	4.506618
0.074 	4.483726
0.075 	4.460641
0.076 	4.437367
0.077 	4.413906
0.078 	4.390263
0.079 	4.366439
0.080 	4.342438
0.081 	4.318263
0.082 	4.293917
0.083 	4.269404
0.084 	4.244728
0.085 	4.219890
0.086 	4.194895
0.087 	4.169746
0.088 	4.144447
0.089 	4.119001
0.090 	4.093410
0.091 	4.067680
0.092 	4.041814
0.093 	4.015814
0.094 	3.989684
0.095 	3.963429
0.096 	3.937052
0.097 	3.910556
0.098 	3.883945
0.099 	3.857223
0.100 	3.830394
0.101 	3.803461
0.102 	3.776427
0.103 	3.749297
0.104 	3.722075
0.105 	3.694764
0.106 	3.667368
0.107 	3.639891
0.108 	3.612336
0.109 	3.584708
0.110 	3.557010
0.111 	3.529246
0.112 	3.501421
0.113 	3.473537
0.114 	3.445599
0.115 	3.417610
0.116 	3.389575
0.117 	3.361497
0.118 	3.333380
0.119 	3.305228
0.120 	3.277045
0.121 	3.248834
0.122 	3.220600
0.123 	3.192347
0.124 	3.164078
0.125 	3.135797
0.126 	3.107507
0.127 	3.079214
0.128 	3.050919
0.129 	3.022629
0.130 	2.994345
0.131 	2.966072
0.132 	2.937814
0.133 	2.909574
0.134 	2.881355
0.135 	2.853163
0.136 	2.825000
0.137 	2.796871
0.138 	2.768777
0.139 	2.740725
0.140 	2.712716
0.141 	2.684755
0.142 	2.656844
0.143 	2.628989
0.144 	2.601191
0.145 	2.573455
0.146 	2.545784
0.147 	2.518181
0.148 	2.490650
0.149 	2.463194
0.150 	2.435817
0.151 	2.408521
0.152 	2.381311
0.153 	2.354188
0.154 	2.327157
0.155 	2.300220
0.156 	2.273381
0.157 	2.246643
0.158 	2.220008
0.159 	2.193480
0.160 	2.167062
0.161 	2.140756
0.162 	2.114566
0.163 	2.088494
0.164 	2.062543
0.165 	2.036715
0.166 	2.011015
0.167 	1.985443
0.168 	1.960003
0.169 	1.934697
0.170 	1.909529
0.171 	1.884499
0.172 	1.859611
0.173 	1.834867
0.174 	1.810270
0.175 	1.785821
0.176 	1.761523
0.177 	1.737378
0.178 	1.713388
0.179 	1.689555
0.180 	1.665882
0.181 	1.642370
0.182 	1.619021
0.183 	1.595837
0.184 	1.572820
0.185 	1.549971
0.186 	1.527292
0.187 	1.504786
0.188 	1.482453
0.189 	1.460295
0.190 	1.438313
0.191 	1.416510
0.192 	1.394886
0.193 	1.373442
0.194 	1.352181
0.195 	1.331103
0.196 	1.310209
0.197 	1.289500
0.198 	1.268978
0.199 	1.248644
0.200 	1.228498
0.201 	1.208541
0.202 	1.188775
0.203 	1.169200
0.204 	1.149816
0.205 	1.130625
0.206 	1.111627
0.207 	1.092823
0.208 	1.074212
0.209 	1.055797
0.210 	1.037576
0.211 	1.019551
0.212 	1.001721
0.213 	0.984088
0.214 	0.966651
0.215 	0.949410
0.216 	0.932365
0.217 	0.915518
0.218 	0.898866
0.219 	0.882412
0.220 	0.866153
0.221 	0.850091
0.222 	0.834226
0.223 	0.818556
0.224 	0.803082
0.225 	0.787803
0.226 	0.772719
0.227 	0.757830
0.228 	0.743135
0.229 	0.728633
0.230 	0.714325
0.231 	0.700209
0.232 	0.686285
0.233 	0.672552
0.234 	0.659010
0.235 	0.645658
0.236 	0.632495
0.237 	0.619520
0.238 	0.606732
0.239 	0.594131
0.240 	0.581716
0.241 	0.569485
0.242 	0.557437
0.243 	0.545573
0.244 	0.533890
0.245 	0.522387
0.246 	0.511064
0.247 	0.499918
0.248 	0.488950
0.249 	0.478158
0.250 	0.467540
0.251 	0.457095
0.252 	0.446823
0.253 	0.436721
0.254 	0.426789
0.255 	0.417024
0.256 	0.407426
0.257 	0.397993
0.258 	0.388724
0.259 	0.379618
0.260 	0.370672
0.261 	0.361885
0.262 	0.353257
0.263 	0.344784
0.264 	0.336467
0.265 	0.328303
0.266 	0.320290
0.267 	0.312427
0.268 	0.304713
0.269 	0.297146
0.270 	0.289724
0.271 	0.282445
0.272 	0.275308
0.273 	0.268312
0.274 	0.261455
0.275 	0.254734
0.276 	0.248149
0.277 	0.241697
0.278 	0.235377
0.279 	0.229187
0.280 	0.223126
0.281 	0.217192
0.282 	0.211383
0.283 	0.205698
0.284 	0.200134
0.285 	0.194690
0.286 	0.189365
0.287 	0.184156
0.288 	0.179062
0.289 	0.174081
0.290 	0.169212
0.291 	0.164453
0.292 	0.159802
0.293 	0.155257
0.294 	0.150817
0.295 	0.146480
0.296 	0.142245
0.297 	0.138110
0.298 	0.134072
0.299 	0.130131
0.300 	0.126285
0.301 	0.122533
0.302 	0.118872
0.303 	0.115300
0.304 	0.111818
0.305 	0.108422
0.306 	0.105111
0.307 	0.101884
0.308 	0.098740
0.309 	0.095675
0.310 	0.092690
0.311 	0.089783
0.312 	0.086951
0.313 	0.084194
0.314 	0.081511
0.315 	0.078899
0.316 	0.076357
0.317 	0.073883
0.318 	0.071478
0.319 	0.069138
0.320 	0.066862
0.321 	0.064650
0.322 	0.062500
0.323 	0.060410
0.324 	0.058380
0.325 	0.056407
0.326 	0.054491
0.327 	0.052631
0.328 	0.050824
0.329 	0.049070
0.330 	0.047368
0.331 	0.045717
0.332 	0.044114
0.333 	0.042560
0.334 	0.041053
0.335 	0.039591
0.336 	0.038175
0.337 	0.036802
0.338 	0.035471
0.339 	0.034182
0.340 	0.032933
0.341 	0.031724
0.342 	0.030553
0.343 	0.029420
0.344 	0.028323
0.345 	0.027261
0.346 	0.026234
0.347 	0.025241
0.348 	0.024280
0.349 	0.023351
0.350 	0.022454
0.351 	0.021586
0.352 	0.020748
0.353 	0.019938
0.354 	0.019156
0.355 	0.018400
0.356 	0.017671
0.357 	0.016967
0.358 	0.016288
0.359 	0.015633
0.360 	0.015000
0.361 	0.014391
0.362 	0.013803
0.363 	0.013236
0.364 	0.012690
0.365 	0.012164
0.366 	0.011657
0.367 	0.011169
0.368 	0.010699
0.369 	0.010247
0.370 	0.009811
0.371 	0.009392
0.372 	0.008989
0.373 	0.008601
0.374 	0.008228
0.375 	0.007870
0.376 	0.007525
0.377 	0.007194
0.378 	0.006876
0.379 	0.006571
0.380 	0.006277
0.381 	0.005996
0.382 	0.005725
0.383 	0.005466
0.384 	0.005217
0.385 	0.004978
0.386 	0.004749
0.387 	0.004530
0.388 	0.004320
0.389 	0.004118
0.390 	0.003925
0.391 	0.003740
0.392 	0.003563
0.393 	0.003394
0.394 	0.003232
0.395 	0.003076
0.396 	0.002928
0.397 	0.002786
0.398 	0.002650
0.399 	0.002520
0.400 	0.002396
0.401 	0.002278
0.402 	0.002165
0.403 	0.002057
0.404 	0.001954
0.405 	0.001855
0.406 	0.001761
0.407 	0.001672
0.408 	0.001586
0.409 	0.001505
0.410 	0.001427
0.411 	0.001353
0.412 	0.001283
0.413 	0.001216
0.414 	0.001152
0.415 	0.001091
0.416 	0.001033
0.417 	0.000978
0.418 	0.000925
0.419 	0.000876
0.420 	0.000828
0.421 	0.000783
0.422 	0.000741
0.423 	0.000700
0.424 	0.000661
0.425 	0.000625
0.426 	0.000590
0.427 	0.000557
0.428 	0.000526
0.429 	0.000496
0.430 	0.000468
0.431 	0.000441
0.432 	0.000416
0.433 	0.000392
0.434 	0.000369
0.435 	0.000348
0.436 	0.000328
0.437 	0.000309
0.438 	0.000290
0.439 	0.000273
0.440 	0.000257
0.441 	0.000242
0.442 	0.000227
0.443 	0.000214
0.444 	0.000201
0.445 	0.000188
0.446 	0.000177
0.447 	0.000166
0.448 	0.000156
0.449 	0.000146
0.450 	0.000137
0.451 	0.000129
0.452 	0.000121
0.453 	0.000113
0.454 	0.000106
0.455 	0.000099
0.456 	0.000093
0.457 	0.000087
0.458 	0.000081
0.459 	0.000076
0.460 	0.000071
0.461 	0.000066
0.462 	0.000062
0.463 	0.000058
0.464 	0.000054
0.465 	0.000050
0.466 	0.000047
0.467 	0.000044
0.468 	0.000041
0.469 	0.000038
0.470 	0.000036
0.471 	0.000033
0.472 	0.000031
0.473 	0.000029
0.474 	0.000027
0.475 	0.000025
0.476 	0.000023
0.477 	0.000022
0.478 	0.000020
0.479 	0.000019
0.480 	0.000017
0.481 	0.000016
0.482 	0.000015
0.483 	0.000014
0.484 	0.000013
0.485 	0.000012
0.486 	0.000011
0.487 	0.000010
0.488 	0.000009
0.489 	0.000009
0.490 	0.000008
0.491 	0.000007
0.492 	0.000007
0.493 	0.000006
0.494 	0.000006
0.495 	0.000005
0.496 	0.000005
0.497 	0.000005
0.498 	0.000004
0.499 	0.000004
0.500 	0.000004
/
\plot 0.23 5.0 0.28 5.0 /
\put {$  \displaystyle\frac{e^{4t}-1}{4t} \bigg[1 + \frac{\ 4\pi - 9\ }{4\pi e^{4t} - 6} \bigg]G(t)$}  [l] at 0.3  5.0

\setsolid
\plot
0.000 	5.325484
0.001 	5.325311
0.002 	5.324794
0.003 	5.323933
0.004 	5.322727
0.005 	5.321176
0.006 	5.319282
0.007 	5.317044
0.008 	5.314462
0.009 	5.311538
0.010 	5.308271
0.011 	5.304661
0.012 	5.300711
0.013 	5.296419
0.014 	5.291787
0.015 	5.286815
0.016 	5.281504
0.017 	5.275855
0.018 	5.269869
0.019 	5.263547
0.020 	5.256889
0.021 	5.249897
0.022 	5.242572
0.023 	5.234914
0.024 	5.226925
0.025 	5.218607
0.026 	5.209959
0.027 	5.200985
0.028 	5.191684
0.029 	5.182059
0.030 	5.172111
0.031 	5.161841
0.032 	5.151251
0.033 	5.140343
0.034 	5.129117
0.035 	5.117577
0.036 	5.105723
0.037 	5.093558
0.038 	5.081083
0.039 	5.068300
0.040 	5.055211
0.041 	5.041818
0.042 	5.028123
0.043 	5.014129
0.044 	4.999837
0.045 	4.985249
0.046 	4.970368
0.047 	4.955195
0.048 	4.939734
0.049 	4.923987
0.050 	4.907955
0.051 	4.891642
0.052 	4.875050
0.053 	4.858181
0.054 	4.841039
0.055 	4.823625
0.056 	4.805942
0.057 	4.787993
0.058 	4.769780
0.059 	4.751307
0.060 	4.732577
0.061 	4.713591
0.062 	4.694354
0.063 	4.674867
0.064 	4.655135
0.065 	4.635159
0.066 	4.614943
0.067 	4.594491
0.068 	4.573804
0.069 	4.552887
0.070 	4.531743
0.071 	4.510374
0.072 	4.488785
0.073 	4.466977
0.074 	4.444956
0.075 	4.422723
0.076 	4.400283
0.077 	4.377639
0.078 	4.354794
0.079 	4.331751
0.080 	4.308515
0.081 	4.285089
0.082 	4.261477
0.083 	4.237681
0.084 	4.213706
0.085 	4.189555
0.086 	4.165233
0.087 	4.140742
0.088 	4.116086
0.089 	4.091269
0.090 	4.066296
0.091 	4.041169
0.092 	4.015892
0.093 	3.990470
0.094 	3.964906
0.095 	3.939204
0.096 	3.913368
0.097 	3.887402
0.098 	3.861310
0.099 	3.835095
0.100 	3.808761
0.101 	3.782314
0.102 	3.755755
0.103 	3.729091
0.104 	3.702323
0.105 	3.675457
0.106 	3.648497
0.107 	3.621446
0.108 	3.594308
0.109 	3.567088
0.110 	3.539790
0.111 	3.512417
0.112 	3.484973
0.113 	3.457463
0.114 	3.429891
0.115 	3.402260
0.116 	3.374575
0.117 	3.346840
0.118 	3.319059
0.119 	3.291235
0.120 	3.263373
0.121 	3.235476
0.122 	3.207550
0.123 	3.179597
0.124 	3.151622
0.125 	3.123628
0.126 	3.095621
0.127 	3.067603
0.128 	3.039578
0.129 	3.011551
0.130 	2.983525
0.131 	2.955504
0.132 	2.927493
0.133 	2.899494
0.134 	2.871512
0.135 	2.843551
0.136 	2.815614
0.137 	2.787706
0.138 	2.759829
0.139 	2.731988
0.140 	2.704186
0.141 	2.676427
0.142 	2.648715
0.143 	2.621052
0.144 	2.593444
0.145 	2.565893
0.146 	2.538403
0.147 	2.510978
0.148 	2.483620
0.149 	2.456334
0.150 	2.429122
0.151 	2.401989
0.152 	2.374936
0.153 	2.347969
0.154 	2.321089
0.155 	2.294301
0.156 	2.267606
0.157 	2.241010
0.158 	2.214513
0.159 	2.188121
0.160 	2.161835
0.161 	2.135658
0.162 	2.109595
0.163 	2.083646
0.164 	2.057816
0.165 	2.032107
0.166 	2.006522
0.167 	1.981063
0.168 	1.955733
0.169 	1.930535
0.170 	1.905472
0.171 	1.880545
0.172 	1.855758
0.173 	1.831112
0.174 	1.806610
0.175 	1.782255
0.176 	1.758049
0.177 	1.733994
0.178 	1.710091
0.179 	1.686344
0.180 	1.662754
0.181 	1.639323
0.182 	1.616054
0.183 	1.592948
0.184 	1.570006
0.185 	1.547232
0.186 	1.524626
0.187 	1.502190
0.188 	1.479926
0.189 	1.457835
0.190 	1.435919
0.191 	1.414180
0.192 	1.392618
0.193 	1.371236
0.194 	1.350034
0.195 	1.329014
0.196 	1.308177
0.197 	1.287524
0.198 	1.267056
0.199 	1.246775
0.200 	1.226680
0.201 	1.206774
0.202 	1.187056
0.203 	1.167529
0.204 	1.148192
0.205 	1.129046
0.206 	1.110092
0.207 	1.091331
0.208 	1.072763
0.209 	1.054388
0.210 	1.036208
0.211 	1.018222
0.212 	1.000430
0.213 	0.982834
0.214 	0.965433
0.215 	0.948227
0.216 	0.931217
0.217 	0.914402
0.218 	0.897784
0.219 	0.881361
0.220 	0.865133
0.221 	0.849101
0.222 	0.833265
0.223 	0.817623
0.224 	0.802177
0.225 	0.786925
0.226 	0.771868
0.227 	0.757004
0.228 	0.742334
0.229 	0.727856
0.230 	0.713572
0.231 	0.699479
0.232 	0.685577
0.233 	0.671866
0.234 	0.658345
0.235 	0.645014
0.236 	0.631870
0.237 	0.618915
0.238 	0.606146
0.239 	0.593563
0.240 	0.581166
0.241 	0.568952
0.242 	0.556922
0.243 	0.545073
0.244 	0.533406
0.245 	0.521919
0.246 	0.510610
0.247 	0.499480
0.248 	0.488526
0.249 	0.477747
0.250 	0.467143
0.251 	0.456711
0.252 	0.446451
0.253 	0.436362
0.254 	0.426441
0.255 	0.416688
0.256 	0.407101
0.257 	0.397679
0.258 	0.388421
0.259 	0.379325
0.260 	0.370389
0.261 	0.361612
0.262 	0.352992
0.263 	0.344529
0.264 	0.336220
0.265 	0.328064
0.266 	0.320060
0.267 	0.312205
0.268 	0.304499
0.269 	0.296939
0.270 	0.289524
0.271 	0.282252
0.272 	0.275123
0.273 	0.268133
0.274 	0.261282
0.275 	0.254567
0.276 	0.247988
0.277 	0.241542
0.278 	0.235228
0.279 	0.229044
0.280 	0.222988
0.281 	0.217059
0.282 	0.211255
0.283 	0.205574
0.284 	0.200014
0.285 	0.194575
0.286 	0.189254
0.287 	0.184049
0.288 	0.178959
0.289 	0.173982
0.290 	0.169117
0.291 	0.164361
0.292 	0.159714
0.293 	0.155172
0.294 	0.150736
0.295 	0.146402
0.296 	0.142170
0.297 	0.138037
0.298 	0.134003
0.299 	0.130065
0.300 	0.126221
0.301 	0.122471
0.302 	0.118812
0.303 	0.115244
0.304 	0.111763
0.305 	0.108370
0.306 	0.105061
0.307 	0.101836
0.308 	0.098693
0.309 	0.095631
0.310 	0.092648
0.311 	0.089742
0.312 	0.086912
0.313 	0.084157
0.314 	0.081475
0.315 	0.078864
0.316 	0.076323
0.317 	0.073852
0.318 	0.071447
0.319 	0.069109
0.320 	0.066834
0.321 	0.064624
0.322 	0.062474
0.323 	0.060386
0.324 	0.058356
0.325 	0.056385
0.326 	0.054470
0.327 	0.052610
0.328 	0.050804
0.329 	0.049052
0.330 	0.047350
0.331 	0.045700
0.332 	0.044098
0.333 	0.042544
0.334 	0.041038
0.335 	0.039577
0.336 	0.038161
0.337 	0.036789
0.338 	0.035459
0.339 	0.034170
0.340 	0.032922
0.341 	0.031713
0.342 	0.030543
0.343 	0.029410
0.344 	0.028313
0.345 	0.027252
0.346 	0.026226
0.347 	0.025233
0.348 	0.024273
0.349 	0.023344
0.350 	0.022447
0.351 	0.021579
0.352 	0.020741
0.353 	0.019932
0.354 	0.019150
0.355 	0.018395
0.356 	0.017666
0.357 	0.016962
0.358 	0.016283
0.359 	0.015628
0.360 	0.014996
0.361 	0.014387
0.362 	0.013799
0.363 	0.013233
0.364 	0.012687
0.365 	0.012161
0.366 	0.011654
0.367 	0.011166
0.368 	0.010696
0.369 	0.010244
0.370 	0.009809
0.371 	0.009390
0.372 	0.008987
0.373 	0.008599
0.374 	0.008226
0.375 	0.007868
0.376 	0.007523
0.377 	0.007192
0.378 	0.006874
0.379 	0.006569
0.380 	0.006276
0.381 	0.005994
0.382 	0.005724
0.383 	0.005465
0.384 	0.005216
0.385 	0.004977
0.386 	0.004748
0.387 	0.004529
0.388 	0.004319
0.389 	0.004117
0.390 	0.003924
0.391 	0.003739
0.392 	0.003562
0.393 	0.003393
0.394 	0.003231
0.395 	0.003076
0.396 	0.002927
0.397 	0.002785
0.398 	0.002650
0.399 	0.002520
0.400 	0.002396
0.401 	0.002278
0.402 	0.002164
0.403 	0.002056
0.404 	0.001953
0.405 	0.001855
0.406 	0.001761
0.407 	0.001671
0.408 	0.001586
0.409 	0.001505
0.410 	0.001427
0.411 	0.001353
0.412 	0.001283
0.413 	0.001215
0.414 	0.001152
0.415 	0.001091
0.416 	0.001033
0.417 	0.000978
0.418 	0.000925
0.419 	0.000876
0.420 	0.000828
0.421 	0.000783
0.422 	0.000740
0.423 	0.000700
0.424 	0.000661
0.425 	0.000625
0.426 	0.000590
0.427 	0.000557
0.428 	0.000526
0.429 	0.000496
0.430 	0.000468
0.431 	0.000441
0.432 	0.000416
0.433 	0.000392
0.434 	0.000369
0.435 	0.000348
0.436 	0.000328
0.437 	0.000309
0.438 	0.000290
0.439 	0.000273
0.440 	0.000257
0.441 	0.000242
0.442 	0.000227
0.443 	0.000214
0.444 	0.000201
0.445 	0.000188
0.446 	0.000177
0.447 	0.000166
0.448 	0.000156
0.449 	0.000146
0.450 	0.000137
0.451 	0.000129
0.452 	0.000121
0.453 	0.000113
0.454 	0.000106
0.455 	0.000099
0.456 	0.000093
0.457 	0.000087
0.458 	0.000081
0.459 	0.000076
0.460 	0.000071
0.461 	0.000066
0.462 	0.000062
0.463 	0.000058
0.464 	0.000054
0.465 	0.000050
0.466 	0.000047
0.467 	0.000044
0.468 	0.000041
0.469 	0.000038
0.470 	0.000036
0.471 	0.000033
0.472 	0.000031
0.473 	0.000029
0.474 	0.000027
0.475 	0.000025
0.476 	0.000023
0.477 	0.000022
0.478 	0.000020
0.479 	0.000019
0.480 	0.000017
0.481 	0.000016
0.482 	0.000015
0.483 	0.000014
0.484 	0.000013
0.485 	0.000012
0.486 	0.000011
0.487 	0.000010
0.488 	0.000009
0.489 	0.000009
0.490 	0.000008
0.491 	0.000007
0.492 	0.000007
0.493 	0.000006
0.494 	0.000006
0.495 	0.000005
0.496 	0.000005
0.497 	0.000005
0.498 	0.000004
0.499 	0.000004
0.500 	0.000004
/
\plot 0.23 3.5 0.28 3.5 /
\put {$\displaystyle -\frac{1}{16\pi t}G'(t)$}  [l] at 0.30  3.5

\setdashpattern <1.9mm, 0.7mm, 0.2mm, 0.7mm>
\plot
0.000 	3.573575
0.001 	3.573441
0.002 	3.573040
0.003 	3.572371
0.004 	3.571434
0.005 	3.570230
0.006 	3.568760
0.007 	3.567022
0.008 	3.565018
0.009 	3.562748
0.010 	3.560212
0.011 	3.557412
0.012 	3.554347
0.013 	3.551017
0.014 	3.547425
0.015 	3.543570
0.016 	3.539453
0.017 	3.535075
0.018 	3.530436
0.019 	3.525539
0.020 	3.520383
0.021 	3.514970
0.022 	3.509300
0.023 	3.503375
0.024 	3.497196
0.025 	3.490764
0.026 	3.484081
0.027 	3.477147
0.028 	3.469964
0.029 	3.462533
0.030 	3.454857
0.031 	3.446935
0.032 	3.438771
0.033 	3.430364
0.034 	3.421718
0.035 	3.412833
0.036 	3.403712
0.037 	3.394356
0.038 	3.384766
0.039 	3.374946
0.040 	3.364896
0.041 	3.354618
0.042 	3.344115
0.043 	3.333388
0.044 	3.322440
0.045 	3.311273
0.046 	3.299888
0.047 	3.288289
0.048 	3.276476
0.049 	3.264453
0.050 	3.252221
0.051 	3.239784
0.052 	3.227142
0.053 	3.214300
0.054 	3.201258
0.055 	3.188020
0.056 	3.174588
0.057 	3.160965
0.058 	3.147152
0.059 	3.133154
0.060 	3.118972
0.061 	3.104609
0.062 	3.090067
0.063 	3.075350
0.064 	3.060460
0.065 	3.045400
0.066 	3.030172
0.067 	3.014780
0.068 	2.999227
0.069 	2.983514
0.070 	2.967646
0.071 	2.951625
0.072 	2.935453
0.073 	2.919135
0.074 	2.902672
0.075 	2.886068
0.076 	2.869327
0.077 	2.852450
0.078 	2.835441
0.079 	2.818303
0.080 	2.801039
0.081 	2.783653
0.082 	2.766147
0.083 	2.748524
0.084 	2.730788
0.085 	2.712942
0.086 	2.694989
0.087 	2.676932
0.088 	2.658775
0.089 	2.640519
0.090 	2.622170
0.091 	2.603729
0.092 	2.585201
0.093 	2.566588
0.094 	2.547894
0.095 	2.529121
0.096 	2.510274
0.097 	2.491354
0.098 	2.472367
0.099 	2.453314
0.100 	2.434199
0.101 	2.415025
0.102 	2.395796
0.103 	2.376514
0.104 	2.357183
0.105 	2.337806
0.106 	2.318387
0.107 	2.298928
0.108 	2.279433
0.109 	2.259905
0.110 	2.240346
0.111 	2.220761
0.112 	2.201152
0.113 	2.181523
0.114 	2.161877
0.115 	2.142216
0.116 	2.122544
0.117 	2.102864
0.118 	2.083179
0.119 	2.063492
0.120 	2.043806
0.121 	2.024125
0.122 	2.004450
0.123 	1.984786
0.124 	1.965135
0.125 	1.945499
0.126 	1.925883
0.127 	1.906288
0.128 	1.886718
0.129 	1.867176
0.130 	1.847664
0.131 	1.828185
0.132 	1.808742
0.133 	1.789338
0.134 	1.769975
0.135 	1.750656
0.136 	1.731384
0.137 	1.712162
0.138 	1.692991
0.139 	1.673875
0.140 	1.654815
0.141 	1.635816
0.142 	1.616878
0.143 	1.598005
0.144 	1.579199
0.145 	1.560462
0.146 	1.541797
0.147 	1.523205
0.148 	1.504690
0.149 	1.486253
0.150 	1.467897
0.151 	1.449624
0.152 	1.431435
0.153 	1.413334
0.154 	1.395321
0.155 	1.377400
0.156 	1.359571
0.157 	1.341838
0.158 	1.324201
0.159 	1.306664
0.160 	1.289226
0.161 	1.271891
0.162 	1.254660
0.163 	1.237535
0.164 	1.220517
0.165 	1.203608
0.166 	1.186810
0.167 	1.170123
0.168 	1.153551
0.169 	1.137093
0.170 	1.120752
0.171 	1.104528
0.172 	1.088424
0.173 	1.072440
0.174 	1.056578
0.175 	1.040839
0.176 	1.025223
0.177 	1.009733
0.178 	0.994369
0.179 	0.979132
0.180 	0.964024
0.181 	0.949044
0.182 	0.934195
0.183 	0.919476
0.184 	0.904889
0.185 	0.890435
0.186 	0.876114
0.187 	0.861926
0.188 	0.847874
0.189 	0.833956
0.190 	0.820175
0.191 	0.806529
0.192 	0.793020
0.193 	0.779649
0.194 	0.766415
0.195 	0.753319
0.196 	0.740362
0.197 	0.727543
0.198 	0.714863
0.199 	0.702322
0.200 	0.689921
0.201 	0.677658
0.202 	0.665536
0.203 	0.653553
0.204 	0.641709
0.205 	0.630005
0.206 	0.618441
0.207 	0.607016
0.208 	0.595730
0.209 	0.584584
0.210 	0.573576
0.211 	0.562708
0.212 	0.551978
0.213 	0.541386
0.214 	0.530932
0.215 	0.520616
0.216 	0.510437
0.217 	0.500395
0.218 	0.490489
0.219 	0.480719
0.220 	0.471084
0.221 	0.461585
0.222 	0.452219
0.223 	0.442988
0.224 	0.433889
0.225 	0.424923
0.226 	0.416089
0.227 	0.407386
0.228 	0.398814
0.229 	0.390371
0.230 	0.382057
0.231 	0.373871
0.232 	0.365813
0.233 	0.357881
0.234 	0.350075
0.235 	0.342394
0.236 	0.334837
0.237 	0.327402
0.238 	0.320090
0.239 	0.312899
0.240 	0.305828
0.241 	0.298877
0.242 	0.292043
0.243 	0.285327
0.244 	0.278728
0.245 	0.272243
0.246 	0.265872
0.247 	0.259615
0.248 	0.253469
0.249 	0.247435
0.250 	0.241510
0.251 	0.235694
0.252 	0.229985
0.253 	0.224383
0.254 	0.218886
0.255 	0.213493
0.256 	0.208204
0.257 	0.203016
0.258 	0.197929
0.259 	0.192941
0.260 	0.188052
0.261 	0.183259
0.262 	0.178563
0.263 	0.173962
0.264 	0.169454
0.265 	0.165038
0.266 	0.160713
0.267 	0.156479
0.268 	0.152333
0.269 	0.148274
0.270 	0.144302
0.271 	0.140415
0.272 	0.136612
0.273 	0.132892
0.274 	0.129253
0.275 	0.125694
0.276 	0.122215
0.277 	0.118813
0.278 	0.115488
0.279 	0.112239
0.280 	0.109063
0.281 	0.105961
0.282 	0.102931
0.283 	0.099972
0.284 	0.097082
0.285 	0.094261
0.286 	0.091507
0.287 	0.088819
0.288 	0.086196
0.289 	0.083637
0.290 	0.081141
0.291 	0.078707
0.292 	0.076333
0.293 	0.074018
0.294 	0.071762
0.295 	0.069562
0.296 	0.067419
0.297 	0.065331
0.298 	0.063298
0.299 	0.061317
0.300 	0.059388
0.301 	0.057510
0.302 	0.055681
0.303 	0.053902
0.304 	0.052171
0.305 	0.050486
0.306 	0.048847
0.307 	0.047254
0.308 	0.045704
0.309 	0.044198
0.310 	0.042733
0.311 	0.041310
0.312 	0.039927
0.313 	0.038584
0.314 	0.037279
0.315 	0.036011
0.316 	0.034781
0.317 	0.033587
0.318 	0.032427
0.319 	0.031302
0.320 	0.030211
0.321 	0.029152
0.322 	0.028125
0.323 	0.027129
0.324 	0.026164
0.325 	0.025228
0.326 	0.024321
0.327 	0.023443
0.328 	0.022592
0.329 	0.021767
0.330 	0.020969
0.331 	0.020196
0.332 	0.019448
0.333 	0.018724
0.334 	0.018024
0.335 	0.017346
0.336 	0.016691
0.337 	0.016057
0.338 	0.015444
0.339 	0.014852
0.340 	0.014279
0.341 	0.013726
0.342 	0.013192
0.343 	0.012676
0.344 	0.012178
0.345 	0.011697
0.346 	0.011232
0.347 	0.010784
0.348 	0.010352
0.349 	0.009935
0.350 	0.009533
0.351 	0.009145
0.352 	0.008771
0.353 	0.008411
0.354 	0.008064
0.355 	0.007729
0.356 	0.007407
0.357 	0.007097
0.358 	0.006798
0.359 	0.006511
0.360 	0.006234
0.361 	0.005968
0.362 	0.005712
0.363 	0.005466
0.364 	0.005229
0.365 	0.005001
0.366 	0.004783
0.367 	0.004572
0.368 	0.004371
0.369 	0.004177
0.370 	0.003990
0.371 	0.003812
0.372 	0.003640
0.373 	0.003476
0.374 	0.003318
0.375 	0.003166
0.376 	0.003021
0.377 	0.002882
0.378 	0.002748
0.379 	0.002620
0.380 	0.002498
0.381 	0.002381
0.382 	0.002268
0.383 	0.002161
0.384 	0.002058
0.385 	0.001959
0.386 	0.001865
0.387 	0.001775
0.388 	0.001689
0.389 	0.001606
0.390 	0.001528
0.391 	0.001453
0.392 	0.001381
0.393 	0.001312
0.394 	0.001247
0.395 	0.001184
0.396 	0.001124
0.397 	0.001068
0.398 	0.001013
0.399 	0.000961
0.400 	0.000912
0.401 	0.000865
0.402 	0.000820
0.403 	0.000778
0.404 	0.000737
0.405 	0.000698
0.406 	0.000661
0.407 	0.000626
0.408 	0.000593
0.409 	0.000561
0.410 	0.000531
0.411 	0.000502
0.412 	0.000475
0.413 	0.000449
0.414 	0.000425
0.415 	0.000401
0.416 	0.000379
0.417 	0.000358
0.418 	0.000338
0.419 	0.000319
0.420 	0.000301
0.421 	0.000284
0.422 	0.000268
0.423 	0.000253
0.424 	0.000238
0.425 	0.000225
0.426 	0.000212
0.427 	0.000199
0.428 	0.000188
0.429 	0.000177
0.430 	0.000166
0.431 	0.000156
0.432 	0.000147
0.433 	0.000138
0.434 	0.000130
0.435 	0.000122
0.436 	0.000115
0.437 	0.000108
0.438 	0.000101
0.439 	0.000095
0.440 	0.000089
0.441 	0.000084
0.442 	0.000079
0.443 	0.000074
0.444 	0.000069
0.445 	0.000065
0.446 	0.000061
0.447 	0.000057
0.448 	0.000053
0.449 	0.000050
0.450 	0.000047
0.451 	0.000044
0.452 	0.000041
0.453 	0.000038
0.454 	0.000036
0.455 	0.000033
0.456 	0.000031
0.457 	0.000029
0.458 	0.000027
0.459 	0.000025
0.460 	0.000024
0.461 	0.000022
0.462 	0.000020
0.463 	0.000019
0.464 	0.000018
0.465 	0.000017
0.466 	0.000015
0.467 	0.000014
0.468 	0.000013
0.469 	0.000012
0.470 	0.000012
0.471 	0.000011
0.472 	0.000010
0.473 	0.000009
0.474 	0.000009
0.475 	0.000008
0.476 	0.000007
0.477 	0.000007
0.478 	0.000006
0.479 	0.000006
0.480 	0.000005
0.481 	0.000005
0.482 	0.000005
0.483 	0.000004
0.484 	0.000004
0.485 	0.000004
0.486 	0.000003
0.487 	0.000003
0.488 	0.000003
0.489 	0.000003
0.490 	0.000002
0.491 	0.000002
0.492 	0.000002
0.493 	0.000002
0.494 	0.000002
0.495 	0.000002
0.496 	0.000002
0.497 	0.000001
0.498 	0.000001
0.499 	0.000001
0.500 	0.000001
/

\plot 0.23 2.0 0.28 2.0 /
\put {$G(t)$}  [l] at 0.30  2.0

\setsolid
\axis bottom label {} shiftedto y=0 ticks numbered from 0.0 to 0.5 by 0.1 unlabeled short from 0.05 to 0.45 by 0.10 /
\axis left   label {} shiftedto x=0 ticks numbered from 1.0 to 6.0 by 1.0 unlabeled short from 0.5  to 5.5  by 1.0 /

\put {$t$}     [r] at 0.495  0.30

\endpicture     \]
\vspace{-8mm}

\caption{Comparison of the three functions \(G(t)\),  \(-G'(t)/t\) and \((\cdots)G(t)\) in \eqref{Jensen_4}.}

\end{figure}


%% file: fig_03.tex




\baselineskip1.8em

\begin{figure} [hbt]

\[  \beginpicture   \setlinear

\setcoordinatesystem units <110mm,4mm>
\setplotarea x from 0.0 to 1.0, y from 0.0 to 15.0



\thicklines


\setdashpattern <1.9mm, 0.7mm, 0.2mm, 0.7mm>
\plot
0.000 	1.490240
0.001 	1.490252
0.002 	1.490287
0.003 	1.490347
0.004 	1.490430
0.005 	1.490538
0.006 	1.490669
0.007 	1.490824
0.008 	1.491003
0.009 	1.491206
0.010 	1.491432
0.011 	1.491683
0.012 	1.491957
0.013 	1.492255
0.014 	1.492577
0.015 	1.492923
0.016 	1.493293
0.017 	1.493687
0.018 	1.494104
0.019 	1.494546
0.020 	1.495011
0.021 	1.495500
0.022 	1.496014
0.023 	1.496551
0.024 	1.497112
0.025 	1.497697
0.026 	1.498306
0.027 	1.498939
0.028 	1.499596
0.029 	1.500277
0.030 	1.500982
0.031 	1.501711
0.032 	1.502464
0.033 	1.503242
0.034 	1.504043
0.035 	1.504868
0.036 	1.505717
0.037 	1.506591
0.038 	1.507488
0.039 	1.508410
0.040 	1.509356
0.041 	1.510325
0.042 	1.511320
0.043 	1.512338
0.044 	1.513380
0.045 	1.514447
0.046 	1.515538
0.047 	1.516653
0.048 	1.517792
0.049 	1.518956
0.050 	1.520144
0.051 	1.521356
0.052 	1.522593
0.053 	1.523854
0.054 	1.525139
0.055 	1.526449
0.056 	1.527783
0.057 	1.529142
0.058 	1.530525
0.059 	1.531933
0.060 	1.533365
0.061 	1.534822
0.062 	1.536303
0.063 	1.537809
0.064 	1.539339
0.065 	1.540894
0.066 	1.542474
0.067 	1.544078
0.068 	1.545707
0.069 	1.547361
0.070 	1.549040
0.071 	1.550743
0.072 	1.552471
0.073 	1.554224
0.074 	1.556002
0.075 	1.557805
0.076 	1.559633
0.077 	1.561485
0.078 	1.563363
0.079 	1.565266
0.080 	1.567193
0.081 	1.569146
0.082 	1.571124
0.083 	1.573127
0.084 	1.575155
0.085 	1.577209
0.086 	1.579287
0.087 	1.581391
0.088 	1.583520
0.089 	1.585675
0.090 	1.587855
0.091 	1.590060
0.092 	1.592291
0.093 	1.594547
0.094 	1.596829
0.095 	1.599136
0.096 	1.601469
0.097 	1.603827
0.098 	1.606212
0.099 	1.608621
0.100 	1.611057
0.101 	1.613518
0.102 	1.616006
0.103 	1.618519
0.104 	1.621058
0.105 	1.623622
0.106 	1.626213
0.107 	1.628830
0.108 	1.631473
0.109 	1.634142
0.110 	1.636837
0.111 	1.639559
0.112 	1.642306
0.113 	1.645080
0.114 	1.647880
0.115 	1.650707
0.116 	1.653560
0.117 	1.656439
0.118 	1.659345
0.119 	1.662278
0.120 	1.665237
0.121 	1.668223
0.122 	1.671235
0.123 	1.674274
0.124 	1.677340
0.125 	1.680433
0.126 	1.683553
0.127 	1.686699
0.128 	1.689873
0.129 	1.693074
0.130 	1.696301
0.131 	1.699556
0.132 	1.702838
0.133 	1.706148
0.134 	1.709484
0.135 	1.712848
0.136 	1.716239
0.137 	1.719658
0.138 	1.723105
0.139 	1.726579
0.140 	1.730080
0.141 	1.733609
0.142 	1.737166
0.143 	1.740751
0.144 	1.744364
0.145 	1.748004
0.146 	1.751673
0.147 	1.755369
0.148 	1.759094
0.149 	1.762846
0.150 	1.766627
0.151 	1.770437
0.152 	1.774274
0.153 	1.778140
0.154 	1.782034
0.155 	1.785957
0.156 	1.789909
0.157 	1.793889
0.158 	1.797898
0.159 	1.801935
0.160 	1.806002
0.161 	1.810097
0.162 	1.814221
0.163 	1.818374
0.164 	1.822557
0.165 	1.826768
0.166 	1.831009
0.167 	1.835279
0.168 	1.839579
0.169 	1.843908
0.170 	1.848266
0.171 	1.852654
0.172 	1.857071
0.173 	1.861519
0.174 	1.865996
0.175 	1.870503
0.176 	1.875040
0.177 	1.879606
0.178 	1.884203
0.179 	1.888830
0.180 	1.893488
0.181 	1.898175
0.182 	1.902893
0.183 	1.907642
0.184 	1.912421
0.185 	1.917230
0.186 	1.922070
0.187 	1.926941
0.188 	1.931843
0.189 	1.936776
0.190 	1.941740
0.191 	1.946734
0.192 	1.951760
0.193 	1.956817
0.194 	1.961906
0.195 	1.967026
0.196 	1.972177
0.197 	1.977360
0.198 	1.982574
0.199 	1.987821
0.200 	1.993099
0.201 	1.998409
0.202 	2.003751
0.203 	2.009124
0.204 	2.014531
0.205 	2.019969
0.206 	2.025440
0.207 	2.030943
0.208 	2.036478
0.209 	2.042046
0.210 	2.047647
0.211 	2.053281
0.212 	2.058947
0.213 	2.064647
0.214 	2.070379
0.215 	2.076144
0.216 	2.081943
0.217 	2.087775
0.218 	2.093641
0.219 	2.099540
0.220 	2.105472
0.221 	2.111439
0.222 	2.117439
0.223 	2.123472
0.224 	2.129540
0.225 	2.135642
0.226 	2.141778
0.227 	2.147949
0.228 	2.154154
0.229 	2.160393
0.230 	2.166667
0.231 	2.172975
0.232 	2.179318
0.233 	2.185696
0.234 	2.192109
0.235 	2.198558
0.236 	2.205041
0.237 	2.211560
0.238 	2.218113
0.239 	2.224703
0.240 	2.231328
0.241 	2.237989
0.242 	2.244685
0.243 	2.251418
0.244 	2.258186
0.245 	2.264991
0.246 	2.271832
0.247 	2.278709
0.248 	2.285623
0.249 	2.292573
0.250 	2.299560
0.251 	2.306584
0.252 	2.313644
0.253 	2.320742
0.254 	2.327877
0.255 	2.335049
0.256 	2.342258
0.257 	2.349505
0.258 	2.356790
0.259 	2.364112
0.260 	2.371472
0.261 	2.378870
0.262 	2.386306
0.263 	2.393780
0.264 	2.401293
0.265 	2.408843
0.266 	2.416433
0.267 	2.424061
0.268 	2.431728
0.269 	2.439434
0.270 	2.447179
0.271 	2.454963
0.272 	2.462786
0.273 	2.470649
0.274 	2.478551
0.275 	2.486493
0.276 	2.494474
0.277 	2.502496
0.278 	2.510557
0.279 	2.518659
0.280 	2.526801
0.281 	2.534984
0.282 	2.543207
0.283 	2.551471
0.284 	2.559775
0.285 	2.568121
0.286 	2.576507
0.287 	2.584935
0.288 	2.593404
0.289 	2.601915
0.290 	2.610467
0.291 	2.619062
0.292 	2.627698
0.293 	2.636376
0.294 	2.645096
0.295 	2.653858
0.296 	2.662663
0.297 	2.671511
0.298 	2.680401
0.299 	2.689335
0.300 	2.698311
0.301 	2.707330
0.302 	2.716393
0.303 	2.725499
0.304 	2.734649
0.305 	2.743843
0.306 	2.753080
0.307 	2.762362
0.308 	2.771688
0.309 	2.781058
0.310 	2.790472
0.311 	2.799931
0.312 	2.809435
0.313 	2.818984
0.314 	2.828579
0.315 	2.838218
0.316 	2.847903
0.317 	2.857633
0.318 	2.867409
0.319 	2.877231
0.320 	2.887099
0.321 	2.897013
0.322 	2.906974
0.323 	2.916981
0.324 	2.927034
0.325 	2.937135
0.326 	2.947282
0.327 	2.957477
0.328 	2.967719
0.329 	2.978008
0.330 	2.988345
0.331 	2.998730
0.332 	3.009163
0.333 	3.019644
0.334 	3.030174
0.335 	3.040752
0.336 	3.051378
0.337 	3.062053
0.338 	3.072778
0.339 	3.083551
0.340 	3.094374
0.341 	3.105246
0.342 	3.116168
0.343 	3.127140
0.344 	3.138162
0.345 	3.149234
0.346 	3.160357
0.347 	3.171530
0.348 	3.182754
0.349 	3.194029
0.350 	3.205354
0.351 	3.216732
0.352 	3.228160
0.353 	3.239640
0.354 	3.251172
0.355 	3.262757
0.356 	3.274393
0.357 	3.286082
0.358 	3.297823
0.359 	3.309617
0.360 	3.321464
0.361 	3.333364
0.362 	3.345318
0.363 	3.357325
0.364 	3.369385
0.365 	3.381500
0.366 	3.393669
0.367 	3.405892
0.368 	3.418169
0.369 	3.430501
0.370 	3.442888
0.371 	3.455331
0.372 	3.467828
0.373 	3.480381
0.374 	3.492990
0.375 	3.505654
0.376 	3.518375
0.377 	3.531152
0.378 	3.543985
0.379 	3.556875
0.380 	3.569822
0.381 	3.582827
0.382 	3.595888
0.383 	3.609007
0.384 	3.622184
0.385 	3.635419
0.386 	3.648711
0.387 	3.662063
0.388 	3.675473
0.389 	3.688941
0.390 	3.702469
0.391 	3.716056
0.392 	3.729702
0.393 	3.743408
0.394 	3.757174
0.395 	3.771000
0.396 	3.784887
0.397 	3.798834
0.398 	3.812841
0.399 	3.826910
0.400 	3.841040
0.401 	3.855231
0.402 	3.869485
0.403 	3.883800
0.404 	3.898177
0.405 	3.912616
0.406 	3.927118
0.407 	3.941683
0.408 	3.956311
0.409 	3.971003
0.410 	3.985758
0.411 	4.000576
0.412 	4.015459
0.413 	4.030406
0.414 	4.045417
0.415 	4.060493
0.416 	4.075634
0.417 	4.090841
0.418 	4.106112
0.419 	4.121450
0.420 	4.136853
0.421 	4.152323
0.422 	4.167859
0.423 	4.183462
0.424 	4.199131
0.425 	4.214868
0.426 	4.230672
0.427 	4.246544
0.428 	4.262484
0.429 	4.278492
0.430 	4.294569
0.431 	4.310714
0.432 	4.326928
0.433 	4.343212
0.434 	4.359565
0.435 	4.375988
0.436 	4.392480
0.437 	4.409043
0.438 	4.425677
0.439 	4.442381
0.440 	4.459157
0.441 	4.476003
0.442 	4.492922
0.443 	4.509912
0.444 	4.526975
0.445 	4.544109
0.446 	4.561317
0.447 	4.578598
0.448 	4.595951
0.449 	4.613379
0.450 	4.630880
0.451 	4.648455
0.452 	4.666105
0.453 	4.683829
0.454 	4.701629
0.455 	4.719503
0.456 	4.737453
0.457 	4.755479
0.458 	4.773581
0.459 	4.791759
0.460 	4.810014
0.461 	4.828346
0.462 	4.846755
0.463 	4.865241
0.464 	4.883806
0.465 	4.902449
0.466 	4.921170
0.467 	4.939970
0.468 	4.958849
0.469 	4.977807
0.470 	4.996845
0.471 	5.015962
0.472 	5.035161
0.473 	5.054439
0.474 	5.073799
0.475 	5.093239
0.476 	5.112762
0.477 	5.132366
0.478 	5.152052
0.479 	5.171820
0.480 	5.191671
0.481 	5.211606
0.482 	5.231624
0.483 	5.251725
0.484 	5.271910
0.485 	5.292180
0.486 	5.312535
0.487 	5.332974
0.488 	5.353499
0.489 	5.374110
0.490 	5.394806
0.491 	5.415589
0.492 	5.436458
0.493 	5.457415
0.494 	5.478458
0.495 	5.499590
0.496 	5.520809
0.497 	5.542117
0.498 	5.563513
0.499 	5.584999
0.500 	5.606574
0.501 	5.628238
0.502 	5.649993
0.503 	5.671837
0.504 	5.693773
0.505 	5.715800
0.506 	5.737918
0.507 	5.760128
0.508 	5.782430
0.509 	5.804825
0.510 	5.827313
0.511 	5.849894
0.512 	5.872568
0.513 	5.895336
0.514 	5.918199
0.515 	5.941157
0.516 	5.964209
0.517 	5.987357
0.518 	6.010601
0.519 	6.033940
0.520 	6.057377
0.521 	6.080910
0.522 	6.104541
0.523 	6.128269
0.524 	6.152095
0.525 	6.176020
0.526 	6.200044
0.527 	6.224167
0.528 	6.248389
0.529 	6.272711
0.530 	6.297134
0.531 	6.321658
0.532 	6.346282
0.533 	6.371008
0.534 	6.395836
0.535 	6.420767
0.536 	6.445800
0.537 	6.470936
0.538 	6.496176
0.539 	6.521520
0.540 	6.546968
0.541 	6.572521
0.542 	6.598179
0.543 	6.623943
0.544 	6.649812
0.545 	6.675788
0.546 	6.701871
0.547 	6.728061
0.548 	6.754359
0.549 	6.780764
0.550 	6.807279
0.551 	6.833902
0.552 	6.860634
0.553 	6.887476
0.554 	6.914429
0.555 	6.941492
0.556 	6.968666
0.557 	6.995952
0.558 	7.023349
0.559 	7.050859
0.560 	7.078482
0.561 	7.106218
0.562 	7.134068
0.563 	7.162032
0.564 	7.190110
0.565 	7.218304
0.566 	7.246613
0.567 	7.275038
0.568 	7.303579
0.569 	7.332237
0.570 	7.361013
0.571 	7.389906
0.572 	7.418918
0.573 	7.448048
0.574 	7.477297
0.575 	7.506666
0.576 	7.536156
0.577 	7.565765
0.578 	7.595496
0.579 	7.625348
0.580 	7.655323
0.581 	7.685420
0.582 	7.715639
0.583 	7.745983
0.584 	7.776450
0.585 	7.807041
0.586 	7.837758
0.587 	7.868600
0.588 	7.899567
0.589 	7.930662
0.590 	7.961883
0.591 	7.993231
0.592 	8.024708
0.593 	8.056312
0.594 	8.088046
0.595 	8.119909
0.596 	8.151902
0.597 	8.184026
0.598 	8.216280
0.599 	8.248666
0.600 	8.281184
0.601 	8.313834
0.602 	8.346617
0.603 	8.379534
0.604 	8.412585
0.605 	8.445771
0.606 	8.479091
0.607 	8.512548
0.608 	8.546140
0.609 	8.579869
0.610 	8.613736
0.611 	8.647740
0.612 	8.681883
0.613 	8.716165
0.614 	8.750586
0.615 	8.785147
0.616 	8.819849
0.617 	8.854691
0.618 	8.889676
0.619 	8.924803
0.620 	8.960072
0.621 	8.995485
0.622 	9.031042
0.623 	9.066743
0.624 	9.102589
0.625 	9.138582
0.626 	9.174720
0.627 	9.211005
0.628 	9.247437
0.629 	9.284018
0.630 	9.320746
0.631 	9.357624
0.632 	9.394652
0.633 	9.431830
0.634 	9.469159
0.635 	9.506640
0.636 	9.544272
0.637 	9.582058
0.638 	9.619996
0.639 	9.658089
0.640 	9.696336
0.641 	9.734738
0.642 	9.773296
0.643 	9.812010
0.644 	9.850882
0.645 	9.889911
0.646 	9.929098
0.647 	9.968444
0.648 	10.007949
0.649 	10.047615
0.650 	10.087441
0.651 	10.127429
0.652 	10.167579
0.653 	10.207892
0.654 	10.248368
0.655 	10.289007
0.656 	10.329812
0.657 	10.370782
0.658 	10.411917
0.659 	10.453220
0.660 	10.494689
0.661 	10.536327
0.662 	10.578133
0.663 	10.620108
0.664 	10.662253
0.665 	10.704569
0.666 	10.747056
0.667 	10.789715
0.668 	10.832547
0.669 	10.875552
0.670 	10.918731
0.671 	10.962084
0.672 	11.005613
0.673 	11.049319
0.674 	11.093201
0.675 	11.137260
0.676 	11.181498
0.677 	11.225914
0.678 	11.270510
0.679 	11.315287
0.680 	11.360244
0.681 	11.405384
0.682 	11.450705
0.683 	11.496210
0.684 	11.541899
0.685 	11.587773
0.686 	11.633832
0.687 	11.680077
0.688 	11.726509
0.689 	11.773129
0.690 	11.819937
0.691 	11.866934
0.692 	11.914121
0.693 	11.961499
0.694 	12.009068
0.695 	12.056829
0.696 	12.104783
0.697 	12.152931
0.698 	12.201273
0.699 	12.249810
0.700 	12.298544
0.701 	12.347474
0.702 	12.396602
0.703 	12.445928
0.704 	12.495453
0.705 	12.545178
0.706 	12.595104
0.707 	12.645232
0.708 	12.695561
0.709 	12.746094
0.710 	12.796831
0.711 	12.847773
0.712 	12.898920
0.713 	12.950273
0.714 	13.001834
0.715 	13.053603
0.716 	13.105580
0.717 	13.157768
0.718 	13.210166
0.719 	13.262775
0.720 	13.315596
0.721 	13.368631
0.722 	13.421879
0.723 	13.475342
0.724 	13.529021
0.725 	13.582916
0.726 	13.637028
0.727 	13.691359
0.728 	13.745909
0.729 	13.800678
0.730 	13.855669
0.731 	13.910881
0.732 	13.966316
0.733 	14.021974
0.734 	14.077856
0.735 	14.133964
0.736 	14.190298
0.737 	14.246859
0.738 	14.303648
0.739 	14.360666
0.740 	14.417914
0.741 	14.475392
0.742 	14.533102
0.743 	14.591044
0.744 	14.649220
0.745 	14.707630
0.746 	14.766276
0.747 	14.825158
0.748 	14.884277
0.749 	14.943634
0.750 	15.003230
/
\plot 0.21 13.0 0.3 13.0 /
\put {$  a_0 \cosh(4t)$}  [l] at 0.32  13

\setsolid
\plot
0.000 	1.490240
0.001 	1.490247
0.002 	1.490270
0.003 	1.490308
0.004 	1.490361
0.005 	1.490429
0.006 	1.490513
0.007 	1.490611
0.008 	1.490725
0.009 	1.490854
0.010 	1.490998
0.011 	1.491158
0.012 	1.491332
0.013 	1.491522
0.014 	1.491726
0.015 	1.491946
0.016 	1.492181
0.017 	1.492431
0.018 	1.492696
0.019 	1.492977
0.020 	1.493272
0.021 	1.493583
0.022 	1.493908
0.023 	1.494249
0.024 	1.494605
0.025 	1.494975
0.026 	1.495361
0.027 	1.495762
0.028 	1.496178
0.029 	1.496609
0.030 	1.497055
0.031 	1.497516
0.032 	1.497992
0.033 	1.498483
0.034 	1.498989
0.035 	1.499510
0.036 	1.500046
0.037 	1.500596
0.038 	1.501162
0.039 	1.501743
0.040 	1.502338
0.041 	1.502949
0.042 	1.503574
0.043 	1.504214
0.044 	1.504869
0.045 	1.505538
0.046 	1.506223
0.047 	1.506922
0.048 	1.507636
0.049 	1.508365
0.050 	1.509109
0.051 	1.509867
0.052 	1.510640
0.053 	1.511428
0.054 	1.512230
0.055 	1.513047
0.056 	1.513879
0.057 	1.514725
0.058 	1.515586
0.059 	1.516462
0.060 	1.517352
0.061 	1.518256
0.062 	1.519175
0.063 	1.520109
0.064 	1.521057
0.065 	1.522020
0.066 	1.522997
0.067 	1.523989
0.068 	1.524995
0.069 	1.526015
0.070 	1.527050
0.071 	1.528099
0.072 	1.529162
0.073 	1.530240
0.074 	1.531332
0.075 	1.532439
0.076 	1.533559
0.077 	1.534694
0.078 	1.535844
0.079 	1.537007
0.080 	1.538185
0.081 	1.539376
0.082 	1.540582
0.083 	1.541802
0.084 	1.543036
0.085 	1.544285
0.086 	1.545547
0.087 	1.546824
0.088 	1.548114
0.089 	1.549418
0.090 	1.550737
0.091 	1.552069
0.092 	1.553416
0.093 	1.554776
0.094 	1.556151
0.095 	1.557539
0.096 	1.558941
0.097 	1.560357
0.098 	1.561787
0.099 	1.563230
0.100 	1.564688
0.101 	1.566159
0.102 	1.567644
0.103 	1.569143
0.104 	1.570656
0.105 	1.572182
0.106 	1.573722
0.107 	1.575276
0.108 	1.576843
0.109 	1.578424
0.110 	1.580019
0.111 	1.581627
0.112 	1.583249
0.113 	1.584885
0.114 	1.586534
0.115 	1.588197
0.116 	1.589873
0.117 	1.591563
0.118 	1.593266
0.119 	1.594983
0.120 	1.596713
0.121 	1.598457
0.122 	1.600214
0.123 	1.601985
0.124 	1.603769
0.125 	1.605567
0.126 	1.607377
0.127 	1.609202
0.128 	1.611039
0.129 	1.612891
0.130 	1.614755
0.131 	1.616633
0.132 	1.618524
0.133 	1.620428
0.134 	1.622346
0.135 	1.624277
0.136 	1.626221
0.137 	1.628179
0.138 	1.630150
0.139 	1.632134
0.140 	1.634131
0.141 	1.636142
0.142 	1.638166
0.143 	1.640203
0.144 	1.642253
0.145 	1.644316
0.146 	1.646393
0.147 	1.648483
0.148 	1.650586
0.149 	1.652702
0.150 	1.654831
0.151 	1.656974
0.152 	1.659129
0.153 	1.661298
0.154 	1.663480
0.155 	1.665675
0.156 	1.667883
0.157 	1.670104
0.158 	1.672339
0.159 	1.674586
0.160 	1.676847
0.161 	1.679121
0.162 	1.681407
0.163 	1.683707
0.164 	1.686020
0.165 	1.688346
0.166 	1.690685
0.167 	1.693038
0.168 	1.695403
0.169 	1.697781
0.170 	1.700173
0.171 	1.702577
0.172 	1.704995
0.173 	1.707425
0.174 	1.709869
0.175 	1.712326
0.176 	1.714796
0.177 	1.717279
0.178 	1.719775
0.179 	1.722284
0.180 	1.724806
0.181 	1.727342
0.182 	1.729890
0.183 	1.732451
0.184 	1.735026
0.185 	1.737614
0.186 	1.740214
0.187 	1.742828
0.188 	1.745455
0.189 	1.748095
0.190 	1.750748
0.191 	1.753414
0.192 	1.756094
0.193 	1.758786
0.194 	1.761492
0.195 	1.764211
0.196 	1.766943
0.197 	1.769688
0.198 	1.772446
0.199 	1.775217
0.200 	1.778002
0.201 	1.780800
0.202 	1.783610
0.203 	1.786435
0.204 	1.789272
0.205 	1.792122
0.206 	1.794986
0.207 	1.797863
0.208 	1.800753
0.209 	1.803657
0.210 	1.806573
0.211 	1.809503
0.212 	1.812447
0.213 	1.815403
0.214 	1.818373
0.215 	1.821356
0.216 	1.824353
0.217 	1.827362
0.218 	1.830386
0.219 	1.833422
0.220 	1.836472
0.221 	1.839535
0.222 	1.842612
0.223 	1.845702
0.224 	1.848806
0.225 	1.851922
0.226 	1.855053
0.227 	1.858197
0.228 	1.861354
0.229 	1.864525
0.230 	1.867709
0.231 	1.870907
0.232 	1.874119
0.233 	1.877344
0.234 	1.880582
0.235 	1.883834
0.236 	1.887100
0.237 	1.890380
0.238 	1.893673
0.239 	1.896979
0.240 	1.900300
0.241 	1.903634
0.242 	1.906982
0.243 	1.910343
0.244 	1.913719
0.245 	1.917108
0.246 	1.920510
0.247 	1.923927
0.248 	1.927358
0.249 	1.930802
0.250 	1.934260
0.251 	1.937732
0.252 	1.941219
0.253 	1.944719
0.254 	1.948232
0.255 	1.951760
0.256 	1.955302
0.257 	1.958858
0.258 	1.962428
0.259 	1.966012
0.260 	1.969610
0.261 	1.973223
0.262 	1.976849
0.263 	1.980490
0.264 	1.984144
0.265 	1.987813
0.266 	1.991497
0.267 	1.995194
0.268 	1.998906
0.269 	2.002632
0.270 	2.006372
0.271 	2.010127
0.272 	2.013896
0.273 	2.017679
0.274 	2.021477
0.275 	2.025290
0.276 	2.029116
0.277 	2.032958
0.278 	2.036814
0.279 	2.040684
0.280 	2.044569
0.281 	2.048469
0.282 	2.052383
0.283 	2.056312
0.284 	2.060256
0.285 	2.064214
0.286 	2.068188
0.287 	2.072176
0.288 	2.076178
0.289 	2.080196
0.290 	2.084229
0.291 	2.088276
0.292 	2.092338
0.293 	2.096416
0.294 	2.100508
0.295 	2.104616
0.296 	2.108738
0.297 	2.112876
0.298 	2.117028
0.299 	2.121196
0.300 	2.125379
0.301 	2.129577
0.302 	2.133791
0.303 	2.138020
0.304 	2.142264
0.305 	2.146523
0.306 	2.150798
0.307 	2.155088
0.308 	2.159394
0.309 	2.163715
0.310 	2.168052
0.311 	2.172404
0.312 	2.176772
0.313 	2.181155
0.314 	2.185554
0.315 	2.189969
0.316 	2.194399
0.317 	2.198846
0.318 	2.203308
0.319 	2.207786
0.320 	2.212279
0.321 	2.216789
0.322 	2.221315
0.323 	2.225856
0.324 	2.230414
0.325 	2.234987
0.326 	2.239577
0.327 	2.244183
0.328 	2.248805
0.329 	2.253443
0.330 	2.258098
0.331 	2.262768
0.332 	2.267456
0.333 	2.272159
0.334 	2.276879
0.335 	2.281615
0.336 	2.286368
0.337 	2.291137
0.338 	2.295923
0.339 	2.300725
0.340 	2.305544
0.341 	2.310380
0.342 	2.315233
0.343 	2.320102
0.344 	2.324988
0.345 	2.329891
0.346 	2.334810
0.347 	2.339747
0.348 	2.344701
0.349 	2.349671
0.350 	2.354659
0.351 	2.359664
0.352 	2.364686
0.353 	2.369725
0.354 	2.374781
0.355 	2.379855
0.356 	2.384946
0.357 	2.390054
0.358 	2.395180
0.359 	2.400323
0.360 	2.405484
0.361 	2.410662
0.362 	2.415858
0.363 	2.421071
0.364 	2.426302
0.365 	2.431551
0.366 	2.436818
0.367 	2.442102
0.368 	2.447404
0.369 	2.452725
0.370 	2.458063
0.371 	2.463419
0.372 	2.468793
0.373 	2.474185
0.374 	2.479596
0.375 	2.485024
0.376 	2.490471
0.377 	2.495936
0.378 	2.501420
0.379 	2.506922
0.380 	2.512442
0.381 	2.517981
0.382 	2.523538
0.383 	2.529114
0.384 	2.534709
0.385 	2.540322
0.386 	2.545954
0.387 	2.551605
0.388 	2.557274
0.389 	2.562963
0.390 	2.568670
0.391 	2.574397
0.392 	2.580142
0.393 	2.585907
0.394 	2.591690
0.395 	2.597493
0.396 	2.603316
0.397 	2.609157
0.398 	2.615018
0.399 	2.620898
0.400 	2.626798
0.401 	2.632717
0.402 	2.638656
0.403 	2.644615
0.404 	2.650593
0.405 	2.656591
0.406 	2.662609
0.407 	2.668646
0.408 	2.674704
0.409 	2.680781
0.410 	2.686879
0.411 	2.692997
0.412 	2.699134
0.413 	2.705292
0.414 	2.711470
0.415 	2.717669
0.416 	2.723888
0.417 	2.730127
0.418 	2.736387
0.419 	2.742667
0.420 	2.748968
0.421 	2.755290
0.422 	2.761632
0.423 	2.767995
0.424 	2.774379
0.425 	2.780784
0.426 	2.787210
0.427 	2.793656
0.428 	2.800124
0.429 	2.806613
0.430 	2.813123
0.431 	2.819655
0.432 	2.826208
0.433 	2.832782
0.434 	2.839377
0.435 	2.845994
0.436 	2.852633
0.437 	2.859293
0.438 	2.865975
0.439 	2.872679
0.440 	2.879405
0.441 	2.886152
0.442 	2.892922
0.443 	2.899713
0.444 	2.906527
0.445 	2.913362
0.446 	2.920220
0.447 	2.927100
0.448 	2.934003
0.449 	2.940928
0.450 	2.947875
0.451 	2.954845
0.452 	2.961838
0.453 	2.968853
0.454 	2.975891
0.455 	2.982951
0.456 	2.990035
0.457 	2.997142
0.458 	3.004271
0.459 	3.011424
0.460 	3.018600
0.461 	3.025799
0.462 	3.033021
0.463 	3.040266
0.464 	3.047536
0.465 	3.054828
0.466 	3.062144
0.467 	3.069484
0.468 	3.076848
0.469 	3.084235
0.470 	3.091646
0.471 	3.099081
0.472 	3.106540
0.473 	3.114023
0.474 	3.121530
0.475 	3.129062
0.476 	3.136617
0.477 	3.144198
0.478 	3.151802
0.479 	3.159431
0.480 	3.167085
0.481 	3.174763
0.482 	3.182466
0.483 	3.190194
0.484 	3.197946
0.485 	3.205724
0.486 	3.213527
0.487 	3.221354
0.488 	3.229207
0.489 	3.237085
0.490 	3.244989
0.491 	3.252918
0.492 	3.260872
0.493 	3.268852
0.494 	3.276858
0.495 	3.284889
0.496 	3.292946
0.497 	3.301029
0.498 	3.309139
0.499 	3.317274
0.500 	3.325435
0.501 	3.333622
0.502 	3.341836
0.503 	3.350076
0.504 	3.358342
0.505 	3.366635
0.506 	3.374955
0.507 	3.383301
0.508 	3.391674
0.509 	3.400074
0.510 	3.408501
0.511 	3.416955
0.512 	3.425436
0.513 	3.433944
0.514 	3.442479
0.515 	3.451042
0.516 	3.459632
0.517 	3.468250
0.518 	3.476895
0.519 	3.485568
0.520 	3.494269
0.521 	3.502998
0.522 	3.511755
0.523 	3.520539
0.524 	3.529352
0.525 	3.538193
0.526 	3.547063
0.527 	3.555960
0.528 	3.564887
0.529 	3.573841
0.530 	3.582825
0.531 	3.591837
0.532 	3.600878
0.533 	3.609948
0.534 	3.619047
0.535 	3.628175
0.536 	3.637332
0.537 	3.646519
0.538 	3.655735
0.539 	3.664980
0.540 	3.674255
0.541 	3.683560
0.542 	3.692895
0.543 	3.702259
0.544 	3.711653
0.545 	3.721077
0.546 	3.730532
0.547 	3.740016
0.548 	3.749531
0.549 	3.759077
0.550 	3.768653
0.551 	3.778259
0.552 	3.787896
0.553 	3.797564
0.554 	3.807263
0.555 	3.816993
0.556 	3.826754
0.557 	3.836546
0.558 	3.846369
0.559 	3.856224
0.560 	3.866111
0.561 	3.876029
0.562 	3.885978
0.563 	3.895960
0.564 	3.905973
0.565 	3.916018
0.566 	3.926095
0.567 	3.936205
0.568 	3.946347
0.569 	3.956521
0.570 	3.966728
0.571 	3.976967
0.572 	3.987239
0.573 	3.997544
0.574 	4.007882
0.575 	4.018253
0.576 	4.028657
0.577 	4.039094
0.578 	4.049565
0.579 	4.060068
0.580 	4.070606
0.581 	4.081177
0.582 	4.091782
0.583 	4.102421
0.584 	4.113094
0.585 	4.123801
0.586 	4.134542
0.587 	4.145317
0.588 	4.156127
0.589 	4.166971
0.590 	4.177850
0.591 	4.188764
0.592 	4.199713
0.593 	4.210696
0.594 	4.221715
0.595 	4.232769
0.596 	4.243858
0.597 	4.254982
0.598 	4.266142
0.599 	4.277338
0.600 	4.288570
0.601 	4.299837
0.602 	4.311141
0.603 	4.322480
0.604 	4.333856
0.605 	4.345268
0.606 	4.356716
0.607 	4.368202
0.608 	4.379723
0.609 	4.391282
0.610 	4.402878
0.611 	4.414510
0.612 	4.426180
0.613 	4.437887
0.614 	4.449632
0.615 	4.461414
0.616 	4.473233
0.617 	4.485091
0.618 	4.496986
0.619 	4.508919
0.620 	4.520891
0.621 	4.532901
0.622 	4.544949
0.623 	4.557036
0.624 	4.569161
0.625 	4.581325
0.626 	4.593528
0.627 	4.605770
0.628 	4.618051
0.629 	4.630371
0.630 	4.642731
0.631 	4.655130
0.632 	4.667569
0.633 	4.680048
0.634 	4.692567
0.635 	4.705125
0.636 	4.717724
0.637 	4.730363
0.638 	4.743043
0.639 	4.755763
0.640 	4.768524
0.641 	4.781325
0.642 	4.794168
0.643 	4.807052
0.644 	4.819977
0.645 	4.832943
0.646 	4.845951
0.647 	4.859000
0.648 	4.872091
0.649 	4.885224
0.650 	4.898399
0.651 	4.911616
0.652 	4.924876
0.653 	4.938178
0.654 	4.951522
0.655 	4.964910
0.656 	4.978340
0.657 	4.991813
0.658 	5.005329
0.659 	5.018889
0.660 	5.032492
0.661 	5.046138
0.662 	5.059829
0.663 	5.073563
0.664 	5.087341
0.665 	5.101163
0.666 	5.115030
0.667 	5.128941
0.668 	5.142896
0.669 	5.156896
0.670 	5.170941
0.671 	5.185032
0.672 	5.199167
0.673 	5.213347
0.674 	5.227573
0.675 	5.241845
0.676 	5.256163
0.677 	5.270526
0.678 	5.284935
0.679 	5.299391
0.680 	5.313893
0.681 	5.328441
0.682 	5.343036
0.683 	5.357678
0.684 	5.372367
0.685 	5.387103
0.686 	5.401886
0.687 	5.416717
0.688 	5.431595
0.689 	5.446521
0.690 	5.461495
0.691 	5.476517
0.692 	5.491587
0.693 	5.506706
0.694 	5.521873
0.695 	5.537089
0.696 	5.552354
0.697 	5.567667
0.698 	5.583030
0.699 	5.598442
0.700 	5.613904
0.701 	5.629415
0.702 	5.644977
0.703 	5.660588
0.704 	5.676249
0.705 	5.691961
0.706 	5.707723
0.707 	5.723536
0.708 	5.739399
0.709 	5.755314
0.710 	5.771280
0.711 	5.787297
0.712 	5.803365
0.713 	5.819485
0.714 	5.835657
0.715 	5.851882
0.716 	5.868158
0.717 	5.884486
0.718 	5.900868
0.719 	5.917301
0.720 	5.933788
0.721 	5.950328
0.722 	5.966921
0.723 	5.983567
0.724 	6.000267
0.725 	6.017021
0.726 	6.033828
0.727 	6.050690
0.728 	6.067606
0.729 	6.084576
0.730 	6.101602
0.731 	6.118682
0.732 	6.135816
0.733 	6.153007
0.734 	6.170252
0.735 	6.187553
0.736 	6.204910
0.737 	6.222323
0.738 	6.239792
0.739 	6.257317
0.740 	6.274898
0.741 	6.292537
0.742 	6.310232
0.743 	6.327984
0.744 	6.345794
0.745 	6.363660
0.746 	6.381585
0.747 	6.399567
0.748 	6.417608
0.749 	6.435706
0.750 	6.453863
0.751 	6.472078
0.752 	6.490353
0.753 	6.508686
0.754 	6.527078
0.755 	6.545530
0.756 	6.564041
0.757 	6.582612
0.758 	6.601243
0.759 	6.619935
0.760 	6.638686
0.761 	6.657498
0.762 	6.676371
0.763 	6.695305
0.764 	6.714300
0.765 	6.733356
0.766 	6.752474
0.767 	6.771654
0.768 	6.790896
0.769 	6.810200
0.770 	6.829566
0.771 	6.848995
0.772 	6.868487
0.773 	6.888041
0.774 	6.907659
0.775 	6.927341
0.776 	6.947086
0.777 	6.966895
0.778 	6.986768
0.779 	7.006705
0.780 	7.026707
0.781 	7.046774
0.782 	7.066905
0.783 	7.087102
0.784 	7.107364
0.785 	7.127692
0.786 	7.148085
0.787 	7.168545
0.788 	7.189071
0.789 	7.209663
0.790 	7.230322
0.791 	7.251048
0.792 	7.271841
0.793 	7.292701
0.794 	7.313629
0.795 	7.334625
0.796 	7.355689
0.797 	7.376821
0.798 	7.398022
0.799 	7.419291
0.800 	7.440630
0.801 	7.462037
0.802 	7.483514
0.803 	7.505061
0.804 	7.526677
0.805 	7.548364
0.806 	7.570121
0.807 	7.591949
0.808 	7.613847
0.809 	7.635817
0.810 	7.657858
0.811 	7.679970
0.812 	7.702154
0.813 	7.724411
0.814 	7.746739
0.815 	7.769140
0.816 	7.791614
0.817 	7.814161
0.818 	7.836781
0.819 	7.859474
0.820 	7.882241
0.821 	7.905083
0.822 	7.927998
0.823 	7.950988
0.824 	7.974053
0.825 	7.997193
0.826 	8.020408
0.827 	8.043698
0.828 	8.067064
0.829 	8.090506
0.830 	8.114025
0.831 	8.137620
0.832 	8.161291
0.833 	8.185040
0.834 	8.208866
0.835 	8.232770
0.836 	8.256751
0.837 	8.280811
0.838 	8.304948
0.839 	8.329165
0.840 	8.353460
0.841 	8.377834
0.842 	8.402288
0.843 	8.426821
0.844 	8.451435
0.845 	8.476128
0.846 	8.500902
0.847 	8.525757
0.848 	8.550693
0.849 	8.575710
0.850 	8.600808
0.851 	8.625989
0.852 	8.651251
0.853 	8.676596
0.854 	8.702024
0.855 	8.727534
0.856 	8.753128
0.857 	8.778805
0.858 	8.804566
0.859 	8.830411
0.860 	8.856341
0.861 	8.882355
0.862 	8.908454
0.863 	8.934638
0.864 	8.960907
0.865 	8.987263
0.866 	9.013704
0.867 	9.040232
0.868 	9.066847
0.869 	9.093548
0.870 	9.120337
0.871 	9.147213
0.872 	9.174177
0.873 	9.201230
0.874 	9.228370
0.875 	9.255599
0.876 	9.282918
0.877 	9.310326
0.878 	9.337823
0.879 	9.365410
0.880 	9.393088
0.881 	9.420856
0.882 	9.448714
0.883 	9.476664
0.884 	9.504706
0.885 	9.532839
0.886 	9.561064
0.887 	9.589382
0.888 	9.617793
0.889 	9.646296
0.890 	9.674893
0.891 	9.703583
0.892 	9.732368
0.893 	9.761246
0.894 	9.790219
0.895 	9.819288
0.896 	9.848451
0.897 	9.877710
0.898 	9.907065
0.899 	9.936516
0.900 	9.966064
0.901 	9.995708
0.902 	10.025450
0.903 	10.055289
0.904 	10.085226
0.905 	10.115262
0.906 	10.145396
0.907 	10.175628
0.908 	10.205960
0.909 	10.236392
0.910 	10.266923
0.911 	10.297554
0.912 	10.328286
0.913 	10.359119
0.914 	10.390053
0.915 	10.421089
0.916 	10.452227
0.917 	10.483467
0.918 	10.514809
0.919 	10.546255
0.920 	10.577803
0.921 	10.609456
0.922 	10.641212
0.923 	10.673073
0.924 	10.705039
0.925 	10.737109
0.926 	10.769285
0.927 	10.801567
0.928 	10.833955
0.929 	10.866450
0.930 	10.899051
0.931 	10.931760
0.932 	10.964576
0.933 	10.997500
0.934 	11.030533
0.935 	11.063674
0.936 	11.096924
0.937 	11.130284
0.938 	11.163754
0.939 	11.197334
0.940 	11.231024
0.941 	11.264825
0.942 	11.298738
0.943 	11.332762
0.944 	11.366899
0.945 	11.401148
0.946 	11.435509
0.947 	11.469984
0.948 	11.504573
0.949 	11.539275
0.950 	11.574092
/

\plot 0.21 11.0 0.3 11.0 /
\put {$  a(t)$}  [l] at 0.32  11

\setdashes
\plot
0.000 	1.490240
0.001 	1.490244
0.002 	1.490256
0.003 	1.490276
0.004 	1.490303
0.005 	1.490339
0.006 	1.490383
0.007 	1.490434
0.008 	1.490494
0.009 	1.490562
0.010 	1.490637
0.011 	1.490721
0.012 	1.490812
0.013 	1.490911
0.014 	1.491019
0.015 	1.491134
0.016 	1.491257
0.017 	1.491388
0.018 	1.491528
0.019 	1.491675
0.020 	1.491830
0.021 	1.491993
0.022 	1.492164
0.023 	1.492343
0.024 	1.492530
0.025 	1.492725
0.026 	1.492928
0.027 	1.493138
0.028 	1.493357
0.029 	1.493584
0.030 	1.493819
0.031 	1.494062
0.032 	1.494312
0.033 	1.494571
0.034 	1.494838
0.035 	1.495113
0.036 	1.495395
0.037 	1.495686
0.038 	1.495985
0.039 	1.496292
0.040 	1.496606
0.041 	1.496929
0.042 	1.497260
0.043 	1.497598
0.044 	1.497945
0.045 	1.498300
0.046 	1.498663
0.047 	1.499034
0.048 	1.499413
0.049 	1.499800
0.050 	1.500195
0.051 	1.500598
0.052 	1.501009
0.053 	1.501428
0.054 	1.501855
0.055 	1.502290
0.056 	1.502733
0.057 	1.503185
0.058 	1.503644
0.059 	1.504112
0.060 	1.504587
0.061 	1.505071
0.062 	1.505563
0.063 	1.506063
0.064 	1.506571
0.065 	1.507087
0.066 	1.507611
0.067 	1.508143
0.068 	1.508683
0.069 	1.509232
0.070 	1.509789
0.071 	1.510353
0.072 	1.510926
0.073 	1.511508
0.074 	1.512097
0.075 	1.512694
0.076 	1.513300
0.077 	1.513913
0.078 	1.514535
0.079 	1.515165
0.080 	1.515804
0.081 	1.516450
0.082 	1.517105
0.083 	1.517768
0.084 	1.518439
0.085 	1.519118
0.086 	1.519806
0.087 	1.520501
0.088 	1.521205
0.089 	1.521918
0.090 	1.522638
0.091 	1.523367
0.092 	1.524104
0.093 	1.524849
0.094 	1.525603
0.095 	1.526365
0.096 	1.527135
0.097 	1.527913
0.098 	1.528700
0.099 	1.529495
0.100 	1.530299
0.101 	1.531110
0.102 	1.531930
0.103 	1.532759
0.104 	1.533596
0.105 	1.534441
0.106 	1.535294
0.107 	1.536156
0.108 	1.537027
0.109 	1.537905
0.110 	1.538792
0.111 	1.539688
0.112 	1.540592
0.113 	1.541504
0.114 	1.542425
0.115 	1.543354
0.116 	1.544292
0.117 	1.545238
0.118 	1.546193
0.119 	1.547156
0.120 	1.548128
0.121 	1.549108
0.122 	1.550097
0.123 	1.551094
0.124 	1.552100
0.125 	1.553114
0.126 	1.554137
0.127 	1.555168
0.128 	1.556208
0.129 	1.557257
0.130 	1.558314
0.131 	1.559380
0.132 	1.560454
0.133 	1.561537
0.134 	1.562628
0.135 	1.563729
0.136 	1.564838
0.137 	1.565955
0.138 	1.567081
0.139 	1.568216
0.140 	1.569360
0.141 	1.570512
0.142 	1.571674
0.143 	1.572843
0.144 	1.574022
0.145 	1.575209
0.146 	1.576405
0.147 	1.577610
0.148 	1.578824
0.149 	1.580046
0.150 	1.581277
0.151 	1.582518
0.152 	1.583766
0.153 	1.585024
0.154 	1.586291
0.155 	1.587566
0.156 	1.588851
0.157 	1.590144
0.158 	1.591446
0.159 	1.592757
0.160 	1.594077
0.161 	1.595406
0.162 	1.596744
0.163 	1.598091
0.164 	1.599447
0.165 	1.600812
0.166 	1.602186
0.167 	1.603569
0.168 	1.604961
0.169 	1.606362
0.170 	1.607772
0.171 	1.609191
0.172 	1.610620
0.173 	1.612057
0.174 	1.613504
0.175 	1.614959
0.176 	1.616424
0.177 	1.617898
0.178 	1.619381
0.179 	1.620874
0.180 	1.622375
0.181 	1.623886
0.182 	1.625406
0.183 	1.626935
0.184 	1.628474
0.185 	1.630022
0.186 	1.631579
0.187 	1.633145
0.188 	1.634721
0.189 	1.636306
0.190 	1.637901
0.191 	1.639505
0.192 	1.641118
0.193 	1.642740
0.194 	1.644372
0.195 	1.646014
0.196 	1.647665
0.197 	1.649325
0.198 	1.650995
0.199 	1.652675
0.200 	1.654364
0.201 	1.656062
0.202 	1.657770
0.203 	1.659488
0.204 	1.661215
0.205 	1.662951
0.206 	1.664698
0.207 	1.666454
0.208 	1.668219
0.209 	1.669995
0.210 	1.671780
0.211 	1.673574
0.212 	1.675379
0.213 	1.677193
0.214 	1.679017
0.215 	1.680851
0.216 	1.682694
0.217 	1.684547
0.218 	1.686410
0.219 	1.688283
0.220 	1.690166
0.221 	1.692059
0.222 	1.693962
0.223 	1.695874
0.224 	1.697797
0.225 	1.699729
0.226 	1.701671
0.227 	1.703624
0.228 	1.705586
0.229 	1.707558
0.230 	1.709541
0.231 	1.711533
0.232 	1.713536
0.233 	1.715549
0.234 	1.717572
0.235 	1.719605
0.236 	1.721648
0.237 	1.723701
0.238 	1.725765
0.239 	1.727839
0.240 	1.729923
0.241 	1.732017
0.242 	1.734122
0.243 	1.736237
0.244 	1.738362
0.245 	1.740497
0.246 	1.742643
0.247 	1.744800
0.248 	1.746967
0.249 	1.749144
0.250 	1.751332
0.251 	1.753530
0.252 	1.755738
0.253 	1.757958
0.254 	1.760187
0.255 	1.762428
0.256 	1.764678
0.257 	1.766940
0.258 	1.769212
0.259 	1.771495
0.260 	1.773788
0.261 	1.776092
0.262 	1.778407
0.263 	1.780733
0.264 	1.783069
0.265 	1.785416
0.266 	1.787774
0.267 	1.790143
0.268 	1.792523
0.269 	1.794913
0.270 	1.797315
0.271 	1.799727
0.272 	1.802150
0.273 	1.804585
0.274 	1.807030
0.275 	1.809486
0.276 	1.811954
0.277 	1.814432
0.278 	1.816922
0.279 	1.819422
0.280 	1.821934
0.281 	1.824457
0.282 	1.826991
0.283 	1.829537
0.284 	1.832093
0.285 	1.834661
0.286 	1.837240
0.287 	1.839831
0.288 	1.842433
0.289 	1.845046
0.290 	1.847671
0.291 	1.850307
0.292 	1.852954
0.293 	1.855613
0.294 	1.858284
0.295 	1.860966
0.296 	1.863659
0.297 	1.866364
0.298 	1.869081
0.299 	1.871809
0.300 	1.874549
0.301 	1.877301
0.302 	1.880065
0.303 	1.882840
0.304 	1.885627
0.305 	1.888425
0.306 	1.891236
0.307 	1.894058
0.308 	1.896893
0.309 	1.899739
0.310 	1.902597
0.311 	1.905467
0.312 	1.908349
0.313 	1.911243
0.314 	1.914150
0.315 	1.917068
0.316 	1.919998
0.317 	1.922941
0.318 	1.925895
0.319 	1.928862
0.320 	1.931841
0.321 	1.934832
0.322 	1.937836
0.323 	1.940852
0.324 	1.943880
0.325 	1.946921
0.326 	1.949974
0.327 	1.953039
0.328 	1.956117
0.329 	1.959208
0.330 	1.962310
0.331 	1.965426
0.332 	1.968554
0.333 	1.971695
0.334 	1.974848
0.335 	1.978014
0.336 	1.981193
0.337 	1.984384
0.338 	1.987588
0.339 	1.990805
0.340 	1.994035
0.341 	1.997278
0.342 	2.000534
0.343 	2.003802
0.344 	2.007084
0.345 	2.010378
0.346 	2.013686
0.347 	2.017006
0.348 	2.020340
0.349 	2.023687
0.350 	2.027047
0.351 	2.030420
0.352 	2.033807
0.353 	2.037206
0.354 	2.040619
0.355 	2.044046
0.356 	2.047485
0.357 	2.050938
0.358 	2.054405
0.359 	2.057885
0.360 	2.061378
0.361 	2.064885
0.362 	2.068406
0.363 	2.071940
0.364 	2.075488
0.365 	2.079049
0.366 	2.082625
0.367 	2.086214
0.368 	2.089816
0.369 	2.093433
0.370 	2.097063
0.371 	2.100708
0.372 	2.104366
0.373 	2.108038
0.374 	2.111724
0.375 	2.115425
0.376 	2.119139
0.377 	2.122867
0.378 	2.126610
0.379 	2.130367
0.380 	2.134138
0.381 	2.137923
0.382 	2.141722
0.383 	2.145536
0.384 	2.149365
0.385 	2.153207
0.386 	2.157064
0.387 	2.160936
0.388 	2.164822
0.389 	2.168723
0.390 	2.172638
0.391 	2.176568
0.392 	2.180513
0.393 	2.184472
0.394 	2.188446
0.395 	2.192435
0.396 	2.196439
0.397 	2.200458
0.398 	2.204491
0.399 	2.208540
0.400 	2.212604
0.401 	2.216682
0.402 	2.220776
0.403 	2.224885
0.404 	2.229009
0.405 	2.233148
0.406 	2.237302
0.407 	2.241472
0.408 	2.245657
0.409 	2.249858
0.410 	2.254074
0.411 	2.258305
0.412 	2.262552
0.413 	2.266814
0.414 	2.271092
0.415 	2.275386
0.416 	2.279695
0.417 	2.284020
0.418 	2.288361
0.419 	2.292718
0.420 	2.297090
0.421 	2.301478
0.422 	2.305883
0.423 	2.310303
0.424 	2.314739
0.425 	2.319192
0.426 	2.323660
0.427 	2.328145
0.428 	2.332646
0.429 	2.337163
0.430 	2.341696
0.431 	2.346246
0.432 	2.350812
0.433 	2.355395
0.434 	2.359994
0.435 	2.364609
0.436 	2.369241
0.437 	2.373890
0.438 	2.378556
0.439 	2.383238
0.440 	2.387937
0.441 	2.392652
0.442 	2.397385
0.443 	2.402134
0.444 	2.406901
0.445 	2.411684
0.446 	2.416485
0.447 	2.421302
0.448 	2.426137
0.449 	2.430989
0.450 	2.435858
0.451 	2.440745
0.452 	2.445649
0.453 	2.450570
0.454 	2.455509
0.455 	2.460465
0.456 	2.465438
0.457 	2.470430
0.458 	2.475439
0.459 	2.480465
0.460 	2.485510
0.461 	2.490572
0.462 	2.495652
0.463 	2.500750
0.464 	2.505866
0.465 	2.511000
0.466 	2.516151
0.467 	2.521322
0.468 	2.526510
0.469 	2.531716
0.470 	2.536941
0.471 	2.542184
0.472 	2.547445
0.473 	2.552725
0.474 	2.558023
0.475 	2.563340
0.476 	2.568676
0.477 	2.574030
0.478 	2.579402
0.479 	2.584794
0.480 	2.590204
0.481 	2.595633
0.482 	2.601081
0.483 	2.606548
0.484 	2.612035
0.485 	2.617540
0.486 	2.623064
0.487 	2.628607
0.488 	2.634170
0.489 	2.639752
0.490 	2.645354
0.491 	2.650975
0.492 	2.656615
0.493 	2.662275
0.494 	2.667954
0.495 	2.673653
0.496 	2.679372
0.497 	2.685111
0.498 	2.690869
0.499 	2.696648
0.500 	2.702446
0.501 	2.708264
0.502 	2.714102
0.503 	2.719961
0.504 	2.725839
0.505 	2.731738
0.506 	2.737658
0.507 	2.743597
0.508 	2.749557
0.509 	2.755538
0.510 	2.761539
0.511 	2.767560
0.512 	2.773603
0.513 	2.779666
0.514 	2.785749
0.515 	2.791854
0.516 	2.797980
0.517 	2.804126
0.518 	2.810294
0.519 	2.816483
0.520 	2.822693
0.521 	2.828924
0.522 	2.835176
0.523 	2.841450
0.524 	2.847745
0.525 	2.854062
0.526 	2.860401
0.527 	2.866761
0.528 	2.873142
0.529 	2.879546
0.530 	2.885971
0.531 	2.892418
0.532 	2.898887
0.533 	2.905378
0.534 	2.911891
0.535 	2.918427
0.536 	2.924984
0.537 	2.931564
0.538 	2.938166
0.539 	2.944791
0.540 	2.951438
0.541 	2.958108
0.542 	2.964800
0.543 	2.971515
0.544 	2.978253
0.545 	2.985013
0.546 	2.991797
0.547 	2.998603
0.548 	3.005433
0.549 	3.012286
0.550 	3.019161
0.551 	3.026061
0.552 	3.032983
0.553 	3.039929
0.554 	3.046898
0.555 	3.053891
0.556 	3.060908
0.557 	3.067948
0.558 	3.075012
0.559 	3.082100
0.560 	3.089211
0.561 	3.096347
0.562 	3.103507
0.563 	3.110691
0.564 	3.117899
0.565 	3.125131
0.566 	3.132388
0.567 	3.139669
0.568 	3.146975
0.569 	3.154305
0.570 	3.161660
0.571 	3.169040
0.572 	3.176444
0.573 	3.183874
0.574 	3.191328
0.575 	3.198807
0.576 	3.206312
0.577 	3.213841
0.578 	3.221396
0.579 	3.228977
0.580 	3.236583
0.581 	3.244214
0.582 	3.251871
0.583 	3.259553
0.584 	3.267262
0.585 	3.274996
0.586 	3.282756
0.587 	3.290542
0.588 	3.298354
0.589 	3.306192
0.590 	3.314057
0.591 	3.321948
0.592 	3.329865
0.593 	3.337809
0.594 	3.345779
0.595 	3.353776
0.596 	3.361800
0.597 	3.369850
0.598 	3.377927
0.599 	3.386032
0.600 	3.394163
0.601 	3.402322
0.602 	3.410508
0.603 	3.418721
0.604 	3.426962
0.605 	3.435230
0.606 	3.443525
0.607 	3.451849
0.608 	3.460200
0.609 	3.468579
0.610 	3.476986
0.611 	3.485421
0.612 	3.493884
0.613 	3.502375
0.614 	3.510895
0.615 	3.519442
0.616 	3.528019
0.617 	3.536624
0.618 	3.545257
0.619 	3.553920
0.620 	3.562611
0.621 	3.571331
0.622 	3.580080
0.623 	3.588858
0.624 	3.597665
0.625 	3.606502
0.626 	3.615368
0.627 	3.624264
0.628 	3.633189
0.629 	3.642143
0.630 	3.651128
0.631 	3.660142
0.632 	3.669187
0.633 	3.678261
0.634 	3.687365
0.635 	3.696500
0.636 	3.705665
0.637 	3.714860
0.638 	3.724086
0.639 	3.733343
0.640 	3.742630
0.641 	3.751948
0.642 	3.761297
0.643 	3.770677
0.644 	3.780088
0.645 	3.789531
0.646 	3.799004
0.647 	3.808509
0.648 	3.818046
0.649 	3.827614
0.650 	3.837214
0.651 	3.846846
0.652 	3.856509
0.653 	3.866205
0.654 	3.875932
0.655 	3.885692
0.656 	3.895484
0.657 	3.905309
0.658 	3.915166
0.659 	3.925056
0.660 	3.934979
0.661 	3.944934
0.662 	3.954922
0.663 	3.964944
0.664 	3.974998
0.665 	3.985086
0.666 	3.995207
0.667 	4.005362
0.668 	4.015550
0.669 	4.025772
0.670 	4.036028
0.671 	4.046318
0.672 	4.056641
0.673 	4.066999
0.674 	4.077391
0.675 	4.087818
0.676 	4.098279
0.677 	4.108774
0.678 	4.119304
0.679 	4.129869
0.680 	4.140469
0.681 	4.151104
0.682 	4.161774
0.683 	4.172479
0.684 	4.183219
0.685 	4.193996
0.686 	4.204807
0.687 	4.215655
0.688 	4.226538
0.689 	4.237457
0.690 	4.248412
0.691 	4.259403
0.692 	4.270431
0.693 	4.281495
0.694 	4.292595
0.695 	4.303733
0.696 	4.314907
0.697 	4.326117
0.698 	4.337365
0.699 	4.348650
0.700 	4.359972
0.701 	4.371332
0.702 	4.382729
0.703 	4.394163
0.704 	4.405636
0.705 	4.417146
0.706 	4.428694
0.707 	4.440280
0.708 	4.451905
0.709 	4.463567
0.710 	4.475269
0.711 	4.487008
0.712 	4.498787
0.713 	4.510604
0.714 	4.522461
0.715 	4.534356
0.716 	4.546291
0.717 	4.558265
0.718 	4.570278
0.719 	4.582331
0.720 	4.594424
0.721 	4.606557
0.722 	4.618730
0.723 	4.630942
0.724 	4.643196
0.725 	4.655489
0.726 	4.667823
0.727 	4.680198
0.728 	4.692613
0.729 	4.705069
0.730 	4.717567
0.731 	4.730105
0.732 	4.742685
0.733 	4.755306
0.734 	4.767969
0.735 	4.780674
0.736 	4.793421
0.737 	4.806209
0.738 	4.819040
0.739 	4.831913
0.740 	4.844828
0.741 	4.857786
0.742 	4.870787
0.743 	4.883830
0.744 	4.896916
0.745 	4.910046
0.746 	4.923219
0.747 	4.936435
0.748 	4.949695
0.749 	4.962998
0.750 	4.976345
0.751 	4.989736
0.752 	5.003172
0.753 	5.016651
0.754 	5.030175
0.755 	5.043744
0.756 	5.057357
0.757 	5.071015
0.758 	5.084718
0.759 	5.098466
0.760 	5.112260
0.761 	5.126099
0.762 	5.139983
0.763 	5.153913
0.764 	5.167889
0.765 	5.181912
0.766 	5.195980
0.767 	5.210095
0.768 	5.224256
0.769 	5.238464
0.770 	5.252718
0.771 	5.267020
0.772 	5.281369
0.773 	5.295765
0.774 	5.310208
0.775 	5.324699
0.776 	5.339238
0.777 	5.353825
0.778 	5.368459
0.779 	5.383142
0.780 	5.397873
0.781 	5.412653
0.782 	5.427482
0.783 	5.442359
0.784 	5.457285
0.785 	5.472261
0.786 	5.487286
0.787 	5.502360
0.788 	5.517484
0.789 	5.532658
0.790 	5.547882
0.791 	5.563156
0.792 	5.578480
0.793 	5.593855
0.794 	5.609281
0.795 	5.624757
0.796 	5.640284
0.797 	5.655863
0.798 	5.671493
0.799 	5.687174
0.800 	5.702907
0.801 	5.718692
0.802 	5.734529
0.803 	5.750418
0.804 	5.766360
0.805 	5.782354
0.806 	5.798400
0.807 	5.814500
0.808 	5.830653
0.809 	5.846858
0.810 	5.863118
0.811 	5.879431
0.812 	5.895797
0.813 	5.912218
0.814 	5.928693
0.815 	5.945222
0.816 	5.961806
0.817 	5.978444
0.818 	5.995137
0.819 	6.011885
0.820 	6.028689
0.821 	6.045547
0.822 	6.062462
0.823 	6.079432
0.824 	6.096458
0.825 	6.113541
0.826 	6.130679
0.827 	6.147874
0.828 	6.165126
0.829 	6.182435
0.830 	6.199801
0.831 	6.217224
0.832 	6.234705
0.833 	6.252243
0.834 	6.269839
0.835 	6.287493
0.836 	6.305206
0.837 	6.322977
0.838 	6.340806
0.839 	6.358694
0.840 	6.376642
0.841 	6.394648
0.842 	6.412714
0.843 	6.430840
0.844 	6.449025
0.845 	6.467271
0.846 	6.485577
0.847 	6.503943
0.848 	6.522370
0.849 	6.540857
0.850 	6.559406
0.851 	6.578016
0.852 	6.596687
0.853 	6.615420
0.854 	6.634214
0.855 	6.653071
0.856 	6.671990
0.857 	6.690972
0.858 	6.710016
0.859 	6.729123
0.860 	6.748294
0.861 	6.767527
0.862 	6.786824
0.863 	6.806185
0.864 	6.825610
0.865 	6.845099
0.866 	6.864652
0.867 	6.884270
0.868 	6.903953
0.869 	6.923701
0.870 	6.943514
0.871 	6.963393
0.872 	6.983337
0.873 	7.003347
0.874 	7.023423
0.875 	7.043566
0.876 	7.063775
0.877 	7.084051
0.878 	7.104394
0.879 	7.124805
0.880 	7.145282
0.881 	7.165828
0.882 	7.186441
0.883 	7.207123
0.884 	7.227873
0.885 	7.248692
0.886 	7.269579
0.887 	7.290536
0.888 	7.311561
0.889 	7.332657
0.890 	7.353822
0.891 	7.375057
0.892 	7.396363
0.893 	7.417739
0.894 	7.439186
0.895 	7.460703
0.896 	7.482292
0.897 	7.503952
0.898 	7.525684
0.899 	7.547488
0.900 	7.569365
0.901 	7.591313
0.902 	7.613334
0.903 	7.635428
0.904 	7.657596
0.905 	7.679836
0.906 	7.702151
0.907 	7.724539
0.908 	7.747001
0.909 	7.769538
0.910 	7.792149
0.911 	7.814836
0.912 	7.837597
0.913 	7.860434
0.914 	7.883346
0.915 	7.906335
0.916 	7.929399
0.917 	7.952540
0.918 	7.975758
0.919 	7.999053
0.920 	8.022424
0.921 	8.045874
0.922 	8.069401
0.923 	8.093006
0.924 	8.116689
0.925 	8.140451
0.926 	8.164291
0.927 	8.188211
0.928 	8.212210
0.929 	8.236289
0.930 	8.260447
0.931 	8.284686
0.932 	8.309004
0.933 	8.333404
0.934 	8.357885
0.935 	8.382446
0.936 	8.407090
0.937 	8.431815
0.938 	8.456622
0.939 	8.481511
0.940 	8.506483
0.941 	8.531538
0.942 	8.556676
0.943 	8.581898
0.944 	8.607203
0.945 	8.632592
0.946 	8.658066
0.947 	8.683624
0.948 	8.709267
0.949 	8.734995
0.950 	8.760808
0.951 	8.786708
0.952 	8.812693
0.953 	8.838765
0.954 	8.864923
0.955 	8.891168
0.956 	8.917501
0.957 	8.943921
0.958 	8.970428
0.959 	8.997024
0.960 	9.023708
0.961 	9.050481
0.962 	9.077343
0.963 	9.104294
0.964 	9.131335
0.965 	9.158466
0.966 	9.185687
0.967 	9.212998
0.968 	9.240401
0.969 	9.267894
0.970 	9.295479
0.971 	9.323155
0.972 	9.350924
0.973 	9.378785
0.974 	9.406739
0.975 	9.434786
0.976 	9.462926
0.977 	9.491159
0.978 	9.519487
0.979 	9.547909
0.980 	9.576426
0.981 	9.605037
0.982 	9.633744
0.983 	9.662546
0.984 	9.691445
0.985 	9.720439
0.986 	9.749530
0.987 	9.778718
0.988 	9.808003
0.989 	9.837386
0.990 	9.866867
0.991 	9.896445
0.992 	9.926123
0.993 	9.955899
0.994 	9.985774
0.995 	10.015749
0.996 	10.045824
0.997 	10.075999
0.998 	10.106275
0.999 	10.136652
/
\plot 0.21 9.0 0.3 9.0 /
\put {$  a_0 \sinh(4t)/(4t)$}  [l] at 0.32  9

\setdots
\plot 
0.0000	1.49024
0.0010	1.49025
0.0020	1.49026
0.0030	1.49029
0.0040	1.49034
0.0050	1.49039
0.0060	1.49045
0.0070	1.49053
0.0080	1.49062
0.0090	1.49072
0.0100	1.49084
0.0110	1.49096
0.0120	1.49110
0.0130	1.49125
0.0140	1.49141
0.0150	1.49158
0.0160	1.49177
0.0170	1.49196
0.0180	1.49217
0.0190	1.49239
0.0200	1.49263
0.0210	1.49287
0.0220	1.49313
0.0230	1.49339
0.0240	1.49367
0.0250	1.49397
0.0260	1.49427
0.0270	1.49459
0.0280	1.49492
0.0290	1.49526
0.0300	1.49561
0.0310	1.49597
0.0320	1.49635
0.0330	1.49674
0.0340	1.49714
0.0350	1.49755
0.0360	1.49797
0.0370	1.49841
0.0380	1.49886
0.0390	1.49932
0.0400	1.49979
0.0410	1.50027
0.0420	1.50077
0.0430	1.50128
0.0440	1.50180
0.0450	1.50233
0.0460	1.50287
0.0470	1.50343
0.0480	1.50400
0.0490	1.50458
0.0500	1.50517
0.0510	1.50577
0.0520	1.50639
0.0530	1.50702
0.0540	1.50766
0.0550	1.50831
0.0560	1.50898
0.0570	1.50965
0.0580	1.51034
0.0590	1.51104
0.0600	1.51176
0.0610	1.51248
0.0620	1.51322
0.0630	1.51397
0.0640	1.51473
0.0650	1.51550
0.0660	1.51629
0.0670	1.51709
0.0680	1.51790
0.0690	1.51872
0.0700	1.51956
0.0710	1.52040
0.0720	1.52126
0.0730	1.52213
0.0740	1.52302
0.0750	1.52391
0.0760	1.52482
0.0770	1.52574
0.0780	1.52667
0.0790	1.52762
0.0800	1.52857
0.0810	1.52954
0.0820	1.53052
0.0830	1.53152
0.0840	1.53252
0.0850	1.53354
0.0860	1.53457
0.0870	1.53561
0.0880	1.53667
0.0890	1.53774
0.0900	1.53882
0.0910	1.53991
0.0920	1.54101
0.0930	1.54213
0.0940	1.54326
0.0950	1.54440
0.0960	1.54556
0.0970	1.54672
0.0980	1.54790
0.0990	1.54909
0.1000	1.55030
0.1010	1.55151
0.1020	1.55274
0.1030	1.55398
0.1040	1.55524
0.1050	1.55650
0.1060	1.55778
0.1070	1.55907
0.1080	1.56038
0.1090	1.56169
0.1100	1.56302
0.1110	1.56436
0.1120	1.56572
0.1130	1.56708
0.1140	1.56846
0.1150	1.56985
0.1160	1.57126
0.1170	1.57268
0.1180	1.57411
0.1190	1.57555
0.1200	1.57700
0.1210	1.57847
0.1220	1.57995
0.1230	1.58145
0.1240	1.58295
0.1250	1.58447
0.1260	1.58600
0.1270	1.58755
0.1280	1.58910
0.1290	1.59067
0.1300	1.59226
0.1310	1.59385
0.1320	1.59546
0.1330	1.59708
0.1340	1.59872
0.1350	1.60036
0.1360	1.60202
0.1370	1.60370
0.1380	1.60538
0.1390	1.60708
0.1400	1.60879
0.1410	1.61052
0.1420	1.61226
0.1430	1.61401
0.1440	1.61577
0.1450	1.61755
0.1460	1.61934
0.1470	1.62114
0.1480	1.62296
0.1490	1.62479
0.1500	1.62663
0.1510	1.62848
0.1520	1.63035
0.1530	1.63223
0.1540	1.63413
0.1550	1.63604
0.1560	1.63796
0.1570	1.63989
0.1580	1.64184
0.1590	1.64380
0.1600	1.64578
0.1610	1.64777
0.1620	1.64977
0.1630	1.65178
0.1640	1.65381
0.1650	1.65585
0.1660	1.65791
0.1670	1.65997
0.1680	1.66206
0.1690	1.66415
0.1700	1.66626
0.1710	1.66838
0.1720	1.67052
0.1730	1.67267
0.1740	1.67483
0.1750	1.67700
0.1760	1.67919
0.1770	1.68140
0.1780	1.68361
0.1790	1.68584
0.1800	1.68809
0.1810	1.69035
0.1820	1.69262
0.1830	1.69490
0.1840	1.69720
0.1850	1.69952
0.1860	1.70184
0.1870	1.70418
0.1880	1.70654
0.1890	1.70891
0.1900	1.71129
0.1910	1.71368
0.1920	1.71609
0.1930	1.71852
0.1940	1.72096
0.1950	1.72341
0.1960	1.72587
0.1970	1.72835
0.1980	1.73085
0.1990	1.73336
0.2000	1.73588
0.2010	1.73841
0.2020	1.74096
0.2030	1.74353
0.2040	1.74611
0.2050	1.74870
0.2060	1.75131
0.2070	1.75393
0.2080	1.75656
0.2090	1.75921
0.2100	1.76188
0.2110	1.76456
0.2120	1.76725
0.2130	1.76996
0.2140	1.77268
0.2150	1.77542
0.2160	1.77817
0.2170	1.78093
0.2180	1.78371
0.2190	1.78651
0.2200	1.78932
0.2210	1.79214
0.2220	1.79498
0.2230	1.79783
0.2240	1.80070
0.2250	1.80358
0.2260	1.80648
0.2270	1.80939
0.2280	1.81231
0.2290	1.81526
0.2300	1.81821
0.2310	1.82118
0.2320	1.82417
0.2330	1.82717
0.2340	1.83018
0.2350	1.83321
0.2360	1.83626
0.2370	1.83932
0.2380	1.84240
0.2390	1.84549
0.2400	1.84859
0.2410	1.85171
0.2420	1.85485
0.2430	1.85800
0.2440	1.86117
0.2450	1.86435
0.2460	1.86755
0.2470	1.87076
0.2480	1.87398
0.2490	1.87723
0.2500	1.88049
0.2510	1.88376
0.2520	1.88705
0.2530	1.89035
0.2540	1.89367
0.2550	1.89701
0.2560	1.90036
0.2570	1.90373
0.2580	1.90711
0.2590	1.91051
0.2600	1.91392
0.2610	1.91735
0.2620	1.92079
0.2630	1.92426
0.2640	1.92773
0.2650	1.93122
0.2660	1.93473
0.2670	1.93826
0.2680	1.94180
0.2690	1.94535
0.2700	1.94892
0.2710	1.95251
0.2720	1.95612
0.2730	1.95974
0.2740	1.96337
0.2750	1.96703
0.2760	1.97069
0.2770	1.97438
0.2780	1.97808
0.2790	1.98180
0.2800	1.98553
0.2810	1.98928
0.2820	1.99305
0.2830	1.99683
0.2840	2.00063
0.2850	2.00444
0.2860	2.00828
0.2870	2.01213
0.2880	2.01599
0.2890	2.01987
0.2900	2.02377
0.2910	2.02769
0.2920	2.03162
0.2930	2.03557
0.2940	2.03953
0.2950	2.04351
0.2960	2.04751
0.2970	2.05153
0.2980	2.05556
0.2990	2.05961
0.3000	2.06368
0.3010	2.06776
0.3020	2.07186
0.3030	2.07598
0.3040	2.08012
0.3050	2.08427
0.3060	2.08844
0.3070	2.09263
0.3080	2.09683
0.3090	2.10105
0.3100	2.10529
0.3110	2.10955
0.3120	2.11382
0.3130	2.11811
0.3140	2.12242
0.3150	2.12675
0.3160	2.13109
0.3170	2.13545
0.3180	2.13983
0.3190	2.14423
0.3200	2.14864
0.3210	2.15308
0.3220	2.15753
0.3230	2.16199
0.3240	2.16648
0.3250	2.17098
0.3260	2.17551
0.3270	2.18005
0.3280	2.18460
0.3290	2.18918
0.3300	2.19377
0.3310	2.19839
0.3320	2.20302
0.3330	2.20767
0.3340	2.21233
0.3350	2.21702
0.3360	2.22172
0.3370	2.22645
0.3380	2.23119
0.3390	2.23595
0.3400	2.24073
0.3410	2.24552
0.3420	2.25034
0.3430	2.25517
0.3440	2.26002
0.3450	2.26490
0.3460	2.26979
0.3470	2.27469
0.3480	2.27962
0.3490	2.28457
0.3500	2.28954
0.3510	2.29452
0.3520	2.29953
0.3530	2.30455
0.3540	2.30959
0.3550	2.31465
0.3560	2.31973
0.3570	2.32483
0.3580	2.32995
0.3590	2.33509
0.3600	2.34025
0.3610	2.34543
0.3620	2.35063
0.3630	2.35584
0.3640	2.36108
0.3650	2.36634
0.3660	2.37161
0.3670	2.37691
0.3680	2.38222
0.3690	2.38756
0.3700	2.39292
0.3710	2.39829
0.3720	2.40369
0.3730	2.40910
0.3740	2.41454
0.3750	2.41999
0.3760	2.42547
0.3770	2.43096
0.3780	2.43648
0.3790	2.44202
0.3800	2.44758
0.3810	2.45315
0.3820	2.45875
0.3830	2.46437
0.3840	2.47001
0.3850	2.47567
0.3860	2.48135
0.3870	2.48705
0.3880	2.49277
0.3890	2.49851
0.3900	2.50428
0.3910	2.51006
0.3920	2.51587
0.3930	2.52170
0.3940	2.52754
0.3950	2.53341
0.3960	2.53930
0.3970	2.54521
0.3980	2.55115
0.3990	2.55710
0.4000	2.56308
0.4010	2.56907
0.4020	2.57509
0.4030	2.58113
0.4040	2.58719
0.4050	2.59328
0.4060	2.59938
0.4070	2.60551
0.4080	2.61166
0.4090	2.61783
0.4100	2.62402
0.4110	2.63024
0.4120	2.63647
0.4130	2.64273
0.4140	2.64901
0.4150	2.65531
0.4160	2.66164
0.4170	2.66799
0.4180	2.67436
0.4190	2.68075
0.4200	2.68717
0.4210	2.69360
0.4220	2.70006
0.4230	2.70655
0.4240	2.71305
0.4250	2.71958
0.4260	2.72613
0.4270	2.73271
0.4280	2.73930
0.4290	2.74592
0.4300	2.75257
0.4310	2.75923
0.4320	2.76592
0.4330	2.77264
0.4340	2.77937
0.4350	2.78613
0.4360	2.79291
0.4370	2.79972
0.4380	2.80655
0.4390	2.81340
0.4400	2.82028
0.4410	2.82718
0.4420	2.83411
0.4430	2.84106
0.4440	2.84803
0.4450	2.85503
0.4460	2.86205
0.4470	2.86909
0.4480	2.87616
0.4490	2.88325
0.4500	2.89037
0.4510	2.89751
0.4520	2.90468
0.4530	2.91187
0.4540	2.91908
0.4550	2.92632
0.4560	2.93359
0.4570	2.94088
0.4580	2.94819
0.4590	2.95553
0.4600	2.96289
0.4610	2.97028
0.4620	2.97770
0.4630	2.98514
0.4640	2.99260
0.4650	3.00009
0.4660	3.00761
0.4670	3.01515
0.4680	3.02271
0.4690	3.03030
0.4700	3.03792
0.4710	3.04556
0.4720	3.05323
0.4730	3.06093
0.4740	3.06865
0.4750	3.07639
0.4760	3.08416
0.4770	3.09196
0.4780	3.09978
0.4790	3.10764
0.4800	3.11551
0.4810	3.12341
0.4820	3.13134
0.4830	3.13930
0.4840	3.14728
0.4850	3.15529
0.4860	3.16333
0.4870	3.17139
0.4880	3.17948
0.4890	3.18760
0.4900	3.19574
0.4910	3.20391
0.4920	3.21211
0.4930	3.22033
0.4940	3.22858
0.4950	3.23686
0.4960	3.24517
0.4970	3.25350
0.4980	3.26187
0.4990	3.27026
0.5000	3.27867
0.5010	3.28712
0.5020	3.29559
0.5030	3.30409
0.5040	3.31262
0.5050	3.32118
0.5060	3.32976
0.5070	3.33838
0.5080	3.34702
0.5090	3.35569
0.5100	3.36439
0.5110	3.37311
0.5120	3.38187
0.5130	3.39065
0.5140	3.39947
0.5150	3.40831
0.5160	3.41718
0.5170	3.42608
0.5180	3.43501
0.5190	3.44397
0.5200	3.45295
0.5210	3.46197
0.5220	3.47102
0.5230	3.48009
0.5240	3.48920
0.5250	3.49833
0.5260	3.50750
0.5270	3.51669
0.5280	3.52592
0.5290	3.53517
0.5300	3.54446
0.5310	3.55377
0.5320	3.56312
0.5330	3.57249
0.5340	3.58190
0.5350	3.59133
0.5360	3.60080
0.5370	3.61030
0.5380	3.61983
0.5390	3.62938
0.5400	3.63897
0.5410	3.64860
0.5420	3.65825
0.5430	3.66793
0.5440	3.67764
0.5450	3.68739
0.5460	3.69717
0.5470	3.70698
0.5480	3.71682
0.5490	3.72669
0.5500	3.73659
0.5510	3.74653
0.5520	3.75649
0.5530	3.76649
0.5540	3.77652
0.5550	3.78659
0.5560	3.79668
0.5570	3.80681
0.5580	3.81697
0.5590	3.82717
0.5600	3.83739
0.5610	3.84765
0.5620	3.85794
0.5630	3.86827
0.5640	3.87863
0.5650	3.88902
0.5660	3.89944
0.5670	3.90990
0.5680	3.92039
0.5690	3.93091
0.5700	3.94147
0.5710	3.95206
0.5720	3.96269
0.5730	3.97335
0.5740	3.98404
0.5750	3.99477
0.5760	4.00553
0.5770	4.01632
0.5780	4.02715
0.5790	4.03802
0.5800	4.04891
0.5810	4.05985
0.5820	4.07082
0.5830	4.08182
0.5840	4.09286
0.5850	4.10393
0.5860	4.11504
0.5870	4.12618
0.5880	4.13736
0.5890	4.14857
0.5900	4.15982
0.5910	4.17111
0.5920	4.18243
0.5930	4.19378
0.5940	4.20517
0.5950	4.21660
0.5960	4.22807
0.5970	4.23957
0.5980	4.25110
0.5990	4.26268
0.6000	4.27429
0.6010	4.28593
0.6020	4.29762
0.6030	4.30934
0.6040	4.32109
0.6050	4.33289
0.6060	4.34472
0.6070	4.35659
0.6080	4.36849
0.6090	4.38044
0.6100	4.39242
0.6110	4.40444
0.6120	4.41649
0.6130	4.42859
0.6140	4.44072
0.6150	4.45289
0.6160	4.46510
0.6170	4.47734
0.6180	4.48963
0.6190	4.50195
0.6200	4.51432
0.6210	4.52672
0.6220	4.53916
0.6230	4.55164
0.6240	4.56415
0.6250	4.57671
0.6260	4.58931
0.6270	4.60194
0.6280	4.61462
0.6290	4.62733
0.6300	4.64009
0.6310	4.65288
0.6320	4.66572
0.6330	4.67859
0.6340	4.69151
0.6350	4.70446
0.6360	4.71746
0.6370	4.73050
0.6380	4.74357
0.6390	4.75669
0.6400	4.76985
0.6410	4.78305
0.6420	4.79629
0.6430	4.80958
0.6440	4.82290
0.6450	4.83626
0.6460	4.84967
0.6470	4.86312
0.6480	4.87661
0.6490	4.89014
0.6500	4.90372
0.6510	4.91733
0.6520	4.93099
0.6530	4.94470
0.6540	4.95844
0.6550	4.97223
0.6560	4.98606
0.6570	4.99993
0.6580	5.01385
0.6590	5.02780
0.6600	5.04181
0.6610	5.05585
0.6620	5.06994
0.6630	5.08408
0.6640	5.09825
0.6650	5.11247
0.6660	5.12674
0.6670	5.14105
0.6680	5.15540
0.6690	5.16980
0.6700	5.18424
0.6710	5.19873
0.6720	5.21326
0.6730	5.22784
0.6740	5.24246
0.6750	5.25713
0.6760	5.27184
0.6770	5.28660
0.6780	5.30141
0.6790	5.31626
0.6800	5.33115
0.6810	5.34610
0.6820	5.36108
0.6830	5.37612
0.6840	5.39120
0.6850	5.40633
0.6860	5.42150
0.6870	5.43672
0.6880	5.45199
0.6890	5.46731
0.6900	5.48267
0.6910	5.49808
0.6920	5.51354
0.6930	5.52904
0.6940	5.54459
0.6950	5.56020
0.6960	5.57584
0.6970	5.59154
0.6980	5.60729
0.6990	5.62308
0.7000	5.63892
0.7010	5.65482
0.7020	5.67076
0.7030	5.68675
0.7040	5.70279
0.7050	5.71887
0.7060	5.73501
0.7070	5.75120
0.7080	5.76744
0.7090	5.78372
0.7100	5.80006
0.7110	5.81645
0.7120	5.83289
0.7130	5.84938
0.7140	5.86592
0.7150	5.88251
0.7160	5.89915
0.7170	5.91584
0.7180	5.93258
0.7190	5.94938
0.7200	5.96623
0.7210	5.98312
0.7220	6.00007
0.7230	6.01708
0.7240	6.03413
0.7250	6.05124
0.7260	6.06840
0.7270	6.08561
0.7280	6.10288
0.7290	6.12019
0.7300	6.13757
0.7310	6.15499
0.7320	6.17247
0.7330	6.19000
0.7340	6.20758
0.7350	6.22522
0.7360	6.24292
0.7370	6.26066
0.7380	6.27847
0.7390	6.29632
0.7400	6.31423
0.7410	6.33220
0.7420	6.35022
0.7430	6.36830
0.7440	6.38643
0.7450	6.40462
0.7460	6.42286
0.7470	6.44116
0.7480	6.45951
0.7490	6.47792
0.7500	6.49639
0.7510	6.51491
0.7520	6.53350
0.7530	6.55213
0.7540	6.57083
0.7550	6.58958
0.7560	6.60839
0.7570	6.62725
0.7580	6.64618
0.7590	6.66516
0.7600	6.68420
0.7610	6.70330
0.7620	6.72246
0.7630	6.74167
0.7640	6.76095
0.7650	6.78028
0.7660	6.79967
0.7670	6.81912
0.7680	6.83863
0.7690	6.85820
0.7700	6.87783
0.7710	6.89753
0.7720	6.91728
0.7730	6.93709
0.7740	6.95696
0.7750	6.97689
0.7760	6.99688
0.7770	7.01694
0.7780	7.03705
0.7790	7.05723
0.7800	7.07747
0.7810	7.09777
0.7820	7.11813
0.7830	7.13855
0.7840	7.15904
0.7850	7.17959
0.7860	7.20020
0.7870	7.22088
0.7880	7.24161
0.7890	7.26241
0.7900	7.28328
0.7910	7.30421
0.7920	7.32520
0.7930	7.34625
0.7940	7.36737
0.7950	7.38856
0.7960	7.40981
0.7970	7.43112
0.7980	7.45250
0.7990	7.47395
0.8000	7.49546
0.8010	7.51703
0.8020	7.53867
0.8030	7.56038
0.8040	7.58215
0.8050	7.60399
0.8060	7.62590
0.8070	7.64787
0.8080	7.66991
0.8090	7.69202
0.8100	7.71419
0.8110	7.73644
0.8120	7.75875
0.8130	7.78112
0.8140	7.80357
0.8150	7.82608
0.8160	7.84867
0.8170	7.87132
0.8180	7.89404
0.8190	7.91683
0.8200	7.93969
0.8210	7.96262
0.8220	7.98562
0.8230	8.00869
0.8240	8.03183
0.8250	8.05504
0.8260	8.07832
0.8270	8.10167
0.8280	8.12509
0.8290	8.14859
0.8300	8.17215
0.8310	8.19579
0.8320	8.21950
0.8330	8.24328
0.8340	8.26713
0.8350	8.29106
0.8360	8.31506
0.8370	8.33913
0.8380	8.36328
0.8390	8.38750
0.8400	8.41179
0.8410	8.43616
0.8420	8.46060
0.8430	8.48511
0.8440	8.50970
0.8450	8.53437
0.8460	8.55911
0.8470	8.58392
0.8480	8.60882
0.8490	8.63378
0.8500	8.65883
0.8510	8.68394
0.8520	8.70914
0.8530	8.73441
0.8540	8.75976
0.8550	8.78519
0.8560	8.81069
0.8570	8.83627
0.8580	8.86193
0.8590	8.88767
0.8600	8.91348
0.8610	8.93938
0.8620	8.96535
0.8630	8.99140
0.8640	9.01754
0.8650	9.04375
0.8660	9.07004
0.8670	9.09641
0.8680	9.12286
0.8690	9.14939
0.8700	9.17600
0.8710	9.20270
0.8720	9.22947
0.8730	9.25633
0.8740	9.28326
0.8750	9.31028
0.8760	9.33739
0.8770	9.36457
0.8780	9.39184
0.8790	9.41919
0.8800	9.44662
0.8810	9.47414
0.8820	9.50174
0.8830	9.52942
0.8840	9.55719
0.8850	9.58505
0.8860	9.61298
0.8870	9.64101
0.8880	9.66912
0.8890	9.69731
0.8900	9.72559
0.8910	9.75396
0.8920	9.78241
0.8930	9.81095
0.8940	9.83957
0.8950	9.86829
0.8960	9.89709
0.8970	9.92598
0.8980	9.95495
0.8990	9.98402
0.9000	10.01317
0.9010	10.04241
0.9020	10.07174
0.9030	10.10116
0.9040	10.13067
0.9050	10.16027
0.9060	10.18996
0.9070	10.21974
0.9080	10.24961
0.9090	10.27957
0.9100	10.30962
0.9110	10.33977
0.9120	10.37000
0.9130	10.40033
0.9140	10.43075
0.9150	10.46126
0.9160	10.49187
0.9170	10.52257
0.9180	10.55336
0.9190	10.58425
0.9200	10.61523
0.9210	10.64630
0.9220	10.67747
0.9230	10.70874
0.9240	10.74010
0.9250	10.77155
0.9260	10.80310
0.9270	10.83475
0.9280	10.86649
0.9290	10.89833
0.9300	10.93027
0.9310	10.96230
0.9320	10.99444
0.9330	11.02667
0.9340	11.05899
0.9350	11.09142
0.9360	11.12395
0.9370	11.15657
0.9380	11.18929
0.9390	11.22212
0.9400	11.25504
0.9410	11.28807
0.9420	11.32119
0.9430	11.35442
0.9440	11.38774
0.9450	11.42117
0.9460	11.45470
0.9470	11.48834
0.9480	11.52207
0.9490	11.55591
0.9500	11.58985
0.9510	11.62389
0.9520	11.65804
0.9530	11.69230
0.9540	11.72665
0.9550	11.76111
0.9560	11.79568
0.9570	11.83035
0.9580	11.86513
0.9590	11.90002
0.9600	11.93501
0.9610	11.97010
0.9620	12.00531
0.9630	12.04062
0.9640	12.07604
0.9650	12.11157
0.9660	12.14720
0.9670	12.18295
0.9680	12.21880
0.9690	12.25476
0.9700	12.29084
0.9710	12.32702
0.9720	12.36331
0.9730	12.39972
0.9740	12.43623
0.9750	12.47286
0.9760	12.50960
0.9770	12.54645
0.9780	12.58341
0.9790	12.62049
0.9800	12.65768
0.9810	12.69498
0.9820	12.73240
0.9830	12.76993
0.9840	12.80758
0.9850	12.84534
0.9860	12.88321
0.9870	12.92120
0.9880	12.95931
0.9890	12.99754
0.9900	13.03588
0.9910	13.07433
0.9920	13.11291
0.9930	13.15160
0.9940	13.19042
0.9950	13.22935
0.9960	13.26840
0.9970	13.30756
0.9980	13.34685
0.9990	13.38626
1.0000	13.42579
/

\plot 0.21 6.5 0.3 6.5 /
\put {$\hat{a}(t)$ }  [l] at 0.32  6.5

\setsolid
\axis bottom label {} shiftedto y=0 ticks numbered from 0.0 to 1.0 by 0.2 unlabeled short from 0.1 to 0.9 by 0.2 /
\axis left   label {} shiftedto x=0 ticks numbered from 0 to 15 by 5 unlabeled short from 1  to 14  by 1 /

\put {$t$}     [r] at 0.995  0.70
\put {$\displaystyle \frac{\ a(t)\ }{16\pi}$}     [l] at 0.02  14.5

\endpicture     \]
\vspace{-8mm}

\caption{The curves of \(a(t)\),  two bounding functions, and the fitted approximation \(\hat{a}(t)\).}

\end{figure}


%% file: fig_04.tex





\begin{figure}[hbt]

\[  \beginpicture   \setlinear

\setcoordinatesystem units <6.4mm,2mm>

\setplotarea x from 65 to 82, y from -30 to 10
\thinlines


\plot
65.39 	9.81695
65.40 	9.56509
65.41 	9.31416
65.42 	9.06414
65.43 	8.81503
65.44 	8.56685
65.45 	8.31960
65.46 	8.07327
65.47 	7.82788
65.48 	7.58342
65.49 	7.33989
65.50 	7.09731
65.51 	6.85567
65.52 	6.61497
65.53 	6.37522
65.54 	6.13642
65.55 	5.89857
65.56 	5.66168
65.57 	5.42574
65.58 	5.19076
65.59 	4.95675
65.60 	4.72369
65.61 	4.49161
65.62 	4.26048
65.63 	4.03033
65.64 	3.80115
65.65 	3.57294
65.66 	3.34571
65.67 	3.11946
65.68 	2.89419
65.69 	2.66989
65.70 	2.44658
65.71 	2.22425
65.72 	2.00291
65.73 	1.78255
65.74 	1.56318
65.75 	1.34480
65.76 	1.12742
65.77 	0.91102
65.78 	0.69562
65.79 	0.48121
65.80 	0.26781
65.81 	0.05539
65.82 	-0.15602
65.83 	-0.36643
65.84 	-0.57584
65.85 	-0.78425
65.86 	-0.99166
65.87 	-1.19807
65.88 	-1.40347
65.89 	-1.60786
65.90 	-1.81125
65.91 	-2.01363
65.92 	-2.21500
65.93 	-2.41537
65.94 	-2.61473
65.95 	-2.81307
65.96 	-3.01041
65.97 	-3.20674
65.98 	-3.40205
65.99 	-3.59636
66.00 	-3.78965
66.01 	-3.98193
66.02 	-4.17319
66.03 	-4.36345
66.04 	-4.55269
66.05 	-4.74092
66.06 	-4.92813
66.07 	-5.11433
66.08 	-5.29952
66.09 	-5.48369
66.10 	-5.66685
66.11 	-5.84900
66.12 	-6.03013
66.13 	-6.21025
66.14 	-6.38935
66.15 	-6.56745
66.16 	-6.74452
66.17 	-6.92059
66.18 	-7.09564
66.19 	-7.26968
66.20 	-7.44271
66.21 	-7.61473
66.22 	-7.78573
66.23 	-7.95573
66.24 	-8.12471
66.25 	-8.29269
66.26 	-8.45965
66.27 	-8.62560
66.28 	-8.79055
66.29 	-8.95449
66.30 	-9.11742
66.31 	-9.27935
66.32 	-9.44027
66.33 	-9.60018
66.34 	-9.75909
66.35 	-9.91700
66.36 	-10.07391
66.37 	-10.22981
66.38 	-10.38471
66.39 	-10.53861
66.40 	-10.69151
66.41 	-10.84341
66.42 	-10.99432
66.43 	-11.14423
66.44 	-11.29315
66.45 	-11.44107
66.46 	-11.58799
66.47 	-11.73393
66.48 	-11.87888
66.49 	-12.02283
66.50 	-12.16580
66.51 	-12.30778
66.52 	-12.44877
66.53 	-12.58878
66.54 	-12.72781
66.55 	-12.86585
66.56 	-13.00292
66.57 	-13.13900
66.58 	-13.27411
66.59 	-13.40824
66.60 	-13.54139
66.61 	-13.67358
66.62 	-13.80479
66.63 	-13.93502
66.64 	-14.06429
66.65 	-14.19260
66.66 	-14.31993
66.67 	-14.44630
66.68 	-14.57172
66.69 	-14.69616
66.70 	-14.81965
66.71 	-14.94218
66.72 	-15.06376
66.73 	-15.18438
66.74 	-15.30405
66.75 	-15.42277
66.76 	-15.54054
66.77 	-15.65736
66.78 	-15.77324
66.79 	-15.88817
66.80 	-16.00216
66.81 	-16.11522
66.82 	-16.22733
66.83 	-16.33851
66.84 	-16.44876
66.85 	-16.55807
66.86 	-16.66645
66.87 	-16.77391
66.88 	-16.88044
66.89 	-16.98605
66.90 	-17.09073
66.91 	-17.19450
66.92 	-17.29734
66.93 	-17.39927
66.94 	-17.50029
66.95 	-17.60040
66.96 	-17.69960
66.97 	-17.79789
66.98 	-17.89528
66.99 	-17.99176
67.00 	-18.08734
67.01 	-18.18203
67.02 	-18.27582
67.03 	-18.36871
67.04 	-18.46072
67.05 	-18.55183
67.06 	-18.64207
67.07 	-18.73141
67.08 	-18.81987
67.09 	-18.90746
67.10 	-18.99416
67.11 	-19.07999
67.12 	-19.16495
67.13 	-19.24904
67.14 	-19.33226
67.15 	-19.41461
67.16 	-19.49611
67.17 	-19.57674
67.18 	-19.65651
67.19 	-19.73543
67.20 	-19.81350
67.21 	-19.89071
67.22 	-19.96708
67.23 	-20.04260
67.24 	-20.11728
67.25 	-20.19112
67.26 	-20.26412
67.27 	-20.33629
67.28 	-20.40762
67.29 	-20.47813
67.30 	-20.54780
67.31 	-20.61665
67.32 	-20.68468
67.33 	-20.75189
67.34 	-20.81829
67.35 	-20.88386
67.36 	-20.94863
67.37 	-21.01259
67.38 	-21.07574
67.39 	-21.13809
67.40 	-21.19963
67.41 	-21.26039
67.42 	-21.32034
67.43 	-21.37950
67.44 	-21.43787
67.45 	-21.49546
67.46 	-21.55225
67.47 	-21.60827
67.48 	-21.66351
67.49 	-21.71797
67.50 	-21.77166
67.51 	-21.82458
67.52 	-21.87673
67.53 	-21.92812
67.54 	-21.97874
67.55 	-22.02861
67.56 	-22.07771
67.57 	-22.12607
67.58 	-22.17367
67.59 	-22.22052
67.60 	-22.26663
67.61 	-22.31200
67.62 	-22.35663
67.63 	-22.40052
67.64 	-22.44367
67.65 	-22.48610
67.66 	-22.52780
67.67 	-22.56877
67.68 	-22.60902
67.69 	-22.64855
67.70 	-22.68736
67.71 	-22.72546
67.72 	-22.76285
67.73 	-22.79954
67.74 	-22.83552
67.75 	-22.87079
67.76 	-22.90536
67.77 	-22.93924
67.78 	-22.97243
67.79 	-23.00492
67.80 	-23.03673
67.81 	-23.06786
67.82 	-23.09830
67.83 	-23.12807
67.84 	-23.15715
67.85 	-23.18557
67.86 	-23.21332
67.87 	-23.24040
67.88 	-23.26681
67.89 	-23.29256
67.90 	-23.31766
67.91 	-23.34211
67.92 	-23.36589
67.93 	-23.38903
67.94 	-23.41153
67.95 	-23.43338
67.96 	-23.45459
67.97 	-23.47517
67.98 	-23.49511
67.99 	-23.51442
68.00 	-23.53310
68.01 	-23.55116
68.02 	-23.56859
68.03 	-23.58541
68.04 	-23.60161
68.05 	-23.61719
68.06 	-23.63217
68.07 	-23.64653
68.08 	-23.66030
68.09 	-23.67346
68.10 	-23.68602
68.11 	-23.69799
68.12 	-23.70937
68.13 	-23.72015
68.14 	-23.73036
68.15 	-23.73997
68.16 	-23.74901
68.17 	-23.75747
68.18 	-23.76535
68.19 	-23.77266
68.20 	-23.77941
68.21 	-23.78559
68.22 	-23.79120
68.23 	-23.79626
68.24 	-23.80076
68.25 	-23.80471
68.26 	-23.80810
68.27 	-23.81095
68.28 	-23.81325
68.29 	-23.81501
68.30 	-23.81624
68.31 	-23.81692
68.32 	-23.81708
68.33 	-23.81670
68.34 	-23.81580
68.35 	-23.81437
68.36 	-23.81242
68.37 	-23.80996
68.38 	-23.80698
68.39 	-23.80348
68.40 	-23.79948
68.41 	-23.79497
68.42 	-23.78996
68.43 	-23.78445
68.44 	-23.77844
68.45 	-23.77193
68.46 	-23.76493
68.47 	-23.75745
68.48 	-23.74948
68.49 	-23.74102
68.50 	-23.73209
68.51 	-23.72268
68.52 	-23.71279
68.53 	-23.70244
68.54 	-23.69161
68.55 	-23.68032
68.56 	-23.66857
68.57 	-23.65636
68.58 	-23.64369
68.59 	-23.63057
68.60 	-23.61699
68.61 	-23.60297
68.62 	-23.58850
68.63 	-23.57359
68.64 	-23.55824
68.65 	-23.54246
68.66 	-23.52624
68.67 	-23.50959
68.68 	-23.49251
68.69 	-23.47501
68.70 	-23.45708
68.71 	-23.43873
68.72 	-23.41997
68.73 	-23.40079
68.74 	-23.38121
68.75 	-23.36121
68.76 	-23.34081
68.77 	-23.32000
68.78 	-23.29880
68.79 	-23.27720
68.80 	-23.25520
68.81 	-23.23281
68.82 	-23.21003
68.83 	-23.18687
68.84 	-23.16333
68.85 	-23.13940
68.86 	-23.11510
68.87 	-23.09042
68.88 	-23.06537
68.89 	-23.03994
68.90 	-23.01416
68.91 	-22.98801
68.92 	-22.96149
68.93 	-22.93462
68.94 	-22.90739
68.95 	-22.87981
68.96 	-22.85188
68.97 	-22.82360
68.98 	-22.79498
68.99 	-22.76601
69.00 	-22.73671
69.01 	-22.70706
69.02 	-22.67708
69.03 	-22.64677
69.04 	-22.61613
69.05 	-22.58517
69.06 	-22.55388
69.07 	-22.52227
69.08 	-22.49034
69.09 	-22.45809
69.10 	-22.42553
69.11 	-22.39266
69.12 	-22.35949
69.13 	-22.32600
69.14 	-22.29221
69.15 	-22.25813
69.16 	-22.22374
69.17 	-22.18906
69.18 	-22.15409
69.19 	-22.11883
69.20 	-22.08328
69.21 	-22.04744
69.22 	-22.01132
69.23 	-21.97492
69.24 	-21.93825
69.25 	-21.90130
69.26 	-21.86407
69.27 	-21.82658
69.28 	-21.78881
69.29 	-21.75079
69.30 	-21.71250
69.31 	-21.67395
69.32 	-21.63514
69.33 	-21.59608
69.34 	-21.55676
69.35 	-21.51720
69.36 	-21.47738
69.37 	-21.43732
69.38 	-21.39702
69.39 	-21.35647
69.40 	-21.31569
69.41 	-21.27467
69.42 	-21.23342
69.43 	-21.19194
69.44 	-21.15023
69.45 	-21.10829
69.46 	-21.06613
69.47 	-21.02374
69.48 	-20.98114
69.49 	-20.93832
69.50 	-20.89528
69.51 	-20.85203
69.52 	-20.80858
69.53 	-20.76491
69.54 	-20.72104
69.55 	-20.67696
69.56 	-20.63269
69.57 	-20.58821
69.58 	-20.54354
69.59 	-20.49867
69.60 	-20.45361
69.61 	-20.40837
69.62 	-20.36293
69.63 	-20.31731
69.64 	-20.27151
69.65 	-20.22552
69.66 	-20.17936
69.67 	-20.13302
69.68 	-20.08650
69.69 	-20.03981
69.70 	-19.99296
69.71 	-19.94593
69.72 	-19.89874
69.73 	-19.85138
69.74 	-19.80387
69.75 	-19.75619
69.76 	-19.70836
69.77 	-19.66037
69.78 	-19.61223
69.79 	-19.56394
69.80 	-19.51549
69.81 	-19.46691
69.82 	-19.41818
69.83 	-19.36930
69.84 	-19.32028
69.85 	-19.27113
69.86 	-19.22183
69.87 	-19.17241
69.88 	-19.12285
69.89 	-19.07316
69.90 	-19.02334
69.91 	-18.97339
69.92 	-18.92332
69.93 	-18.87313
69.94 	-18.82282
69.95 	-18.77239
69.96 	-18.72184
69.97 	-18.67117
69.98 	-18.62040
69.99 	-18.56951
70.00 	-18.51851
70.01 	-18.46741
70.02 	-18.41620
70.03 	-18.36489
70.04 	-18.31347
70.05 	-18.26196
70.06 	-18.21035
70.07 	-18.15864
70.08 	-18.10683
70.09 	-18.05494
70.10 	-18.00295
70.11 	-17.95088
70.12 	-17.89872
70.13 	-17.84647
70.14 	-17.79415
70.15 	-17.74173
70.16 	-17.68924
70.17 	-17.63668
70.18 	-17.58403
70.19 	-17.53131
70.20 	-17.47852
70.21 	-17.42566
70.22 	-17.37273
70.23 	-17.31973
70.24 	-17.26666
70.25 	-17.21353
70.26 	-17.16034
70.27 	-17.10709
70.28 	-17.05378
70.29 	-17.00041
70.30 	-16.94699
70.31 	-16.89351
70.32 	-16.83998
70.33 	-16.78640
70.34 	-16.73276
70.35 	-16.67909
70.36 	-16.62536
70.37 	-16.57159
70.38 	-16.51778
70.39 	-16.46393
70.40 	-16.41004
70.41 	-16.35611
70.42 	-16.30214
70.43 	-16.24814
70.44 	-16.19410
70.45 	-16.14003
70.46 	-16.08594
70.47 	-16.03181
70.48 	-15.97766
70.49 	-15.92348
70.50 	-15.86928
70.51 	-15.81505
70.52 	-15.76080
70.53 	-15.70653
70.54 	-15.65225
70.55 	-15.59795
70.56 	-15.54363
70.57 	-15.48929
70.58 	-15.43495
70.59 	-15.38059
70.60 	-15.32622
70.61 	-15.27185
70.62 	-15.21747
70.63 	-15.16308
70.64 	-15.10869
70.65 	-15.05429
70.66 	-14.99990
70.67 	-14.94550
70.68 	-14.89110
70.69 	-14.83671
70.70 	-14.78232
70.71 	-14.72794
70.72 	-14.67356
70.73 	-14.61919
70.74 	-14.56483
70.75 	-14.51048
70.76 	-14.45614
70.77 	-14.40181
70.78 	-14.34750
70.79 	-14.29320
70.80 	-14.23892
70.81 	-14.18466
70.82 	-14.13041
70.83 	-14.07619
70.84 	-14.02199
70.85 	-13.96781
70.86 	-13.91365
70.87 	-13.85952
70.88 	-13.80542
70.89 	-13.75134
70.90 	-13.69729
70.91 	-13.64328
70.92 	-13.58929
70.93 	-13.53533
70.94 	-13.48141
70.95 	-13.42752
70.96 	-13.37367
70.97 	-13.31986
70.98 	-13.26608
70.99 	-13.21234
71.00 	-13.15864
71.01 	-13.10498
71.02 	-13.05136
71.03 	-12.99779
71.04 	-12.94426
71.05 	-12.89078
71.06 	-12.83734
71.07 	-12.78395
71.08 	-12.73061
71.09 	-12.67732
71.10 	-12.62408
71.11 	-12.57088
71.12 	-12.51774
71.13 	-12.46466
71.14 	-12.41163
71.15 	-12.35865
71.16 	-12.30573
71.17 	-12.25287
71.18 	-12.20007
71.19 	-12.14732
71.20 	-12.09464
71.21 	-12.04201
71.22 	-11.98945
71.23 	-11.93695
71.24 	-11.88452
71.25 	-11.83215
71.26 	-11.77984
71.27 	-11.72761
71.28 	-11.67543
71.29 	-11.62333
71.30 	-11.57130
71.31 	-11.51933
71.32 	-11.46744
71.33 	-11.41562
71.34 	-11.36387
71.35 	-11.31219
71.36 	-11.26059
71.37 	-11.20906
71.38 	-11.15761
71.39 	-11.10624
71.40 	-11.05494
71.41 	-11.00372
71.42 	-10.95258
71.43 	-10.90152
71.44 	-10.85054
71.45 	-10.79964
71.46 	-10.74883
71.47 	-10.69809
71.48 	-10.64744
71.49 	-10.59688
71.50 	-10.54640
71.51 	-10.49600
71.52 	-10.44569
71.53 	-10.39547
71.54 	-10.34534
71.55 	-10.29529
71.56 	-10.24533
71.57 	-10.19547
71.58 	-10.14569
71.59 	-10.09601
71.60 	-10.04641
71.61 	-9.99691
71.62 	-9.94750
71.63 	-9.89819
71.64 	-9.84897
71.65 	-9.79985
71.66 	-9.75082
71.67 	-9.70189
71.68 	-9.65305
71.69 	-9.60431
71.70 	-9.55568
71.71 	-9.50714
71.72 	-9.45869
71.73 	-9.41035
71.74 	-9.36211
71.75 	-9.31397
71.76 	-9.26593
71.77 	-9.21800
71.78 	-9.17016
71.79 	-9.12243
71.80 	-9.07481
71.81 	-9.02728
71.82 	-8.97987
71.83 	-8.93255
71.84 	-8.88535
71.85 	-8.83825
71.86 	-8.79125
71.87 	-8.74437
71.88 	-8.69759
71.89 	-8.65092
71.90 	-8.60436
71.91 	-8.55791
71.92 	-8.51157
71.93 	-8.46534
71.94 	-8.41922
71.95 	-8.37321
71.96 	-8.32731
71.97 	-8.28152
71.98 	-8.23585
71.99 	-8.19029
72.00 	-8.14484
72.01 	-8.09951
72.02 	-8.05429
72.03 	-8.00918
72.04 	-7.96419
72.05 	-7.91931
72.06 	-7.87456
72.07 	-7.82991
72.08 	-7.78538
72.09 	-7.74097
72.10 	-7.69668
72.11 	-7.65251
72.12 	-7.60845
72.13 	-7.56451
72.14 	-7.52069
72.15 	-7.47698
72.16 	-7.43340
72.17 	-7.38994
72.18 	-7.34659
72.19 	-7.30337
72.20 	-7.26027
72.21 	-7.21729
72.22 	-7.17443
72.23 	-7.13168
72.24 	-7.08907
72.25 	-7.04657
72.26 	-7.00420
72.27 	-6.96195
72.28 	-6.91982
72.29 	-6.87781
72.30 	-6.83593
72.31 	-6.79417
72.32 	-6.75254
72.33 	-6.71103
72.34 	-6.66964
72.35 	-6.62838
72.36 	-6.58725
72.37 	-6.54624
72.38 	-6.50535
72.39 	-6.46459
72.40 	-6.42395
72.41 	-6.38344
72.42 	-6.34306
72.43 	-6.30280
72.44 	-6.26267
72.45 	-6.22267
72.46 	-6.18279
72.47 	-6.14304
72.48 	-6.10342
72.49 	-6.06392
72.50 	-6.02455
72.51 	-5.98532
72.52 	-5.94620
72.53 	-5.90722
72.54 	-5.86836
72.55 	-5.82964
72.56 	-5.79104
72.57 	-5.75257
72.58 	-5.71423
72.59 	-5.67601
72.60 	-5.63793
72.61 	-5.59997
72.62 	-5.56215
72.63 	-5.52445
72.64 	-5.48689
72.65 	-5.44945
72.66 	-5.41214
72.67 	-5.37496
72.68 	-5.33792
72.69 	-5.30100
72.70 	-5.26421
72.71 	-5.22755
72.72 	-5.19103
72.73 	-5.15463
72.74 	-5.11836
72.75 	-5.08223
72.76 	-5.04622
72.77 	-5.01034
72.78 	-4.97460
72.79 	-4.93898
72.80 	-4.90350
72.81 	-4.86815
72.82 	-4.83292
72.83 	-4.79783
72.84 	-4.76287
72.85 	-4.72804
72.86 	-4.69334
72.87 	-4.65877
72.88 	-4.62434
72.89 	-4.59003
72.90 	-4.55585
72.91 	-4.52181
72.92 	-4.48789
72.93 	-4.45411
72.94 	-4.42046
72.95 	-4.38694
72.96 	-4.35355
72.97 	-4.32029
72.98 	-4.28716
72.99 	-4.25416
73.00 	-4.22129
73.01 	-4.18856
73.02 	-4.15595
73.03 	-4.12347
73.04 	-4.09113
73.05 	-4.05892
73.06 	-4.02683
73.07 	-3.99488
73.08 	-3.96306
73.09 	-3.93137
73.10 	-3.89981
73.11 	-3.86837
73.12 	-3.83707
73.13 	-3.80590
73.14 	-3.77486
73.15 	-3.74395
73.16 	-3.71317
73.17 	-3.68252
73.18 	-3.65199
73.19 	-3.62160
73.20 	-3.59134
73.21 	-3.56121
73.22 	-3.53120
73.23 	-3.50133
73.24 	-3.47158
73.25 	-3.44197
73.26 	-3.41248
73.27 	-3.38312
73.28 	-3.35390
73.29 	-3.32480
73.30 	-3.29582
73.31 	-3.26698
73.32 	-3.23826
73.33 	-3.20968
73.34 	-3.18122
73.35 	-3.15289
73.36 	-3.12468
73.37 	-3.09661
73.38 	-3.06866
73.39 	-3.04084
73.40 	-3.01315
73.41 	-2.98558
73.42 	-2.95814
73.43 	-2.93083
73.44 	-2.90364
73.45 	-2.87658
73.46 	-2.84965
73.47 	-2.82284
73.48 	-2.79616
73.49 	-2.76960
73.50 	-2.74318
73.51 	-2.71687
73.52 	-2.69069
73.53 	-2.66464
73.54 	-2.63871
73.55 	-2.61291
73.56 	-2.58723
73.57 	-2.56167
73.58 	-2.53624
73.59 	-2.51094
73.60 	-2.48576
73.61 	-2.46070
73.62 	-2.43576
73.63 	-2.41095
73.64 	-2.38626
73.65 	-2.36170
73.66 	-2.33725
73.67 	-2.31293
73.68 	-2.28873
73.69 	-2.26466
73.70 	-2.24070
73.71 	-2.21687
73.72 	-2.19316
73.73 	-2.16957
73.74 	-2.14610
73.75 	-2.12275
73.76 	-2.09953
73.77 	-2.07642
73.78 	-2.05343
73.79 	-2.03057
73.80 	-2.00782
73.81 	-1.98519
73.82 	-1.96268
73.83 	-1.94029
73.84 	-1.91802
73.85 	-1.89587
73.86 	-1.87384
73.87 	-1.85192
73.88 	-1.83013
73.89 	-1.80845
73.90 	-1.78688
73.91 	-1.76544
73.92 	-1.74411
73.93 	-1.72290
73.94 	-1.70180
73.95 	-1.68082
73.96 	-1.65996
73.97 	-1.63921
73.98 	-1.61858
73.99 	-1.59806
74.00 	-1.57766
74.01 	-1.55737
74.02 	-1.53720
74.03 	-1.51714
74.04 	-1.49719
74.05 	-1.47736
74.06 	-1.45764
74.07 	-1.43803
74.08 	-1.41854
74.09 	-1.39916
74.10 	-1.37989
74.11 	-1.36073
74.12 	-1.34169
74.13 	-1.32275
74.14 	-1.30393
74.15 	-1.28522
74.16 	-1.26662
74.17 	-1.24813
74.18 	-1.22975
74.19 	-1.21148
74.20 	-1.19331
74.21 	-1.17526
74.22 	-1.15732
74.23 	-1.13948
74.24 	-1.12175
74.25 	-1.10414
74.26 	-1.08663
74.27 	-1.06922
74.28 	-1.05192
74.29 	-1.03473
74.30 	-1.01765
74.31 	-1.00067
74.32 	-0.98380
74.33 	-0.96704
74.34 	-0.95037
74.35 	-0.93382
74.36 	-0.91737
74.37 	-0.90102
74.38 	-0.88478
74.39 	-0.86864
74.40 	-0.85261
74.41 	-0.83668
74.42 	-0.82085
74.43 	-0.80512
74.44 	-0.78950
74.45 	-0.77398
74.46 	-0.75856
74.47 	-0.74324
74.48 	-0.72802
74.49 	-0.71290
74.50 	-0.69789
74.51 	-0.68298
74.52 	-0.66816
74.53 	-0.65344
74.54 	-0.63882
74.55 	-0.62431
74.56 	-0.60989
74.57 	-0.59557
74.58 	-0.58134
74.59 	-0.56722
74.60 	-0.55319
74.61 	-0.53926
74.62 	-0.52542
74.63 	-0.51169
74.64 	-0.49804
74.65 	-0.48450
74.66 	-0.47105
74.67 	-0.45769
74.68 	-0.44443
74.69 	-0.43126
74.70 	-0.41819
74.71 	-0.40521
74.72 	-0.39233
74.73 	-0.37954
74.74 	-0.36684
74.75 	-0.35423
74.76 	-0.34172
74.77 	-0.32930
74.78 	-0.31696
74.79 	-0.30473
74.80 	-0.29258
74.81 	-0.28052
74.82 	-0.26855
74.83 	-0.25668
74.84 	-0.24489
74.85 	-0.23319
74.86 	-0.22158
74.87 	-0.21006
74.88 	-0.19863
74.89 	-0.18728
74.90 	-0.17603
74.91 	-0.16486
74.92 	-0.15378
74.93 	-0.14278
74.94 	-0.13187
74.95 	-0.12105
74.96 	-0.11031
74.97 	-0.09966
74.98 	-0.08910
74.99 	-0.07861
75.00 	-0.06821
75.01 	-0.05790
75.02 	-0.04767
75.03 	-0.03753
75.04 	-0.02746
75.05 	-0.01748
75.06 	-0.00758
75.07 	0.00223
75.08 	0.01196
75.09 	0.02162
75.10 	0.03119
75.11 	0.04068
75.12 	0.05009
75.13 	0.05942
75.14 	0.06867
75.15 	0.07784
75.16 	0.08693
75.17 	0.09594
75.18 	0.10488
75.19 	0.11373
75.20 	0.12251
75.21 	0.13121
75.22 	0.13984
75.23 	0.14838
75.24 	0.15685
75.25 	0.16525
75.26 	0.17356
75.27 	0.18181
75.28 	0.18997
75.29 	0.19807
75.30 	0.20608
75.31 	0.21403
75.32 	0.22190
75.33 	0.22969
75.34 	0.23741
75.35 	0.24506
75.36 	0.25264
75.37 	0.26015
75.38 	0.26758
75.39 	0.27494
75.40 	0.28223
75.41 	0.28945
75.42 	0.29660
75.43 	0.30367
75.44 	0.31068
75.45 	0.31762
75.46 	0.32449
75.47 	0.33129
75.48 	0.33802
75.49 	0.34468
75.50 	0.35128
75.51 	0.35780
75.52 	0.36426
75.53 	0.37065
75.54 	0.37698
75.55 	0.38324
75.56 	0.38943
75.57 	0.39556
75.58 	0.40162
75.59 	0.40762
75.60 	0.41355
75.61 	0.41941
75.62 	0.42522
75.63 	0.43096
75.64 	0.43663
75.65 	0.44224
75.66 	0.44779
75.67 	0.45328
75.68 	0.45870
75.69 	0.46406
75.70 	0.46936
75.71 	0.47460
75.72 	0.47978
75.73 	0.48489
75.74 	0.48995
75.75 	0.49495
75.76 	0.49988
75.77 	0.50476
75.78 	0.50958
75.79 	0.51434
75.80 	0.51903
75.81 	0.52368
75.82 	0.52826
75.83 	0.53279
75.84 	0.53725
75.85 	0.54167
75.86 	0.54602
75.87 	0.55032
75.88 	0.55456
75.89 	0.55875
75.90 	0.56288
75.91 	0.56695
75.92 	0.57097
75.93 	0.57494
75.94 	0.57885
75.95 	0.58271
75.96 	0.58651
75.97 	0.59026
75.98 	0.59396
75.99 	0.59761
76.00 	0.60120
76.01 	0.60474
76.02 	0.60823
76.03 	0.61166
76.04 	0.61505
76.05 	0.61838
76.06 	0.62166
76.07 	0.62490
76.08 	0.62808
76.09 	0.63121
76.10 	0.63430
76.11 	0.63733
76.12 	0.64031
76.13 	0.64325
76.14 	0.64614
76.15 	0.64898
76.16 	0.65177
76.17 	0.65451
76.18 	0.65721
76.19 	0.65986
76.20 	0.66246
76.21 	0.66502
76.22 	0.66753
76.23 	0.67000
76.24 	0.67242
76.25 	0.67479
76.26 	0.67712
76.27 	0.67941
76.28 	0.68165
76.29 	0.68384
76.30 	0.68600
76.31 	0.68810
76.32 	0.69017
76.33 	0.69219
76.34 	0.69417
76.35 	0.69611
76.36 	0.69800
76.37 	0.69985
76.38 	0.70166
76.39 	0.70343
76.40 	0.70516
76.41 	0.70685
76.42 	0.70850
76.43 	0.71010
76.44 	0.71167
76.45 	0.71320
76.46 	0.71468
76.47 	0.71613
76.48 	0.71754
76.49 	0.71891
76.50 	0.72024
76.51 	0.72154
76.52 	0.72279
76.53 	0.72401
76.54 	0.72519
76.55 	0.72634
76.56 	0.72744
76.57 	0.72852
76.58 	0.72955
76.59 	0.73055
76.60 	0.73151
76.61 	0.73243
76.62 	0.73332
76.63 	0.73418
76.64 	0.73500
76.65 	0.73579
76.66 	0.73654
76.67 	0.73726
76.68 	0.73794
76.69 	0.73859
76.70 	0.73921
76.71 	0.73979
76.72 	0.74034
76.73 	0.74086
76.74 	0.74135
76.75 	0.74180
76.76 	0.74222
76.77 	0.74261
76.78 	0.74298
76.79 	0.74330
76.80 	0.74360
76.81 	0.74387
76.82 	0.74410
76.83 	0.74431
76.84 	0.74448
76.85 	0.74463
76.86 	0.74475
76.87 	0.74483
76.88 	0.74489
76.89 	0.74492
76.90 	0.74492
76.91 	0.74489
76.92 	0.74484
76.93 	0.74475
76.94 	0.74464
76.95 	0.74450
76.96 	0.74434
76.97 	0.74415
76.98 	0.74393
76.99 	0.74368
77.00 	0.74341
77.01 	0.74311
77.02 	0.74279
77.03 	0.74244
77.04 	0.74206
77.05 	0.74166
77.06 	0.74123
77.07 	0.74078
77.08 	0.74031
77.09 	0.73981
77.10 	0.73929
77.11 	0.73874
77.12 	0.73817
77.13 	0.73757
77.14 	0.73696
77.15 	0.73632
77.16 	0.73565
77.17 	0.73497
77.18 	0.73426
77.19 	0.73353
77.20 	0.73278
77.21 	0.73200
77.22 	0.73120
77.23 	0.73039
77.24 	0.72955
77.25 	0.72869
77.26 	0.72780
77.27 	0.72691
77.28 	0.72598
77.29 	0.72504
77.30 	0.72408
77.31 	0.72310
77.32 	0.72210
77.33 	0.72108
77.34 	0.72004
77.35 	0.71898
77.36 	0.71790
77.37 	0.71680
77.38 	0.71569
77.39 	0.71455
77.40 	0.71340
77.41 	0.71223
77.42 	0.71105
77.43 	0.70984
77.44 	0.70862
77.45 	0.70738
77.46 	0.70612
77.47 	0.70485
77.48 	0.70356
77.49 	0.70226
77.50 	0.70093
77.51 	0.69959
77.52 	0.69824
77.53 	0.69687
77.54 	0.69548
77.55 	0.69408
77.56 	0.69267
77.57 	0.69124
77.58 	0.68979
77.59 	0.68833
77.60 	0.68686
77.61 	0.68537
77.62 	0.68387
77.63 	0.68235
77.64 	0.68082
77.65 	0.67927
77.66 	0.67771
77.67 	0.67614
77.68 	0.67456
77.69 	0.67296
77.70 	0.67135
77.71 	0.66972
77.72 	0.66809
77.73 	0.66644
77.74 	0.66478
77.75 	0.66311
77.76 	0.66142
77.77 	0.65973
77.78 	0.65802
77.79 	0.65630
77.80 	0.65457
77.81 	0.65283
77.82 	0.65108
77.83 	0.64932
77.84 	0.64754
77.85 	0.64576
77.86 	0.64396
77.87 	0.64216
77.88 	0.64034
77.89 	0.63852
77.90 	0.63669
77.91 	0.63484
77.92 	0.63299
77.93 	0.63113
77.94 	0.62926
77.95 	0.62738
77.96 	0.62548
77.97 	0.62359
77.98 	0.62168
77.99 	0.61976
78.00 	0.61784
78.01 	0.61591
78.02 	0.61397
78.03 	0.61202
78.04 	0.61007
78.05 	0.60810
78.06 	0.60613
78.07 	0.60416
78.08 	0.60217
78.09 	0.60018
78.10 	0.59818
78.11 	0.59617
78.12 	0.59416
78.13 	0.59214
78.14 	0.59012
78.15 	0.58808
78.16 	0.58605
78.17 	0.58400
78.18 	0.58195
78.19 	0.57990
78.20 	0.57783
78.21 	0.57577
78.22 	0.57369
78.23 	0.57162
78.24 	0.56953
78.25 	0.56745
78.26 	0.56535
78.27 	0.56325
78.28 	0.56115
78.29 	0.55904
78.30 	0.55693
78.31 	0.55482
78.32 	0.55270
78.33 	0.55057
78.34 	0.54844
78.35 	0.54631
78.36 	0.54417
78.37 	0.54203
78.38 	0.53989
78.39 	0.53774
78.40 	0.53559
78.41 	0.53344
78.42 	0.53128
78.43 	0.52912
78.44 	0.52696
78.45 	0.52479
78.46 	0.52262
78.47 	0.52045
78.48 	0.51828
78.49 	0.51610
78.50 	0.51392
78.51 	0.51174
78.52 	0.50956
78.53 	0.50737
78.54 	0.50519
78.55 	0.50300
78.56 	0.50081
78.57 	0.49861
78.58 	0.49642
78.59 	0.49422
78.60 	0.49203
78.61 	0.48983
78.62 	0.48763
78.63 	0.48543
78.64 	0.48323
78.65 	0.48103
78.66 	0.47882
78.67 	0.47662
78.68 	0.47442
78.69 	0.47221
78.70 	0.47001
78.71 	0.46780
78.72 	0.46560
78.73 	0.46339
78.74 	0.46118
78.75 	0.45898
78.76 	0.45677
78.77 	0.45456
78.78 	0.45236
78.79 	0.45015
78.80 	0.44795
78.81 	0.44575
78.82 	0.44354
78.83 	0.44134
78.84 	0.43913
78.85 	0.43693
78.86 	0.43473
78.87 	0.43253
78.88 	0.43033
78.89 	0.42813
78.90 	0.42593
78.91 	0.42374
78.92 	0.42155
78.93 	0.41935
78.94 	0.41716
78.95 	0.41497
78.96 	0.41278
78.97 	0.41060
78.98 	0.40841
78.99 	0.40623
79.00 	0.40405
79.01 	0.40187
79.02 	0.39969
79.03 	0.39751
79.04 	0.39534
79.05 	0.39317
79.06 	0.39100
79.07 	0.38883
79.08 	0.38667
79.09 	0.38451
79.10 	0.38235
79.11 	0.38019
79.12 	0.37804
79.13 	0.37589
79.14 	0.37374
79.15 	0.37159
79.16 	0.36945
79.17 	0.36731
79.18 	0.36517
79.19 	0.36304
79.20 	0.36091
79.21 	0.35878
79.22 	0.35666
79.23 	0.35454
79.24 	0.35242
79.25 	0.35030
79.26 	0.34820
79.27 	0.34609
79.28 	0.34398
79.29 	0.34188
79.30 	0.33979
79.31 	0.33769
79.32 	0.33561
79.33 	0.33352
79.34 	0.33144
79.35 	0.32936
79.36 	0.32729
79.37 	0.32522
79.38 	0.32316
79.39 	0.32109
79.40 	0.31904
79.41 	0.31698
79.42 	0.31493
79.43 	0.31289
79.44 	0.31085
79.45 	0.30881
79.46 	0.30678
79.47 	0.30476
79.48 	0.30273
79.49 	0.30071
79.50 	0.29870
79.51 	0.29669
79.52 	0.29469
79.53 	0.29269
79.54 	0.29069
79.55 	0.28870
79.56 	0.28672
79.57 	0.28474
79.58 	0.28276
79.59 	0.28079
79.60 	0.27883
79.61 	0.27687
79.62 	0.27491
79.63 	0.27296
79.64 	0.27102
79.65 	0.26908
79.66 	0.26714
79.67 	0.26521
79.68 	0.26329
79.69 	0.26137
79.70 	0.25945
79.71 	0.25755
79.72 	0.25564
79.73 	0.25375
79.74 	0.25185
79.75 	0.24997
79.76 	0.24809
79.77 	0.24621
79.78 	0.24434
79.79 	0.24247
79.80 	0.24061
79.81 	0.23876
79.82 	0.23691
79.83 	0.23507
79.84 	0.23323
79.85 	0.23140
79.86 	0.22958
79.87 	0.22776
79.88 	0.22595
79.89 	0.22414
79.90 	0.22234
79.91 	0.22054
79.92 	0.21875
79.93 	0.21697
79.94 	0.21519
79.95 	0.21342
79.96 	0.21165
79.97 	0.20989
79.98 	0.20814
79.99 	0.20639
80.00 	0.20465
80.01 	0.20291
80.02 	0.20118
80.03 	0.19945
80.04 	0.19774
80.05 	0.19603
80.06 	0.19432
80.07 	0.19262
80.08 	0.19093
80.09 	0.18924
80.10 	0.18756
80.11 	0.18589
80.12 	0.18422
80.13 	0.18256
80.14 	0.18090
80.15 	0.17925
80.16 	0.17761
80.17 	0.17598
80.18 	0.17435
80.19 	0.17272
80.20 	0.17111
80.21 	0.16949
80.22 	0.16789
80.23 	0.16629
80.24 	0.16470
80.25 	0.16312
80.26 	0.16154
80.27 	0.15996
80.28 	0.15840
80.29 	0.15684
80.30 	0.15529
80.31 	0.15374
80.32 	0.15220
80.33 	0.15067
80.34 	0.14914
80.35 	0.14762
80.36 	0.14611
80.37 	0.14460
80.38 	0.14310
80.39 	0.14161
80.40 	0.14012
80.41 	0.13864
80.42 	0.13716
80.43 	0.13570
80.44 	0.13424
80.45 	0.13278
80.46 	0.13133
80.47 	0.12989
80.48 	0.12846
80.49 	0.12703
80.50 	0.12561
80.51 	0.12420
80.52 	0.12279
80.53 	0.12139
80.54 	0.11999
80.55 	0.11860
80.56 	0.11722
80.57 	0.11585
80.58 	0.11448
80.59 	0.11312
80.60 	0.11176
80.61 	0.11041
80.62 	0.10907
80.63 	0.10773
80.64 	0.10641
80.65 	0.10508
80.66 	0.10377
80.67 	0.10246
80.68 	0.10116
80.69 	0.09986
80.70 	0.09858
80.71 	0.09729
80.72 	0.09602
80.73 	0.09475
80.74 	0.09348
80.75 	0.09223
80.76 	0.09098
80.77 	0.08974
80.78 	0.08850
80.79 	0.08727
80.80 	0.08605
80.81 	0.08483
80.82 	0.08363
80.83 	0.08242
80.84 	0.08123
80.85 	0.08004
80.86 	0.07886
80.87 	0.07768
80.88 	0.07651
80.89 	0.07535
80.90 	0.07419
80.91 	0.07304
80.92 	0.07190
80.93 	0.07076
80.94 	0.06963
80.95 	0.06850
80.96 	0.06739
80.97 	0.06628
80.98 	0.06517
80.99 	0.06407
81.00 	0.06298
81.01 	0.06190
81.02 	0.06082
81.03 	0.05975
81.04 	0.05868
81.05 	0.05762
81.06 	0.05657
81.07 	0.05552
81.08 	0.05448
81.09 	0.05345
81.10 	0.05242
81.11 	0.05140
81.12 	0.05038
81.13 	0.04938
81.14 	0.04838
81.15 	0.04738
81.16 	0.04639
81.17 	0.04541
81.18 	0.04443
81.19 	0.04346
81.20 	0.04250
81.21 	0.04154
81.22 	0.04059
81.23 	0.03964
81.24 	0.03871
81.25 	0.03777
81.26 	0.03685
81.27 	0.03593
81.28 	0.03501
81.29 	0.03411
81.30 	0.03320
81.31 	0.03231
81.32 	0.03142
81.33 	0.03054
81.34 	0.02966
81.35 	0.02879
81.36 	0.02792
81.37 	0.02707
81.38 	0.02621
81.39 	0.02537
81.40 	0.02453
81.41 	0.02369
81.42 	0.02286
81.43 	0.02204
81.44 	0.02122
81.45 	0.02041
81.46 	0.01961
81.47 	0.01881
81.48 	0.01802
81.49 	0.01723
81.50 	0.01645
81.51 	0.01568
81.52 	0.01490
81.53 	0.01414
81.54 	0.01339
81.55 	0.01263
81.56 	0.01189
81.57 	0.01115
81.58 	0.01041
81.59 	0.00969
81.60 	0.00896
81.61 	0.00824
81.62 	0.00753
81.63 	0.00683
81.64 	0.00613
81.65 	0.00543
81.66 	0.00475
81.67 	0.00406
81.68 	0.00338
81.69 	0.00271
81.70 	0.00205
81.71 	0.00139
81.72 	0.00073
81.73 	0.00008
81.74 	-0.00056
81.75 	-0.00120
81.76 	-0.00183
81.77 	-0.00246
81.78 	-0.00308
81.79 	-0.00370
81.80 	-0.00431
81.81 	-0.00492
81.82 	-0.00552
81.83 	-0.00612
81.84 	-0.00670
81.85 	-0.00729
81.86 	-0.00787
81.87 	-0.00844
81.88 	-0.00901
81.89 	-0.00958
81.90 	-0.01013
81.91 	-0.01069
81.92 	-0.01124
81.93 	-0.01178
81.94 	-0.01232
81.95 	-0.01285
81.96 	-0.01338
81.97 	-0.01390
81.98 	-0.01442
81.99 	-0.01493
/

\setsolid

\axis bottom label {} shiftedto y=0.0  ticks numbered from 65 to 82 by 5     unlabeled short from   65 to 82 by 1 /
\axis left label   {} shiftedto x=65   ticks numbered from -30.0 to 10.0 by 10   unlabeled short from   -25 to 5 by 10 /

\put {$U(y)\times 10^{10}$ (Riemann)}   [l] at  65.7 9.8
\put {$y$}   [r] at  82 2.0
\put {$U''(y) < 0$} [l] at 73.0 5.0
\arrow <1,5mm>   [0.25,0.75] from  74.5 3.0 to 75.05 0.5

\endpicture     \]
\vspace{-3mm}

\caption{A section of the curve \(U(y)\) for the Riemann \(\xi(s)\) function.}

\end{figure}


%% file: fig_06.tex




\baselineskip1.8em

\begin{figure}[bht]

\[  \beginpicture   \setlinear

\setcoordinatesystem units <1cm,1cm>
\setplotarea x from 0 to 10, y from 0 to 6.1
\thinlines

\arrow <1,5mm>   [0.25,0.75] from  1.0 0.0 to 11.0 0.0
\arrow <1,5mm>   [0.25,0.75] from  5.0 0.0 to  5.0 6.3

\setdashes
\setdashpattern <1mm, 0.7mm, 1mm, 0.7mm>
\plot 5.3 0.0 5.3 5.5 /
\plot 4.7 0.0 4.7 5.5 /
\put {$x$}     [r] at 11.0  0.25
\put {$y$}     [r] at 5.4  6.3

\put { $\bullet$} [r] at 5.1 1.0
\put { $\bullet$} [r] at 5.1 3.0
\put { $\bullet$} [r] at 5.1 5.0

\put {0}     [r] at 5.1  -0.20

\setsolid
\arrow <1,5mm>   [0.25,0.75] from  6.20 0.45 to  5.3 0.45
\arrow <1,5mm>   [0.25,0.75] from  4.09 0.45 to  5.0 0.45
\put {$\varepsilon(y)$}     [c] at 5.9 0.7

\put {$ v_x(0,y) \not= 0$}   [r] at 3.0  1.4
\arrow <1,5mm>   [0.25,0.75] from  3.2 1.4 to  5.0 1.7

\put {$ v(x,y) \!=\! v_x(0,y)x\! +\!  {\mathcal O}(x^2)\!  \not=\! 0$}   [l] at 6.0  1.4
\arrow <1,5mm>   [0.25,0.75] from  5.9 1.5 to  5.18 1.7

\put {$ v_x(0,y_m) = 0$}   [r] at 3.0  4.3
\arrow <1,5mm>   [0.25,0.75] from  3.2 4.3 to  4.9 3.06

\put {$ u_y(0,y_m) \!=\!  0$}   [l] at 6.0  4.6
\put {$ u(x, y) \approx (u_{\max} , u_{\min}) \not=0$ }   [l] at 6.05  3.8
\arrow <1,5mm>   [0.25,0.75] from  6.0 3.8 to  5.18 2.99

\setdashpattern <1mm, 0.7mm, 1mm, 0.7mm>
\plot 
5.35	3
5.349309355	3.021976682
5.347240145	3.043866632
5.343800538	3.06558346
5.339004106	3.087041461
5.332869781	3.108155948
5.32542177	3.128843593
5.316689468	3.149022752
5.306707338	3.168613786
5.295514774	3.187539378
5.283155948	3.205724838
5.269679635	3.223098396
5.25513902	3.239591487
5.239591487	3.25513902
5.223098396	3.269679635
5.205724838	3.283155948
5.187539378	3.295514774
5.168613786	3.306707338
5.149022752	3.316689468
5.128843593	3.32542177
5.108155948	3.332869781
5.087041461	3.339004106
5.06558346	3.343800538
5.043866632	3.347240145
5.021976682	3.349309355
5	3.35
4.978023318	3.349309355
4.956133368	3.347240145
4.93441654	3.343800538
4.912958539	3.339004106
4.891844052	3.332869781
4.871156407	3.32542177
4.850977248	3.316689468
4.831386214	3.306707338
4.812460622	3.295514774
4.794275162	3.283155948
4.776901604	3.269679635
4.760408513	3.25513902
4.74486098	3.239591487
4.730320365	3.223098396
4.716844052	3.205724838
4.704485226	3.187539378
4.693292662	3.168613786
4.683310532	3.149022752
4.67457823	3.128843593
4.667130219	3.108155948
4.660995894	3.087041461
4.656199462	3.06558346
4.652759855	3.043866632
4.650690645	3.021976682
4.65	3
4.650690645	2.978023318
4.652759855	2.956133368
4.656199462	2.93441654
4.660995894	2.912958539
4.667130219	2.891844052
4.67457823	2.871156407
4.683310532	2.850977248
4.693292662	2.831386214
4.704485226	2.812460622
4.716844052	2.794275162
4.730320365	2.776901604
4.74486098	2.760408513
4.760408513	2.74486098
4.776901604	2.730320365
4.794275162	2.716844052
4.812460622	2.704485226
4.831386214	2.693292662
4.850977248	2.683310532
4.871156407	2.67457823
4.891844052	2.667130219
4.912958539	2.660995894
4.93441654	2.656199462
4.956133368	2.652759855
4.978023318	2.650690645
5	2.65
5.021976682	2.650690645
5.043866632	2.652759855
5.06558346	2.656199462
5.087041461	2.660995894
5.108155948	2.667130219
5.128843593	2.67457823
5.149022752	2.683310532
5.168613786	2.693292662
5.187539378	2.704485226
5.205724838	2.716844052
5.223098396	2.730320365
5.239591487	2.74486098
5.25513902	2.760408513
5.269679635	2.776901604
5.283155948	2.794275162
5.295514774	2.812460622
5.306707338	2.831386214
5.316689468	2.850977248
5.32542177	2.871156407
5.332869781	2.891844052
5.339004106	2.912958539
5.343800538	2.93441654
5.347240145	2.956133368
5.349309355	2.978023318
5.35	3
/

\endpicture     \]
\vspace{-6mm}
\caption{Two \(\varepsilon\)-domains \(0<|x|<\varepsilon(y)\) attached to both sides of the critical line, in each of which either \(u\neq 0\) or \(v\neq 0\).}

\end{figure}


%% file: fig_07.tex





\begin{figure}[bht]

\[  \beginpicture   \setlinear

\setcoordinatesystem units <10mm,10mm>
\setplotarea x from -1 to 10, y from 0 to 6.1
\thinlines

\arrow <1,5mm>   [0.25,0.75] from  0.0 0.0 to 10.0 0.0
\arrow <1,5mm>   [0.25,0.75] from  0.0 0.0 to  0.0 6.3

\plot
0.0 	1.50000
0.2 	1.49956
0.4 	1.49824
0.6 	1.49604
0.8 	1.49296
1.0 	1.48900
1.2 	1.48416
1.4 	1.47844
1.6 	1.47184
1.8 	1.46436
2.0 	1.45600
2.2 	1.44676
2.4 	1.43664
2.6 	1.42564
2.8 	1.41376
3.0 	1.40100
3.2 	1.38736
3.4 	1.37284
3.6 	1.35744
3.8 	1.34116
4.0 	1.32400
4.2 	1.30596
4.4 	1.28704
4.6 	1.26724
4.8 	1.24656
5.0 	1.22500
5.2 	1.20256
5.4 	1.17924
5.6 	1.15504
5.8 	1.12996
6.0 	1.10400
6.2 	1.07716
6.4 	1.04944
6.6 	1.02084
6.8 	0.99136
7.0 	0.96100
7.2 	0.92976
7.4 	0.89764
7.6 	0.86464
7.8 	0.83076
8.0 	0.79600
8.2 	0.76036
8.4 	0.72384
8.6 	0.68644
8.8 	0.64816
9.0 	0.60900
9.2 	0.56896
9.4 	0.52804
9.6 	0.48624
9.8 	0.44356
 /

\plot
0.0 	3.50000
0.2 	3.49954
0.4 	3.49816
0.6 	3.49586
0.8 	3.49264
1.0 	3.48850
1.2 	3.48344
1.4 	3.47746
1.6 	3.47056
1.8 	3.46274
2.0 	3.45400
2.2 	3.44434
2.4 	3.43376
2.6 	3.42226
2.8 	3.40984
3.0 	3.39650
3.2 	3.38224
3.4 	3.36706
3.6 	3.35096
3.8 	3.33394
4.0 	3.31600
4.2 	3.29714
4.4 	3.27736
4.6 	3.25666
4.8 	3.23504
5.0 	3.21250
5.2 	3.18904
5.4 	3.16466
5.6 	3.13936
5.8 	3.11314
6.0 	3.08600
6.2 	3.05794
6.4 	3.02896
6.6 	2.99906
6.8 	2.96824
7.0 	2.93650
7.2 	2.90384
7.4 	2.87026
7.6 	2.83576
7.8 	2.80034
8.0 	2.76400
8.2 	2.72674
8.4 	2.68856
8.6 	2.64946
8.8 	2.60944
9.0 	2.56850
9.2 	2.52664
9.4 	2.48386
9.6 	2.44016
9.8 	2.39554
/

\plot
0.0 	5.50000
0.2 	5.49956
0.4 	5.49824
0.6 	5.49604
0.8 	5.49296
1.0 	5.48900
1.2 	5.48416
1.4 	5.47844
1.6 	5.47184
1.8 	5.46436
2.0 	5.45600
2.2 	5.44676
2.4 	5.43664
2.6 	5.42564
2.8 	5.41376
3.0 	5.40100
3.2 	5.38736
3.4 	5.37284
3.6 	5.35744
3.8 	5.34116
4.0 	5.32400
4.2 	5.30596
4.4 	5.28704
4.6 	5.26724
4.8 	5.24656
5.0 	5.22500
5.2 	5.20256
5.4 	5.17924
5.6 	5.15504
5.8 	5.12996
6.0 	5.10400
6.2 	5.07716
6.4 	5.04944
6.6 	5.02084
6.8 	4.99136
7.0 	4.96100
7.2 	4.92976
7.4 	4.89764
7.6 	4.86464
7.8 	4.83076
8.0 	4.79600
8.2 	4.76036
8.4 	4.72384
8.6 	4.68644
8.8 	4.64816
9.0 	4.60900
9.2 	4.56896
9.4 	4.52804
9.6 	4.48624
9.8 	4.44356
/

\setsolid
\arrow <1,5mm>   [0.25,0.75] from  5.99 3.09 to  7.35 2.95
\arrow <1,5mm>   [0.25,0.75] from  5.99 3.11 to  6.18 4.4

\put {$0$}     [r] at -0.1  -0.20

\put {$-u_y = v_x < 0$}     [r] at -0.1  4.4
\put {$-u_y = v_x > 0$}     [r] at -0.1  2.3

\put {$u_{\max}$} [r] at -0.1 5.5
\put {$u_{\min}$} [r] at -0.1 3.5

\put {$y =\varphi_{2k+1}(x)$} [r] at 4.6 5.7
\put {$y =\varphi_{2k}(x)$} [r] at 5.2 3.65
\put {$y =\varphi_{2k-1}(x)$} [r] at 6.5 1.6

\put {$\Omega_{2k+1}: v < 0$} [l] at 0.5 4.4
\put {$\Omega_{2k}: v > 0$} [l] at 0.5 2.3

\put {${\bf e}_1$} [l] at 7.2 3.25
\put {${\bf e}_2$} [l] at 6.3 4.45

\put {$x$} [l] at 9.5 -0.25
\put {$y$} [r] at 0.43 6.3

\endpicture     \]
\vspace{-6mm}

\caption{Case (i): the boundary curves of the subregions do not intersect or branch.}

\end{figure}


%% file: fig_8.tex





\begin{figure}

\[  \beginpicture   \setlinear

\setcoordinatesystem units <180mm,0.5mm>

\setplotarea x from 2.9 to 3.5, y from -80 to 80
\thinlines


\plot
2.970 	56.93271
2.980 	35.70751
2.981 	33.71522
2.982 	31.74601
2.983 	29.79974
2.984 	27.87626
2.985 	25.97542
2.986 	24.09708
2.987 	22.24110
2.988 	20.40734
2.989 	18.59564
2.990 	16.80587
2.991 	15.03789
2.992 	13.29155
2.993 	11.56671
2.994 	9.86324
2.995 	8.18100
2.996 	6.51985
2.997 	4.87964
2.998 	3.26025
2.999 	1.66153
3.000 	0.08336
3.001 	-1.47442
3.002 	-3.01192
3.003 	-4.52928
3.004 	-6.02665
3.005 	-7.50414
3.006 	-8.96190
3.007 	-10.40005
3.008 	-11.81873
3.009 	-13.21806
3.010 	-14.59818
3.011 	-15.95922
3.012 	-17.30129
3.013 	-18.62455
3.014 	-19.92910
3.015 	-21.21508
3.016 	-22.48261
3.017 	-23.73182
3.018 	-24.96283
3.019 	-26.17577
3.020 	-27.37077
3.021 	-28.54794
3.022 	-29.70741
3.023 	-30.84931
3.024 	-31.97374
3.025 	-33.08085
3.026 	-34.17073
3.027 	-35.24352
3.028 	-36.29934
3.029 	-37.33830
3.030 	-38.36052
3.031 	-39.36613
3.032 	-40.35523
3.033 	-41.32794
3.034 	-42.28439
3.035 	-43.22468
3.036 	-44.14893
3.037 	-45.05726
3.038 	-45.94977
3.039 	-46.82659
3.040 	-47.68783
3.041 	-48.53359
3.042 	-49.36399
3.043 	-50.17915
3.044 	-50.97917
3.045 	-51.76415
3.046 	-52.53423
3.047 	-53.28949
3.048 	-54.03005
3.049 	-54.75602
3.050 	-55.46751
3.051 	-56.16462
3.052 	-56.84746
3.053 	-57.51613
3.054 	-58.17075
3.055 	-58.81141
3.056 	-59.43822
3.057 	-60.05128
3.058 	-60.65071
3.059 	-61.23659
3.060 	-61.80904
3.061 	-62.36815
3.062 	-62.91403
3.063 	-63.44678
3.064 	-63.96650
3.065 	-64.47328
3.066 	-64.96724
3.067 	-65.44846
3.068 	-65.91704
3.069 	-66.37309
3.070 	-66.81670
3.071 	-67.24797
3.072 	-67.66699
3.073 	-68.07386
3.074 	-68.46868
3.075 	-68.85154
3.076 	-69.22253
3.077 	-69.58176
3.078 	-69.92931
3.079 	-70.26527
3.080 	-70.58975
3.081 	-70.90283
3.082 	-71.20460
3.083 	-71.49515
3.084 	-71.77459
3.085 	-72.04298
3.086 	-72.30044
3.087 	-72.54704
3.088 	-72.78288
3.089 	-73.00804
3.090 	-73.22261
3.091 	-73.42668
3.092 	-73.62033
3.093 	-73.80366
3.094 	-73.97675
3.095 	-74.13969
3.096 	-74.29255
3.097 	-74.43543
3.098 	-74.56841
3.099 	-74.69158
3.100 	-74.80501
3.101 	-74.90880
3.102 	-75.00302
3.103 	-75.08776
3.104 	-75.16309
3.105 	-75.22911
3.106 	-75.28588
3.107 	-75.33350
3.108 	-75.37204
3.109 	-75.40158
3.110 	-75.42220
3.111 	-75.43398
3.112 	-75.43700
3.113 	-75.43134
3.114 	-75.41707
3.115 	-75.39427
3.116 	-75.36302
3.117 	-75.32340
3.118 	-75.27548
3.119 	-75.21933
3.120 	-75.15504
3.121 	-75.08267
3.122 	-75.00230
3.123 	-74.91401
3.124 	-74.81787
3.125 	-74.71395
3.126 	-74.60232
3.127 	-74.48306
3.128 	-74.35624
3.129 	-74.22193
3.130 	-74.08020
3.131 	-73.93112
3.132 	-73.77477
3.133 	-73.61121
3.134 	-73.44052
3.135 	-73.26275
3.136 	-73.07799
3.137 	-72.88630
3.138 	-72.68775
3.139 	-72.48240
3.140 	-72.27032
3.141 	-72.05158
3.142 	-71.82625
3.143 	-71.59440
3.144 	-71.35608
3.145 	-71.11137
3.146 	-70.86032
3.147 	-70.60301
3.148 	-70.33950
3.149 	-70.06986
3.150 	-69.79413
3.151 	-69.51240
3.152 	-69.22472
3.153 	-68.93116
3.154 	-68.63177
3.155 	-68.32663
3.156 	-68.01578
3.157 	-67.69930
3.158 	-67.37724
3.159 	-67.04966
3.160 	-66.71663
3.161 	-66.37820
3.162 	-66.03443
3.163 	-65.68538
3.164 	-65.33112
3.165 	-64.97170
3.166 	-64.60717
3.167 	-64.23759
3.168 	-63.86303
3.169 	-63.48354
3.170 	-63.09918
3.171 	-62.70999
3.172 	-62.31605
3.173 	-61.91740
3.174 	-61.51410
3.175 	-61.10620
3.176 	-60.69377
3.177 	-60.27685
3.178 	-59.85549
3.179 	-59.42976
3.180 	-58.99970
3.181 	-58.56537
3.182 	-58.12682
3.183 	-57.68410
3.184 	-57.23727
3.185 	-56.78638
3.186 	-56.33148
3.187 	-55.87261
3.188 	-55.40984
3.189 	-54.94321
3.190 	-54.47277
3.191 	-53.99857
3.192 	-53.52067
3.193 	-53.03910
3.194 	-52.55393
3.195 	-52.06519
3.196 	-51.57294
3.197 	-51.07723
3.198 	-50.57809
3.199 	-50.07559
3.200 	-49.56977
3.201 	-49.06067
3.202 	-48.54834
3.203 	-48.03284
3.204 	-47.51419
3.205 	-46.99246
3.206 	-46.46767
3.207 	-45.93989
3.208 	-45.40916
3.209 	-44.87551
3.210 	-44.33900
3.211 	-43.79967
3.212 	-43.25756
3.213 	-42.71272
3.214 	-42.16518
3.215 	-41.61500
3.216 	-41.06221
3.217 	-40.50685
3.218 	-39.94898
3.219 	-39.38862
3.220 	-38.82583
3.221 	-38.26064
3.222 	-37.69310
3.223 	-37.12324
3.224 	-36.55111
3.225 	-35.97674
3.226 	-35.40018
3.227 	-34.82147
3.228 	-34.24064
3.229 	-33.65773
3.230 	-33.07279
3.231 	-32.48585
3.232 	-31.89695
3.233 	-31.30613
3.234 	-30.71342
3.235 	-30.11887
3.236 	-29.52251
3.237 	-28.92437
3.238 	-28.32450
3.239 	-27.72294
3.240 	-27.11971
3.241 	-26.51485
3.242 	-25.90841
3.243 	-25.30041
3.244 	-24.69089
3.245 	-24.07988
3.246 	-23.46743
3.247 	-22.85356
3.248 	-22.23832
3.249 	-21.62172
3.250 	-21.00382
3.251 	-20.38463
3.252 	-19.76421
3.253 	-19.14256
3.254 	-18.51974
3.255 	-17.89578
3.256 	-17.27070
3.257 	-16.64454
3.258 	-16.01733
3.259 	-15.38910
3.260 	-14.75988
3.261 	-14.12971
3.262 	-13.49862
3.263 	-12.86664
3.264 	-12.23379
3.265 	-11.60011
3.266 	-10.96563
3.267 	-10.33038
3.268 	-9.69438
3.269 	-9.05768
3.270 	-8.42029
3.271 	-7.78225
3.272 	-7.14358
3.273 	-6.50432
3.274 	-5.86449
3.275 	-5.22412
3.276 	-4.58324
3.277 	-3.94187
3.278 	-3.30005
3.279 	-2.65781
3.280 	-2.01516
3.281 	-1.37213
3.282 	-0.72876
3.283 	-0.08507
3.284 	0.55891
3.285 	1.20317
3.286 	1.84766
3.287 	2.49238
3.288 	3.13728
3.289 	3.78235
3.299 	10.23668
3.309 	16.68108
3.319 	23.09331
3.329 	29.45292
3.339 	35.74118
3.349 	41.94097
3.359 	48.03672
3.369 	54.01432
3.379 	59.86102
3.389 	65.56542
3.399 	71.11733
3.409 	76.50775
3.419 	81.72876
/

\setsolid

\axis bottom label {} shiftedto y=0.0  ticks numbered from 2.9 to 3.5 by 0.1     unlabeled short from   2.95 to 3.45 by 0.10 /
\axis left label   {} shiftedto x=2.9  ticks numbered from -80.0 to 80.0 by 20   unlabeled short from   -70 to 70 by 20 /

\put {$U_1(y)\times 10^6$ (P{\'o}lya) }   [l] at  2.91 80.0
\put {$y$}   [r] at  3.5 7.0
\put {$U''_1(y) > 0$} [c] at 3.23 20.0
\arrow <1,5mm>   [0.25,0.75] from  3.23 15.0 to 3.28 2

\endpicture     \]
\vspace{-5mm}

\caption{A section of the curve \(U_1(y)\) for P\'olya's function, with \(\alpha=0.99\).}

\end{figure}
